\documentclass[11pt]{article}
\usepackage[centertags]{amsmath}
\usepackage{amsmath}
\usepackage{amssymb}
\usepackage{amsthm}
\usepackage{latexsym}
\usepackage[T1]{fontenc}
\usepackage[latin1]{inputenc}
\usepackage{fancyhdr}
\DeclareMathOperator*{\esssup}{ess\,sup}
\DeclareMathOperator*{\essinf}{ess\,inf}
\numberwithin{equation}{section}

\newcommand{\U}{{\cal U}}
\newcommand{\V}{{\cal V}}

\newcommand{\be}{\begin{equation}}
\newcommand{\ee}{\end{equation}}

\newtheorem{definition}{Definition}[section]
\newtheorem{theorem}{Theorem}[section]
\newtheorem{proposition}{Proposition}[section]
\newtheorem{lemma}{Lemma}[section]
\newtheorem{remark}{Remark}[section]

\newcommand{\eps}{\varepsilon}

\def \esssup {\mbox{ess sup}}

\def \cp {{\cal P}}

\def \R{\mathbb{R}}
\def \P{\mathbb{P}}
\def \Q{\mathbb{Q}}

\def \E{\mathbb{E}}

\def \bf{\textbf}
\def \it{\textit}

\def \ni {\noindent}

\def \F{\mathcal{F}}

\def \mx{\mbox}

\usepackage{amssymb}
\usepackage{amsthm}
\usepackage{amsmath}
\usepackage{latexsym}
\usepackage[T1]{fontenc}
\usepackage[latin1]{inputenc}

\def \R{\mathbb{R}}
\def \P{\mathbb{P}}

\def \E{\mathbb{E}}

\def \bf{\textbf}
\def \it{\textit}

\def \ni {\noindent}

\def \F{\mathcal{F}}

\numberwithin{equation}{section}

\def \R{\mathbb{R}}
\def \P{\mathbb{P}}

\def \E{\mathbb{E}}

\def \bf{\textbf}
\def \it{\textit}

\def \ni {\noindent}

\def \F{\mathcal{F}}

\def \o {\omega}
\def \O {\Omega}
\def \ms {\medskip}

\def \ed {\end{document}}
\def \cc {{\cal C}}

\def \fr {\forall}
\def \t {\tau}

\def \tx {{t,x}}

\def \rw {\rightarrow}
\def \ttrw {(t',x')\rw (t,x)}
\def \esssup {\mbox{ess sup}}

\def \cp {{\cal P}}

\def \R{\mathbb{R}}
\def \P{\mathbb{P}}

\def \E{\mathbb{E}}

\def \bf{\textbf}
\def \it{\textit}

\def \ni {\noindent}

\def \F{\mathcal{F}}

\def \mx{\mbox}

\usepackage{amssymb}
\usepackage{amsthm}
\usepackage{amsmath}
\usepackage{latexsym}
\usepackage[T1]{fontenc}
\usepackage[latin1]{inputenc}

\def \R{\mathbb{R}}
\def \P{\mathbb{P}}

\def \E{\mathbb{E}}

\def \bf{\textbf}
\def \it{\textit}

\def \ni {\noindent}

\def \F{\mathcal{F}}

\numberwithin{equation}{section}

\def \R{\mathbb{R}}
\def \P{\mathbb{P}}

\def \E{\mathbb{E}}

\def \bf{\textbf}
\def \it{\textit}

\def \ni {\noindent}

\def \F{\mathcal{F}}

\def \o {\omega}
\def \O {\Omega}
\def \ms {\medskip}

\def \ed {\end{document}}
\def \cc {{\cal C}}

\def \fr {\forall}
\def \t {\tau}

\def \tx {{t,x}}

\def \rw {\rightarrow}
\def \ttrw {(t',x')\rw (t,x)}
\begin{document}
% % % % % % % % % % % % % % % % % % % % % % % % % % % % % % % % % % % % % % % % % % % % % % % % % % % % % % % % % % % % % % % % % % % % % % % % % % % % % % % % %

\title{Mixed Zero-Sum Stochastic Differential Game and Doubly Reflected BSDEs with  a Specific Generator}

% % % % % % % % % % % % % % % % % % % % % % % % % % % % % % % % % % % % % % % % % % % % % % % % % % % % % % % % % % % % % % % % % % % % % % % % % % % % % % % % %

\author{ Brahim EL ASRI \thanks{Universit\'e Ibn Zohr, Equipe. Aide \'a la decision,
ENSA, B.P.  1136, Agadir, Maroc. e-mail: b.elasri@uiz.ac.ma } \,\,\, and \, Nacer OURKIYA \thanks{Universit\'e Ibn Zohr, Equipe. Aide \'a la decision,
ENSA, B.P.  1136, Agadir, Maroc. e-mail: nacer.ourkiya@edu.uiz.ac.ma.}}
\date{}
\maketitle

\begin{abstract}
This paper studies the mixed zero-sum stochastic differential game problem. We allow the functionals and dynamics to be  of polynomial
growth. The problem is formulated as an extended  doubly reflected BSDEs with  a specific generator. We  show the existence  of solution for this doubly reflected BSDEs and we  prove the existence of a saddle-point of the game.  Moreover,  in the Markovian framework we prove that the value function is the unique viscosity solution of the associated Hamilton-Jacobi-Bellman equation.
\end{abstract}

\noindent {${Keywords}:$} Reflected Backward stochastic differential
equations; Mixed stochastic control; control-stopping problem;  Hamilton-Jacobi-Bellman equation; Viscosity solution.
\ms

\medskip
% % % % % % % % % % % % % % % % % % % % % % % % % % % % % % % % % % % % % % % % % % % % % % % % % % % % % % % % % % % % % % % % % % % % % % % % % % % % % % % % %
\section{Introduction}
In this paper, we study a  mixed zero-sum stochastic differential game in finite horizon which can be formulated as follows.
We consider  the dynamics of a stochastic system is described by a functional differential equation of the form
\begin{equation*}
dx_t=f(t,x,u_t,v_t)dt+\sigma(t,x)dB^{(u,v)}_t;\qquad t\in[0,T]\, \mbox{and}\, x_0\in \mathbb{R}^d \,\mbox{fixed};
\end{equation*}
where $B^{(u,v)}_t$ is a Brownian motion under a probability $\P^{(u,v)}$.\\
$(u_t)_{t\leq T}$ and $(v_t)_{t\leq T}$ are  the stochastic processes and stand for, respectively, the intervention functions of two agents $C_1$ and $C_2$ on that system. The agent $C_1$  chooses the stopping time $\tau$ and the agent $C_2$ chooses the stopping time $\sigma$. The game stops when one player decides to stop, that is, at the stopping time $\tau\wedge\sigma$, or $\tau=\sigma=T$. The actions of the
agents are not free and their advantages are antagonistic, which means that between
$C_1$ and $C_2$ there is a payoff $J(u, \tau; v, \sigma)$ which is a cost for $C_1$ and a reward for $C_2$, $J(u, \tau; v, \sigma)$ given
by:
$$J(u,\tau;v,\sigma)=\mathbb{E}^{(u,v)}\biggl[\int^{\tau \wedge\sigma}_{0} \Gamma(s,x,u_s,v_s)ds+U_{\tau}1_{[\tau<\sigma]} +L_{\sigma}1_{[\sigma\leq\tau<T]}+g(x)1_{[\tau \wedge\sigma=T]}\biggr],$$
where $\Gamma$ is the instantaneous reward of the controller, $L$ (resp. $U$) is the reward if the controller decided
to stop at $\sigma$ (resp. $\tau$) before the terminal time $T$ and $g$ is the reward if he decides to stay until $T$. The payoff depends on the process $(x_t )_{t\leq T}$ and is the  accumulate of the instantaneous payoff $\Gamma$ over time  interval $[0,\tau \wedge\sigma]$, the reward of the controller  and terminal payoffs $g$. Therefore the agent $C_1$ aims at minimizing $J(u,\tau;v,\sigma)$ while $C_2$ aims at maximizing the same payoff. In the particular case of agents who have non control actions, the mixed game is just the well known Dynkin game (see e.g. $\cite{CK}$).
\par The problem we are interested in is finding a saddle-point for the payoff functional $J(u,\tau;v,\sigma)$, i.e., a quadruple  $(u^*,\tau^*;v^*,\sigma^*)$ which satisfies that for any $(u,\tau)$ and $(v,\sigma)$, we have
\begin{equation*}
 J(u^*,\tau^*;v,\sigma)\leq  J(u^*,\tau^*;v^*,\sigma^*)\leq  J(u^*,\tau^*;v,\sigma).
\end{equation*}
\par In comparison with the classical game problem, the one that we consider here is of mixed type since the controllers have two actions: continuous control and stopping.
\par It is well-known that the nonanticipative  mixed zero-sum stochastic differential game is related to doubly  reflected
backward stochastic differential equations (DRBSDEs for short) (see e.g. \cite {BHM, CK, H2006,HM, HH2006, HL2000, HRZ} ) of the following type
\begin{equation}\label{eq1}
\begin{cases}
\begin{aligned}
Y^{*}_t = g(x) &+ \int_{t}^{T} H^*(s,x,Z^{*}) \, \mathrm{d}{s} + (K^{*+}_T - K^{*+}_t)\\ & - ( K^{*-}_T - K^{*-}_t) -  \int_{t}^{T}Z^{*}_s \, \mathrm{d}{B_s};
\end{aligned}\\
\forall t \leq T, \quad L_t \leq Y^{*}_t \leq U_t,\\
\int_{0}^{T} (Y^{*}_s - L_s) \, \mathrm{d}{K^{*+}_s} = 0 \; and \; \int_{0}^{T} (U_s-Y^{*}_s) \, \mathrm{d}{K^{*-}_s} = 0,
\end{cases}
\end{equation}
where $$H^*(t,x,z)=\displaystyle\inf_{u\in A}\displaystyle\sup_{v\in B}H(t,x,z,u,v) = \displaystyle\sup_{v\in B}\displaystyle\inf_{u\in A}H(t,x,z,u,v),\,\, \mbox{with}$$
$$H(t,x,z,u,v):=z\sigma^{-1}(t,x)f(t,x,u,v)+\Gamma(t,x,u,v).$$
The function $H^*$ is the Hamiltonian of the problem, $A$ and $B$  the sets of values of the controls. $H^*$ associated with the game, and $(Y^{*},Z^{*}, K^{*+} ,K^{*-} )$ is the solution of the doubly reflected BSDE  in $U$ (upper) and $L$ (lower) and whose generator and terminal value are respectively $H^*$ and $g(x)$.

\par In the case of Markovian setting, the  mixed zero-sum stochastic differential game problem is also related to the
following double obstacle variational inequality (DOVI in short):
\begin{equation}\label{dovi0}
\begin{cases}
\begin{aligned}
&\min\biggl\{u(t,x)-h(t,x),\max
\biggl[-\frac{\partial u}{\partial t}(t,x)-\mathcal{L}u(t,x)\\ &\;\; \qquad-H^*(t,x,\sigma(t,x)\nabla u(t,x)),u(t,x)-h'(t,x)\biggr]\biggr\}=0 ,\;(t,x)\in[0,T[ \times \mathbb{R} ^d;\\
&u(T,x)=g(x),\qquad x\in\mathbb{R}^d,
\end{aligned}
\end{cases}
\end{equation}
where $h$, $h'$ are continuous and of polynomial growth, $\mathcal{L}$ an operator associated with a diffusion process.

\par In Hamadène and Lepeltier \cite{HL2000} and Hamadène \cite{H2006}, the authors have studied this problem when
 $\sigma^{-1}f$ is bounded then obviously the function $H^*$ is Lipschitz w.r.t.
$z$ and consequently the existence and uniqueness of a solution for the above  DRBSDEs  (\ref{eq1}) is
obtained from the classical results of \cite{CK, lp, HM, contcoef, contcoef04}. The PDIs approach turns out to study and to solve, in some sense,  double obstacle variational inequality (\ref{dovi0}). Amongst the papers which consider the same problem as ours, and in the framework of viscosity solutions approach, the most elaborated works are \cite{lp, HM}. The authors show existence and uniqueness of a solution for (\ref{dovi0}) when $H^*$ is Lipschitz w.r.t.
$z$.

\par Therefore, the novelty of this paper lies in the fact that we investigate the solution to the mixed zero-sum stochastic differential game problem, using doubly reflected BSDE  and  viscosity solutions approach when $\sigma^{-1}f$ is not bounded
and satisfies a linear growth condition.\\

\par This paper is organized as follows: in Section $2$, we give some notations and setting of the mixed zero-sum stochastic differential game problem.  In section $3$, we consider doubly reflected BSDEs with  a specific generator whose solution provides appropriate
estimates for the solution of the double barrier BSDE (\ref{eq1}). In section $4$, we show that equation
(\ref{eq1}) has a unique solution $(Y^*, Z^*, K^{*+},K^{-*})$. Further we show existence  of the value and of a saddle-point for the mixed zero-sum stochastic differential game problem. Finally in section $5$, the problem is formulated within a Markovian framework, we give connection with doubly reflected BSDEs (\ref{eq1})  with double obstacle variational inequality (\ref{dovi0}). We show existence and uniqueness of (\ref{dovi0}) in  viscosity sense.
\section{Notations and Setting of the mixed zero-sum stochastic differential game problem}

Let $\Omega=\cc\left([0,T];\R^d\right)$ be the space of $\R^d$-valued continuous function on $[0,T]$ endowed with the metric of uniform convergence on $[0,T]$.
Denote by $\F$ the Borel $\sigma$-field over $\Omega$. Next for $\o\in \O$ and $t\le T$, let us set $||\omega||_t:=\sup_{0\leq s\leq t}|\omega_s|$. Let $x:=(x_s)_{s\le T}$ be the coordinate process on $\O$, i.e., $x_s(\o)=\o(s)$ and denote by $(\F^0_t:=\sigma(x_s,s\leq t))_{t\in [0,T]}$, the filtration on $\O$ generated by $x$. Moreover

\begin{itemize}
\item[-] $\mathcal{P}$ be the progressively measurable $\sigma-$algebra on $[0,T]\times\Omega.$
\item[-] For any stopping time $\tau\in[0,T],$ $\mathcal{T}_{\tau}$ denotes the set of all stopping times $\theta$ such that $\tau\leq\theta\leq T.$
\item[-] By convention $\inf\{\emptyset\}=+\infty.$
\item[-]$\mathcal{L}^{q}$  is the set of $\mathcal{F}_{T}$-measurable random variable $\xi:\Omega\longmapsto\mathbb{R}$ with $\mathbb{E}[|\xi|^{q}]<+\infty,$ $q\geq 1$.
\item[-]$\mathcal{H}^{q,d}$  is the set of $\mathcal{P}$-measurable processes $Z:[0,T]\times\Omega\longmapsto\mathbb{R}^d$ with $\mathbb{E}\left[\left(\int^{T}_{0}|Z_s|^2ds\right)^{\frac{q}{2}}\right]<+\infty.$
\item[-]$\mathcal{S}^{q}$  is the set of $\mathcal{P}$-measurable and continuous processes $Y:[0,T]\times\Omega\longmapsto\mathbb{R}$  with\\ $\mathbb{E}\left[\sup_{0 \leq t \leq T}|Y_{t}|^q\right]<\infty.$
\item[-]$\mathcal{S}_{ci}^q$  is the set of continuous, increasing and $\mathcal{P}$-measurable processes $A:[0,T]\times\Omega\longmapsto [0,\infty)$ with $A(0)=0$ and $\mathbb{E}[(A_{T})^{q}]<+\infty.$
\end{itemize}
\par We now introduce the mixed stochastic differential game problem.
\par Let $\sigma$ be a function from $[0,T]\times \Omega$ into $\mathbb{R}^{d\times d}$ such that:\\
$[\mathbf{A1}]:$  $\sigma$ is $\mathcal{P}$-measurable.\\
$[\mathbf{A2}]:$   There exists a constant $C > 0$ such that:
\begin{itemize}
\item[(i)] For every $t \in [0,T]$ and $x, x'\in\Omega,$ $|\sigma(t,x)- \sigma(t,x')|\leq C||x- x'||_t,$ with $C$ independent of $t,$ $x$ and $x'$.
\item[(ii)]  $\sigma$ is bounded, invertible and its inverse $\sigma^{-1}$ is bounded.
\end{itemize}
Let $\P$ be a probability measure on $\Omega$ such that $(\Omega, \P)$ carries a $d$-dimensional Brownian motion $(B_t)_{0\leq t \leq T}$ and for any  $x_0 \in \R^d$, the process $(x_t)_{0\leq t \leq T}$ is the unique solution of the
following stochastic differential equation:
 $$dx_t=\sigma(t,x)dB_t, \;\;x_0\in\mathbb{R}^d.$$
Moreover, for any $p \geq 1,$
\begin{equation}\label{estmx}
\mathbb{E}[||x||^{p}_{T}]\leq C_p,
\end{equation}
where $C_p$ depends only on $p,$ $T,$ the initial value $x_0$ and the linear growth constant of $\sigma$ $($see, Karatzas and Shreve $\cite{KS},$ $p.\;306)$.\\
Let $A$ (resp., $B$) be a compact metric space and $\U$ (resp., $\V$) be the space of $\mathcal{P}$-measurable processes $u := (u_t )_{t\leq T}$ (resp., $v:= (v_t )_{t\leq T}$ ) with values in $A$ (resp., $B$). Let $f$ and $\Gamma$ such that
\begin{itemize}
\item[i)] $f:[0,T]\times\Omega\times A \times B \longmapsto\mathbb{R}^{d}$ is $\mathcal{P}\otimes\mathcal{B}(A\times B)$-measurable, for all $(t,x)\in[0,T]\times\Omega$, $f(t,x,.,.)$ is continuous on $A\times B$, and there exists a real constant $C$ such that:
\begin{equation}
\mathbf{[H1]:}\,\,|f(t,x,u,v)|\leq C(1+||x||_{t}), \qquad \forall (t,x,u,v)\in[0,T]\times\Omega\times A \times B .
\end{equation}
\item[ii)] $\Gamma:[0,T]\times\Omega\times A \times B \longmapsto\mathbb{R}$ is $\mathcal{P}\otimes\mathcal{B}(A\times B)$-measurable,  for all $(t,x)\in[0,T]\times\Omega$, $\Gamma(t,x,.,.)$ is continuous on $A\times B$. In addition there exists  a positive constants  $C>0$ and $p$ such that:
\begin{equation}
\mathbf{[H2]:}\,\,|\Gamma(t,x,u,v)|\leq C(1+||x||_{t}^{p}), \qquad \forall (t,x,u,v)\in[0,T]\times\Omega\times A \times B .
\end{equation}
\end{itemize}
For any $(u,v)\in \U \times \V$, we define a probability $\P^{(u,v)}$ on $(\Omega,\mathcal{F})$ by
\begin{equation*}
\frac{\mathrm{d}{\P^{(u,v)}}}{\mathrm{d}{\P}}:= \exp\left\{\int_{0}^{T}\sigma^{-1}(s,x)f(s,x,u_s,v_s)\, \mathrm{d}{B_s}-\frac{1}{2}\int_{0}^{T}||\sigma^{-1}(s,x)f(s,x,u_s,v_s)||^2 \, \mathrm{d}{s} \right\}.
\end{equation*}
Under $\P^{(u,v)}$ the process $(x_t )_{t\leq T}$ verifies: $\forall t\leq T$,
\begin{equation}\label{eds}
x_t = x_0+ \int_{0}^{t} f(s,x,u_s,v_s) \, \mathrm{d}{s}  +  \int_{0}^{t}\sigma(s,x) \, \mathrm{d}{B^{(u,v)}_s},
\end{equation}
where the process $(B^{(u,v)}_t := B_t-\int_{0}^{t}\sigma(s,x)^{-1}f(s,x,u_s,v_s) \, \mathrm{d}{s})_{t\leq T},$ is a Brownian motion under $\P^{(u,v)}$. It means that $(x_t)_ {t\leq T}$ is a weak solution for the functional stochastic differential equation $(\ref{eds})$.
\par The process $x=(x_t)_{t\leq T}$ stands for the dynamics of the system when it is not controlled. On that system, two agents $C_1$ and $C_2$ intervene. A control action for $C_1$ (resp., $C_2$ ) is a process $u=(u_t)_{t\leq T}$ (resp., $v=(v_t)_{t\leq T}$ ) which belongs to $\U$ (resp., $\V$).
\par Thereby $\U$ (resp. $\V$) is called the set of admissible controls for $C_1$ (resp. $C_2$ ). When $C_1$ and $C_2$ act with, respectively, $u$ and $v,$ the law of the dynamics of the system is the same as the one of $x$ under $\P^{(u,v)}$. It turns out that the interventions of the agents
give rise to a drift $f(s,x,u_t,v_t)$ in the dynamics.
\par As previously mentioned, an intervention strategy for an agent combines control and stopping. Therefore a strategy for $C_1$ (resp. $C_2$ ) is a pair $(u,\tau)$ (resp. $(v,\sigma)$), where $u$ (resp. $v$) is an admissible control and $\tau$ (resp. $\sigma$) a stopping time. It implies that the agents keep controlling until one of them decides to stop and the system is actually stopped at $\tau\wedge\sigma$. So assume that $C_1$ (resp. $C_2$ ) chooses $(u,\tau)$ (resp. $(v,\sigma)$) as a strategy of intervention which is not free and generates a payoff, which is a cost for $C_1$ and a reward for $C_2$. The expression of that payoff is given by
\begin{equation}
J(u,\tau;v,\sigma)=\mathbb{E}^{(u,v)}\biggl[\int^{\tau \wedge\sigma}_{0} \Gamma(s,x,u_s,v_s)ds+U_{\tau}1_{[\tau<\sigma]} +L_{\sigma}1_{[\sigma\leq\tau<T]}+g(x)1_{[\tau \wedge\sigma=T]}\biggr],
\end{equation}
where
\begin{enumerate}
\item $g: \cc \rightarrow \R$ is a Borel measurable function of polynomial growth, i.e., there exists a constant $C$ and $p$ such that:
\begin{equation}
\mathbf{[H3]:}\,\,|g(x)|\leq C(1+||x||^{p}_{T}), \qquad\qquad \forall x\in\Omega.
\end{equation}
\item  $U:=(U_t )_{t\leq T}$ and $L:=(L_t )_{t\leq T}$ are processes of $\mathcal{S}^{2}$ and satisfy $L_t < U_t,$ for all $t\leq T$ and $L_T\leq g(x)\leq U_T.$ Moreover, There exist a positive constant $C$ and $p>1$ such that $\forall t\in [0,T]$ and $x\in\Omega,$ we have,
\begin{equation}\label{LU}
\mathbf{[H4]:}\,\,
\begin{cases}
|L_t|\leq C(1+||x||^{p}_{t}),\\
\qquad\qquad\qquad\qquad\qquad\qquad \forall (t,x)\in [0,T]\times\Omega.\\
|U_t|\leq C(1+||x||^{p}_{t}).
\end{cases}
\end{equation}
\end{enumerate}

\par The action of $C_1$ (resp. $C_2$) is to minimize (resp. maximize) the payoff $J(u,\tau;v,\sigma)$ whose terms can be understood as:
\begin{itemize}
\item[(i)]  $\Gamma(t,x,u,v)$ is the instantaneous reward (resp. cost) for $C_1$ (resp. $C_2$).
\item[(ii)]  $U_{\tau}$ is the cost (resp. reward) for $C_1$ (resp. $C_2$) if $C_1$ decides to stop first the game.
\item[(iii)] $L_{\sigma}$ is the reward (resp. cost) for $C_2$ (resp. $C_1$) if $C_2$ decides to stop first the game.
\end{itemize}
The problem we are interested in is finding two intervention strategies $(u^* ,\tau^*)$ and $(v^* ,\sigma^* )$, respectively, for $C_1$ and $C_2$ such that for any $(u,\tau)$ and $(v,\sigma)$, we have
\begin{equation}
J(u^*,\tau^*;v,\sigma)\leq J(u^*,\tau^*;v^*,\sigma^*) \leq J(u,\tau;v^*,\sigma^*).
\end{equation}
The quadruple $(u^*,\tau^*;v^*,\sigma^*)$ is called a saddle-point for the  mixed zero-sum stochastic differential game.

\section{Doubly reflected BSDEs with  a specific generator}

Let a function $\varphi$  is defined as follows: For some $C> 0$ and $p>1$, for any $(t,x,z) \in [0,T]\times\cc\times\mathbb{R}^{d},$
\begin{equation}
\varphi(t,x,z)=C(1+||x||_t)|z|+C(1+||x||^{p}_{t}).
\end{equation}

Next let  $\Phi$ be a function from $[0,T]\times\cc\ \times\R^d$ into $\R$ which is $\cp\otimes{\bf B}(\R^d)$-measurable. First we define the notion of a solution of the double barrier reflected BSDE associated with the quadruple $(g, \Phi, L, U)$ which we consider throughout this paper.
\begin{definition}\label{defiLpsol}
	A solution for the double barrier reflected BSDE associated with the generator $\Phi,$ the terminal value $g(x),$ the upper (resp. lower) obstacle $U$ (resp. $L$) is a process $(Y,Z,K^+,K^-):=(Y_t,Z_t,K^{+}_{t},K^{-}_{t})_{t\leq T},$ $\mathcal{P}$-measurable, with values in $\mathbb{R}^{1+d+1+1}$ such that:
	\begin{enumerate}
		\item $Y$ is continuous, $K^+$ and $K^-$ are continuous non-decreasing ($K^+_0=K^-_0=0$), $\P$-a.s. $(Z_t(\omega))_{t\leq T}$ is $dt$-square integrable;
		\item $Y_t = g(x) + \int_{t}^{T} \Phi(s,x,Z_s) \, \mathrm{d}{s} + (K_T^+ - K_t^+) - ( K_T^- - K_t^-) -  \int_{t}^{T}Z_s \, \mathrm{d}{B_s}$, $0 \leq t \leq T$;
		\item $\forall t \leq T, \quad L_t \leq Y_t \leq U_t,
		\;\;\int_{0}^{T} (Y_s - L_s) \, \mathrm{d}{K_s^+} = 0 \;\; and \;\; \int_{0}^{T} (U_s-Y_s) \, \mathrm{d}{K_s^-} = 0$.
	\end{enumerate}
\end{definition}
\par We are now going to introduce a property of exponential martingales which plays an important role in this work. Let $\vartheta$ be a function from
$[0,T]\times \Omega$ into $\R^d$ which is $\cp$-measurable. For $t\le T$ we set:
$$\zeta_t=\exp\biggl\{\int_0^t\vartheta(s,\omega)\sigma^{-1}(s,\omega)dB_s-\frac{1}{2}\int_0^t|
	\vartheta(s,\omega)\sigma^{-1}(s,\omega)|^2ds\biggr\}.$$
We then have:
\begin{lemma}\label{lemma21}${}$ (see, \cite{EHO})\\i) Assume that $\vartheta$ is bounded. Then
	$(\zeta_t)_{t\le T}$ is a martingale such that for any $\ell \ge 1$,
	$$
	\E[(\zeta_T)^\ell]\le C_{\ell,\vartheta, T}
	$$
	where $C_{\ell,\vartheta,T}$ is a constant which depends only on $T$, $\ell$ and the constant of boundedness of $\vartheta \sigma^{-1}$.
	\ms
	
	\ni ii) Assume that $\vartheta$ is of linear growth, i.e.,
	for any $t\le T$, $|\vartheta(t,\omega)|\le c_{\vartheta}(1+|\omega|_t)$. Then there exists a constant $q>1$ and a constant $\varpi$ which depend only on $T$, $\sigma$, $\sigma^{-1}$ and $c_{\vartheta}$ such that
	\begin{equation}\label{probhaus}
	\E[(\zeta_T)^q]\le \varpi.
	\end{equation}
	Moreover $(\zeta_t)_{t\le T}$ is a martingale.
	\ms
	
	\ni iii) In both cases, the measure $\Q$ such that $d\Q=\zeta_T.d\P$ is a probability on $(\Omega,\F)$.
\end{lemma}

We are going to show that the double barrier reflected BSDE associated with $(g, \varphi, L, U)$  has a solution which also satisfies other integrability properties. As it is mentionned previously, $\varphi$ is not a standard generator which does not enter neither in the framework of \cite{CK,lp,H2006,HM, HRZ} nor in the recent paper \cite{marzoug}.
\par We have the following result.
\begin{proposition}\label{prop2.1}
There exist $\mathcal{P}$-measurable processes $(Y,Z,K^+,K^-)$ valued in $\mathbb{R}^{1+d+1+1}$ and a stationary non-decreasing sequence of stopping times $(\tau_k)_{k\geq1}$ such that:
\begin{itemize}
\item[$(i)$] $(Y,Z,K^+,K^-)$ is a solution for the double barrier reflected BSDE associated with $(g, \varphi, L, U)$.
\item[$(ii)$] For any constant $\gamma\geq 1$ and $\tau$ a stopping time valued in $[0,T],$
\begin{equation}\label{Ygamma}
\mathbb{E}[|Y_{\tau}|^{\gamma}]<+\infty.
\end{equation}
\item[$(iii)$] For any $k\geq1$,
\begin{equation}\label{cond3}
\mathbb{E}\left[\sup_{t\leq T}|Y_{t\wedge\tau_k}|^{\gamma}+(K_{\tau_k}^{+})^{\gamma}+(K_{\tau_k}^{-})^{\gamma}+\int_{0}^{\tau_k}|Z_s|^2\mathrm{d}{s}\right]<+\infty.
\end{equation}
\end{itemize}
where $\tau_k$ depends only on $g,$ $L,$ $U$ and $x.$
\end{proposition}
\begin{proof}
See Appendix A.
\end{proof}

\section{ Resolution of the mixed zero-sum stochastic differential game problem}
\qquad Now we use the results of the above section to deal with an application of the double barrier reflected BSDE tool for solving our stochastic mixed game problem. We first start by introducing the objects  which we need in this section.

\par For any
$(t,x,z,u,v)\in [0,T]\times \Omega \times \mathbb{R}^d \times A \times B,$  we define the Hamiltonian associated with this mixed stochastic game problem by $H(t,x,z,u,v)=z\sigma^{-1}(t,x)f(t,x,u,v)+\Gamma(t,x,u,v)$ and we suppose the following assumption, which is called the Isaacs's condition (see, Hamadène and Lepeltier $\cite{HL2000}$ and Hamadène $\cite{H2006}$), holds:
\begin{equation}\label{isac}
\mathbf{[H5]:}\,\, \displaystyle\inf_{u\in A}\displaystyle\sup_{v\in B}H(t,x,z,u,v) = \displaystyle\sup_{v\in B}\displaystyle\inf_{u\in A}H(t,x,z,u,v),\quad \mbox{for any}\,\,
(t,x,z)\in [0,T]\times \cc \times \mathbb{R}^d.
\end{equation}
Under $(\ref{isac})$, through the Benes's theorem (Benes $\cite{benes},$ $1970$), there exists a couple of $\mathcal{P}\otimes\mathcal{B}(\mathbb{R}^{d})$-measurable functions $u^*(t,x,z)$ and $v^*(t,x,z)$ with values respectively in $A$ and $B$ such that $\P$-a.s., for any $(t,z)\in[0,T]\times \mathbb{R}^{d},$ $u\in A$ and $v\in B$,
\begin{equation}\label{H*}
\begin{aligned}
H^*(t,x,z):=H(t,x,z,u^*(t,x,z),v^*(t,x,z)) & = \displaystyle\inf_{u\in A}\displaystyle\sup_{v\in B}H(t,x,z,u,v),\\
                               & = \displaystyle\sup_{v\in B}\displaystyle\inf_{u\in A}H(t,x,z,u,v),
\end{aligned}
\end{equation}
and
\begin{equation}\label{H*H}
H(t,x,z,u^*(t,x,z),v) \leq H(t,x,z,u^*(t,x,z),v^*(t,x,z)) \leq H(t,x,z,u,v^*(t,x,z)).
\end{equation}
\par First we are going to show that the double barrier reflected BSDE associated with  $(g, H^*,L, U)$ has a unique solution which moreover verifies some integrability properties. The following proposition is a step forward the proof of this result.
 \begin{proposition}\label{prop31}
For any $m\geq0$, there exist $\mathcal{P}$-measurable processes $(Y^{*m},Z^{*m},K^{*+,m},K^{*-,m})$ in $\mathbb{R}^{1+d+1+1}$ such that:
\begin{itemize}
\item[$(i)$] For any $\gamma\geq 1$ and any stopping time $\tau\in[0,T],$ such that,
\begin{equation}\label{ym}
\mathbb{E}[|Y^{*m}_{\tau}|^{\gamma}]\leq C, \qquad\forall m\geq1,
\end{equation}
where $C$ is a constant independent of $m$ and $\tau$.
\item[$(ii)$] The quadruple
$(Y^{*m},Z^{*m},K^{*+,m},K^{*-,m})$ is a solution of the double barrier reflected BSDE associated with $(g, H^{*m},L, U)$ where for any $t\leq T,$ $$H^{*m}(t,x,z)=H^{*+}(t,x,z)-H^{*-}(t,x,z)1_{\{||x||_t\leq m\}},$$ with ${H^*}^+$(resp. ${H^*}^-$) is the positive  (resp.  negative)  part of $H^*$.
\end{itemize}
\end{proposition}
\begin{proof}
See Appendix B.
\end{proof}

The following result is related to the caracterization of the value function of the mixed control problem via double barrier reflected BSDE. It also provides a saddle-point of this latter problem.
\begin{theorem}\label{theo31}
Assume that $\mathbf{[H1]}$ to $\mathbf{[H5]}$ are satisfied. Then, there exist $\mathcal{P}$-measurable processes $(Y^{*},Z^{*},K^{*+},K^{*-})$ valued in $\mathbb{R}^{1+d+1+1}$ such that:
\begin{itemize}
	\item[$i)$] For any $\gamma\geq 1$ and any stopping time $\tau\in[0,T]$ such that,
	\begin{equation}\label{majy}
	\mathbb{E}[|Y^{*}_{\tau}|^{\gamma}]<C,
	\end{equation}
	where $C$ is a constant independent from $\tau$.
	\item[$ii)$] $Y^{*}$  and $K^{*\pm}$ are continuous, $K^{*\pm}$ is non-decreasing $(K^{*\pm}_0=0)$ and $\P$-a.s. $(Z^{*}_t)_{t\leq T}$ is $dt$-square integrable.
\item[$iii)$] For any $t\in[0,T],$
\begin{equation}\label{hRBSDE2*}
\begin{cases}
\begin{aligned}
Y^{*}_s = g(x) &+ \int_{t}^{T} H^{*}(s,x,Z^{*}_s) \, \mathrm{d}{s} + (K^{*+}_T - K^{*+}_t) - ( K^{*-}_T - K^{*-}_t) \\ &-  \int_{t}^{T}Z^{*}_s \, \mathrm{d}{B_s};
\end{aligned}\\
\forall t \leq T, \quad L_t \leq Y^{*}_s \leq U_t,\\
\int_{0}^{T} (Y^{*}_s - L_s) \, \mathrm{d}{K^{*+}_s} = 0 \; and \; \int_{0}^{T} (U_s-Y^{*}_s) \, \mathrm{d}{K^{*-}_s} = 0.
\end{cases}
\end{equation}
\item[$iv)$] If $(\underline{Y}, \underline{Z}, \underline{K}^{+}, \underline{K}^{-})$  is another quadruple which satisfies i), ii) and iii), then
$(\underline{Y}, \underline{Z}, \underline{K}^{+}, \underline{K}^{-})=(Y^{*},Z^{*},K^{*+},K^{*-})$, i.e., the solution
of the double barrier reflected BSDE associated with
$(g,H^*, L, U)$ is unique to satisfy $i)-iii).$
\end{itemize}
\end{theorem}
\begin{proof}
From Proposition $\ref{theo31},$ we have that for any $m\geq 1$ there exists a quadruple\\  $(Y^{*m},Z^{*m},K^{*+,m},K^{*-,m})$ that satisfies the following double barrier reflected BSDE:
\begin{equation*}
\begin{cases}
\begin{aligned}
Y^{*m}_s = g(x) &+ \int_{t}^{T} H^{*m}(s,x,Z^{*m}_s) \, \mathrm{d}{s} + (K^{*+,m}_T - K^{*+,m}_t)\\ & - ( K^{*-,m}_T - K^{*-,m}_t) -  \int_{t}^{T}Z^{*m}_s \, \mathrm{d}{B_s},
\end{aligned}\\
\forall t \leq T, \quad L_t \leq Y^{*m}_s \leq U_t,\\
\int_{0}^{T} (Y^{*m}_s - L_s) \, \mathrm{d}{K^{*+,m}_s} = 0 \; and \; \int_{0}^{T} (U_s-Y^{*m}_s) \, \mathrm{d}{K^{*-,m}_s} = 0.
\end{cases}
\end{equation*}
Where, for any $t\leq T,$ $H^{*m}(t,x,z)=H^{*+}(t,x,z)-H^{*-}(t,x,z)1_{\{||x||_t\leq m\}}.$
\par First let us prove $(i)$.\\
From $\eqref{ymym+1}$ we know that for any $m\ge 1$, $L\le Y^{*,m+1}\le Y^{*,m}\le Y$. So, let us set for $t\le T$,
$$Y^*_t=\lim_m Y^{*,m}_t.$$
Therefore $L\le Y^*\le Y$ and then by  $\eqref{estmx},$ $\eqref{LU}$ and $\eqref{Ygamma}$, $Y^*$ verifies $\eqref{majy}.$
\par Next let $(\tau^*_k)_{k\geq1}$ be the sequence of stopping times defined as follows:
$$ \forall k\geq1,~~~\tau^{*}_k:=\inf\{t\geq0, |Y_t|+||x||_t+|L_t|+|U_t|\geq \zeta^{*} +k \}\wedge T,$$
where $\zeta^{*}=||x_0||+|L_0|+|U_0|+|Y_0|$.\\
Thanks to the continuity of $Y$, $L$, $U$ and $x$ the sequence of stopping time $(\tau^*_k)_{k\geq1}$ is increasing of stationary type that converges to $T$.\\
Under $(\ref{ymym+1})$ and as $\forall\gamma\geq1,$ $\mathbb{E}[(K_{\tau^{*}_k}^{*+,1})^{\gamma}+(K_{\tau^{*}_k}^{-})^{\gamma}]<\infty,$ then $\P$-a.s., for any $t\leq T$ and $m\geq1,$ the sequence $(K_t^{*+,m})_{m\geq 1}$ (resp. $(K_t^{*-,m})_{m\geq 1}$) converges to $K_t^{*+}$ (resp., $K_t^{*-})$.
In addition, the process  $K^{*+}:=(K_t^{*+})_{t\leq T}$ (resp. $K^{*-}:=(K_t^{*-})_{t\leq T}$) is non-decreasing lower and (resp. upper) semi-continuous and $\mathbb{E}[(K_{\tau^{*}_k}^{*+})^{\gamma}]<\infty$ (resp. $\mathbb{E}[(K_{\tau^{*}_k}^{*-})^{\gamma}]<\infty$).
\par
In the same way as in the proof of the Proposition $\ref{prop31},$ there exists a constant $C_k$ that depends on $k$ such that:
\begin{equation}\label{zmtau*}
\mathbb{E}\left[\int^{\tau^*_k}_{0}|Z^{*m}_s|^2ds\right]<C_k,\qquad \forall m\geq1.
\end{equation}
This inequality follows in a classic way after using Itô's formula with $(Y^{*m}_{t\wedge\tau^{*}_k})^2$
and taking into account $\eqref{estmx}$. Next we apply again Itô's formula with $(Y^{*m}_{t\wedge\tau^{*}_k}-Y^{*p}_{t\wedge\tau^{*}_k})^2$ and using standard calculations we obtain: for any $t\leq T,$ and $\forall k\geq1$,
 \begin{equation}\label{cvguni}
 \mathbb{E}\left[\sup_{t\leq T}|Y^{*m}_{t\wedge\tau_k}-Y^{*p}_{t\wedge\tau_k}|^2 +\int_{t\wedge\tau_k}^{\tau_k}|Z^{*m}_s-Z^{*p}_s|^2 \, \mathrm{d}{s}\right]\longrightarrow 0, \;\; as \,\, m,\, p \longrightarrow +\infty.
 \end{equation}
 It follows that the process $(Y^*_{t\wedge\tau^*_k})_{t\leq T}$ is continuous for any $k\geq1.$ As $(\tau^*_k)_{k\geq 1}$ is a stationary sequence then the process $(Y^*_{t})_{t\leq T}$ is continuous. On the other hand the sequence $((Z^{*m}_t1_{\{t\leq \tau^*_k\}})_{t\leq T})_{m\geq 1}$ converges in $\mathcal{H}^{2,d}$ to a process which we denote $(Z_t^{*k})_{t\leq T} .$ In addition it satisfies, for any $k\geq1$, $Z_{t\wedge\tau^{*}_k}^{*k+1}=Z_{t\wedge\tau^{*}_k}^{*k},$ $dt\otimes d\P-a.s.$.\\
 Now by $(\ref{zmtau*})$ and $(\ref{cvguni})$, for any $k \geq1$, the sequence of process
  $((H^{*m}(t,x,Z^{*m}_t)1_{\{t\leq \tau^*_k\}})_{t\leq T})_{m\geq 1}$ converges in $\mathcal{H}^{2,d}$ to $(H^{*}(t,x,Z^{*k}_t)1_{\{t\leq \tau^*_k\}})_{t\leq T}$. Therefore for any $k\geq1$ we have,
 \begin{equation}\label{eqzk*}
 \begin{aligned}
 Y^{*}_{t\wedge\tau^{*}_k} = & Y^{*}_{\tau^{*}_k} + \int_{t\wedge\tau^{*}_k}^{\tau^{*}_k} H^{*}(s,x,Z^{*k,}_s) \, \mathrm{d}{s} + (K_{\tau^{*}_k}^{*+} - K_{t\wedge\tau^{*}_k}^{*+}) \\
                  &  - ( K_{\tau^{*}_k}^{*-} - K_{t\wedge\tau^{*}_k}^{*-}) -  \int_{t\wedge\tau^{*}_k}^{\tau^{*}_k}Z^{*k}_s \, \mathrm{d}{B_s}, \qquad \forall t\leq T.
 \end{aligned}
 \end{equation}
 Now $(\ref{eqzk*})$ implies also that,
 \begin{equation*}
 Y^{*}_{t\wedge\tau^{*}_k} = Y^{*}_{0} -\int_{0}^{t\wedge\tau^{*}_k} H^{*}(s,x,Z^{*k}_s) \, \mathrm{d}{s}- K_{t\wedge\tau^{*}_k}^{*+}+K_{t\wedge\tau^{*}_k}^{*-}+  \int_{0}^{t\wedge\tau^{*}_k}Z^{*k,m}_s \, \mathrm{d}{B_s}, \qquad \forall t\leq T.
 \end{equation*}
 And then,
 \begin{equation*}
 K_{t\wedge\tau^{*}_k}^{*+}= K_{t\wedge\tau^{*}_k}^{*-}+Y^{*}_{0} -\int_{0}^{t\wedge\tau^{*}_k} H^{*}(s,x,Z^{*k}_s) \, \mathrm{d}{s}+  \int_{0}^{t\wedge\tau^{*}_k}Z^{*k}_s \, \mathrm{d}{B_s}, \qquad \forall t\leq T.
 \end{equation*}
 As $K^{*+}$ is lower semi-continuous and $K^{*-}$ is upper semi-continuous. It means that $K^{*+}$ and $K^{*-}$ are lower and upper semi-continuous in the same time then the processes $(K_{t\wedge\tau^{*}_k}^{*+})_{t\leq T}$ and $(K_{t\wedge\tau^{*}_k}^{*-})_{t\leq T}$ are continuous. Henceforth $K^{*+}$ and $K^{*-}$ are continuous since $(\tau^{*}_k)_{k\geq1}$ is of stationary type. In addition, from Dini's Theorem, the sequences $(K^{*+,m})_{m\geq1}$ and $(K^{*-,m})_{m\geq1}$ converge $\P$-a.s. uniformly to $K^{*+}$ and $K^{*-}$ respectively.
 \par Now for any $t\leq T$, let us set
 \begin{equation}
 Z^{*}_t= Z^{*1}_t1_{\{t\leq\tau^{*}_1\}}+\sum_{k\geq 2} Z^{*k}_t1_{\{\tau^{*}_{k-1}<t\leq\tau^{*}_{k}\}}.
 \end{equation}
 As $\mathbb{E}\left[\int^{\tau^{*}_k}_{0}|Z^{*k}_s|^2ds\right]<+\infty,$ for any $k\geq1$ and the sequence $(\tau^{*}_k)_{k\geq1}$ is of stationary type then $\int^{T}_{0}|Z^{*}_s|^2ds<+\infty,$ $\P$-a.s.. On the other hand, with the definition of $Z^*$ and $(\ref{eqzk*})$ we have,
 \begin{equation}\label{equtau*}
 \begin{aligned}
 Y^{*}_{t\wedge\tau^{*}_k} = & Y^{*}_{\tau^{*}_k} + \int_{t\wedge\tau^{*}_k}^{\tau^{*}_k} H^{*}(s,x,Z^{*}_s) \, \mathrm{d}{s} + (K_{\tau^{*}_k}^{*+} - K_{t\wedge\tau^{*}_k}^{*+}) \\
                  &  - ( K_{\tau^{*}_k}^{*-} - K_{t\wedge\tau^{*}_k}^{*-}) -  \int_{t\wedge\tau^{*}_k}^{\tau^{*}_k}Z^{*}_s \, \mathrm{d}{B_s}, \qquad \forall t\leq T.
 \end{aligned}
 \end{equation}
 Now taking $k$ great enough in $(\ref{equtau*})$ yields, $\forall t\leq T $
 \begin{equation*}
 Y^{*}_t = g(x) + \int_{t}^{T} H^{*}(s,x,Z^{*}_s) \, \mathrm{d}{s} + (K^{*+}_T- K^{*+}_t) - (K^{*-}_T- K^{*-}_t) -  \int_{t}^{T}Z^{*}_s \, \mathrm{d}{B_s}.
 \end{equation*}
 \par Finally, from the uniform convergence of $Y^{*m}$, $K^{*+,m}$ and  $K^{*-,m}$ we have:
 $$\int_{0}^{T} (Y^{*}_s - L_s) \, \mathrm{d}{K_s^{*+}} = \int_{0}^{T} (U_s-Y^{*}_s) \, \mathrm{d}{K_s^{*-}} = 0.$$
 Which completes the proof of $i)$, $ii)$ and $iii)$.
 \par
  We will now prove $iv)$. Let $(\underline{Y}, \underline{Z}, \underline{K}^{+}, \underline{K}^{-})$ be another quadruple which satisfies $i)$, $ii)$ and $iii)$. Then, using Itô's formula with $(Y^{*}_{t}-\underline{Y}_{t})^{2}$ yields, for any $t\leq T$,
 \begin{equation}\label{unis}
 \begin{aligned}
 (Y^{*}_{t}-\underline{Y}_{t})^{2}&=-2\int_{t}^{T}(Y^*_s-\underline{Y}_s)(Z^*_s-\underline{Z}_s)dB_s\\ & +2\int_{t}^{T}(Y^{*}_s-\underline{Y}_s)(H^{*}(s,x,Z_s^{*})-H^{*}(s,x,\underline{Z}_s))ds\\ &  +2\int_{t}^{T}(Y^{*}_s-\underline{Y}_s)d(K^{*+}_s-\underline{K}^+_s-K^{*-}_s+\underline{K}^-_s)-\int_{t}^{T}|Z^{*}_s-\underline{Z}_s|^2 ds.
 \end{aligned}
 \end{equation}
 Next let $\P^*$ be the probability, equivalent to $\P$, defined as follows:
\begin{equation*}
d\P^*:=\exp\biggl\{\int_0^T \Delta H^*(s)dB_s-\frac{1}{2}\int_0^T||\Delta H^*(s)||^2ds \biggr\}d\P,
\end{equation*}
where $$\Delta H^*(s):=\frac{H^{*}(s,x,Z_s^{*})-H^{*}(s,x,\underline{Z}_s)}{Z^*_s - \underline{Z}_s}1_{\{Z^*_s - \underline{Z}_s\neq 0\}}, $$ is a $\P$-measurable, $\mathbb{R}^{d}$-valued stochastic process such that
$$\forall s\leq T,\; H^{*}(s,x,Z_s^{*})-H^{*}(s,x,\underline{Z}_s) = \Delta H^*(s)(Z^*_s - \underline{Z}_s) .$$
As $$| H^{*}(s,x,z)-H^{*}(s,x,\underline{z})| \leq C(1+||x||_s) |z-\underline{z}|,$$
then $$\forall s\leq T,\;\|\Delta H^*(s)\|\leq C(1+||x||_s) .$$
It means that $\P^*$ is actually a probability equivalent to $\P$  (by Lemma \ref{lemma21}). Next for $k\ge 1$, let $\t_k$ be the following stopping time:
$$\tau_k:=\inf\{t\geq0, |Y^*_t|+|\underline{Y}_t|+\int_0^t|Z^*_s|ds+\int_0^t|\underline{Z}_s|ds \ge  k+|Y^*_0|+|\underline{Y}_0| \}\wedge T.$$
As $(Y^*_t-\underline{Y}_t)d(K^{*+}_t-\underline{K}^+_t-K^{*-}_t+\underline{K}^-_t)\leq 0,$ then going back to $\eqref{unis}$ to obtain:
$$(Y^{*}_{t\wedge \tau_k}-\underline{Y}_{t \wedge \tau_k})^{2}\leq (Y^{*}_{\tau_k}-\underline{Y}_{\tau_k})^{2} +2 \int_{t\wedge \tau_k}^{\tau_k}(\underline{Y}_s-Y^*_s)(Z^*_s-\underline{Z}_s)dB^*_s ,$$
where
$
(B^*_t:=B_t-\int_0^t\Delta H^*(s)ds)_{t\le T}$
is a Brownian motion under $\P^*$. Thus for any $t\le T$ and $k\ge 1$,
 \begin{equation}\label{unisyy*}
 \mathbb{E}^{\P^*}[(Y^{*}_{t\wedge \tau_k}-\underline{Y}_{t \wedge \tau_k})^{2}]\leq \mathbb{E}^{\P^*}[(Y^{*}_{\tau_k}-\underline{Y}_{\tau_k})^{2}].
 \end{equation}
 But for any $\gamma \ge 1$ and $\tau$ stopping time, $
 \E[|Y^*_\tau|^{\gamma}+|\underline{Y}_\tau|^{\gamma}]\leq C$ and by Lemma $\ref{lemma21}$ there exists a constant $p>1$ such that $\E[L_T^p]<\infty$. Then
 there exists a constant $C$ such that for any stopping time $\t$, $\E^{\P^*}[(Y^{*}_{\tau}-\underline{Y}_{\tau})^{2}]\leq C. $ Consequently the process $(Y^{*}-\underline{Y})^2$ is of class [D] under the
 probability $\P^*$. Therefore (one can see e.g. \cite{DMM1}, Theorem 21, pp. 36)
 $$
 \E^{\P^*}[(Y^{*}_{\tau_k}-\underline{Y}_{\tau_k})^{2}]\rightarrow_k \E^{\P^*}[(Y^{*}_{T}-\underline{Y}_{T})^{2}]=0.
 $$
 Going back now to \eqref{unisyy*}, using Fatou's Lemma to obtain
 $\E^{\P^*}[(Y^{*}_{t}-\underline{Y}_{t})^{2}]=0$ for any $t\le T$. It implies that $Y^*=\underline{Y}$, $\P^*$ and $\P$-a.s. since the probabilities are equivalent. Thus we have also $Z^*=\underline{Z}$, $K^{*+}=\underline{K}^{+}$ and $K^{*-}=\underline{K}^{-}$, i.e. uniqueness.
 \\ The proof is now complete.
 \end{proof}
\par We are now ready to give the main result of this paper.
\begin{theorem}
Assume that $\mathbf{[H1]}$ to $\mathbf{[H5]}$ hold true and let $(Y^*,Z^*,K^{*+},K^{*-})$ be the solution of the double barrier reflected BSDE associated with $(g(x),H^*,L,U).$ Set $u^*:=(u^*(t,x,Z^*))_{t\leq T},$ $v^*:=(v^*(t,x,Z^*))_{t\leq T},$ $\sigma^*:=\inf\{t\geq0, Y^*_t=L_t\}\wedge T$ and finally $\tau^*:=\inf\{t\geq0, Y^*_t=U_t\}\wedge T,$ then,
\begin{equation}\label{y0j}
Y_0^*=J(u^*,\tau^*;v^*,\sigma^*).
\end{equation}
 In addition $(u^*,\tau^*;v^*,\sigma^*)$ is a saddle-point strategy for the mixed zero-sum stochastic differential game.
 \end{theorem}
 \begin{proof}
 According to Theorem $\ref{theo31}$ the double barrier reflected BSDE associated with $(g,H^*,L,U)$ has a solution $(Y_t^*,Z_t^*,K_t^{*+},K_t^{*-})_{t\leq T}$. Then for any $t\leq T$ we have,
 \begin{equation*}
 Y^*_t = g(x) + \int_{t}^{T} H^*(s,x,Z^*_s) \, \mathrm{d}{s} + (K^{*+}_T- K^{*+}_t) - (K^{*-}_T- K^{*-}_t) -  \int_{t}^{T}Z^*_s \, \mathrm{d}{B_s}.
 \end{equation*}
 Next for $k\geq 1$, let us set
 $$\delta_k=\inf\{s\geq 0, |L_s|+|U_s|+|Y^*_s|+\int_{0}^{s}|Z^*_r|^{2}\mathrm{d}{r}\geq \widetilde{\theta}+k \}\wedge T.$$
 where $\widetilde{\theta}=|L_0|+|U_0|+|Y^*_0|.$\\
 Now since $Y_0^*$ is $\mathcal{F}_0$-measurable,  it is a deterministic constant and then,
 \begin{equation*}
 \begin{aligned}
 Y_0^*&=\mathbb{E}^{(u^*,v^*)}[ Y_0^*],\\ & = \mathbb{E}^{(u^*,v^*)}\biggl[Y^*_{\tau^*\wedge\sigma^*\wedge\delta_k} +\int_{0}^{\tau^*\wedge\sigma^*\wedge\delta_k} H^*(s,x,Z^{*}_s) \, \mathrm{d}{s}+ K_{\tau^*\wedge\sigma^*\wedge\delta_k}^{*+}-K_{\tau^*\wedge\sigma^*\wedge\delta_k}^{*-}\biggr.\\&\biggl. \qquad\qquad\quad-  \int_{0}^{\tau^*\wedge\sigma^*\wedge\delta_k}Z^{*}_s \, \mathrm{d}{B_s}\biggr], \\ & = \mathbb{E}^{(u^*,v^*)}\biggl[Y^*_{\tau^*\wedge\sigma^*\wedge\delta_k} +\int_{0}^{\tau^*\wedge\sigma^*\wedge\delta_k} \Gamma(s,x,Z^{*}_s,u^*,v^*) \, \mathrm{d}{s}+ K_{\tau^*\wedge\sigma^*\wedge\delta_k}^{*+}-K_{\tau^*\wedge\sigma^*\wedge\delta_k}^{*-}\biggr.\\&\biggl. \qquad\qquad\quad-  \int_{0}^{\tau^*\wedge\sigma^*\wedge\delta_k}Z^{*}_s \, \mathrm{d}{B_s}^{(u^*,v^*)}\biggr].
 \end{aligned}
 \end{equation*}
The processes $K^{*+}$ and $K^{*-}$ are increasing only when $Y^*_t=L_t$ and $Y^*_t=U_t$ respectively. Hence they do not increase between $0$ and $\tau^*\wedge\sigma^*$ and then $K_{\tau^*\wedge\sigma^*\wedge\delta_k}^{*+}=K_{\tau^*\wedge\sigma^*\wedge\delta_k}^{*-}=0.$
On the other hand, the process $(\int_{0}^{t}Z^*_{s}\mathrm{d}{B_s}^{(u^*,v^*)})_{t\leq \delta_k}$ is a $\P^{(u^*,v^*)}$-martingale.  Then for any $k\geq1$, we have
\begin{equation*}
 Y_0^*=\mathbb{E}^{(u^*,v^*)}\biggl[Y^*_{\tau^*\wedge\sigma^*\wedge\delta_k} +\int_{0}^{\tau^*\wedge\sigma^*\wedge\delta_k} \Gamma(s,x,Z^{*}_s,u^*,v^*) \, \mathrm{d}{s}\biggr].
 \end{equation*}
 The sequence of stopping times $(\delta_k)_{k\geq1}$ is increasing, of stationary type converging to $T$, then by taking the limit when $k\longrightarrow +\infty$ we deduce that,
 \begin{equation*}
  Y_0^*=\mathbb{E}^{(u^*,v^*)}\biggl[Y^*_{\tau^*\wedge\sigma^*} +\int_{0}^{\tau^*\wedge\sigma^*} \Gamma(s,x,Z^{*}_s,u^*,v^*) \, \mathrm{d}{s}\biggr].
  \end{equation*}
  And as $$Y^*_{\tau^*\wedge\sigma^*}=U_{\tau^*}1_{[\tau^*<\sigma^*]}+L_{\sigma^*}1_{[\sigma^*\leq\tau^*<T]}+g(x)1_{[\tau^*=\sigma^*=T]},\quad \P^{(u^*,v^*)}-a.s..$$
  Then
  \begin{equation*}
  Y_0^*=J(u^*,\tau^*;v^*,\sigma^*).
  \end{equation*}
  Next let $v=(v_t)_{t\leq T}$ be an admissible control for $C_2$ and let $\sigma$ be a stopping time. Let us show that $Y_0^*\geq J(u^*,\tau^*;v,\sigma).$ Since $\P$ and $\P^{(u^*,v)}$ are equivalent probabilities on $(\Omega,\mathcal{F})$, we have,
  \begin{equation*}
   \begin{aligned}
   Y_0^*&=\mathbb{E}^{(u^*,v)}[ Y_0^*],\\ & = \mathbb{E}^{(u^*,v)}\biggl[Y^*_{\tau^*\wedge\sigma\wedge\delta_k} +\int_{0}^{\tau^*\wedge\sigma\wedge\delta_k} H^*(s,x,Z^{*}_s) \, \mathrm{d}{s}+ K_{\tau^*\wedge\sigma\wedge\delta_k}^{*+}-K_{\tau^*\wedge\sigma\wedge\delta_k}^{*-}\biggr.\\&\biggl. \qquad\qquad\quad-  \int_{0}^{\tau^*\wedge\sigma\wedge\delta_k}Z^{*}_s \, \mathrm{d}{B_s}\biggr], \\ & = \mathbb{E}^{(u^*,v)}\biggl[Y^*_{\tau^*\wedge\sigma\wedge\delta_k} +\int_{0}^{\tau^*\wedge\sigma\wedge\delta_k} \Gamma(s,x,Z^{*}_s,u^*,v) \, \mathrm{d}{s}+ K_{\tau^*\wedge\sigma\wedge\delta_k}^{*+}-K_{\tau^*\wedge\sigma\wedge\delta_k}^{*-}\biggr.\\&\biggl. \qquad\qquad-  \int_{0}^{\tau^*\wedge\sigma\wedge\delta_k}Z^{*}_s \, \mathrm{d}{B_s}^{(u^*,v)}+\int_{0}^{\tau^*\wedge\sigma\wedge\delta_k} (H^*(s,x,Z^{*}_s)-H(s,x,Z^{*}_s,u^*,v)) \, \mathrm{d}{s}\biggr].
   \end{aligned}
   \end{equation*}
   But $\P^{(u^*,v)}$-a.s., $K_{\tau^*\wedge\sigma^*\wedge\delta_k}^{*+}\geq0$ and through $(\ref{H*H})$ we have, $H^*(t,x,Z^{*}_t)-H(t,x,Z^{*}_t,u^*,v)\geq0,$ for any $t\leq T.$ On the other hand, the process $(\int_{0}^{t}Z^*_{s} \, \mathrm{d}{B_s}^{(u^*,v)})_{t\leq \delta_k}$ is a $\P^{(u^*,v)}$-martingale. Then
   \begin{equation*}
    Y_0^*\geq\mathbb{E}^{(u^*,v)}\biggl[Y^*_{\tau^*\wedge\sigma\wedge\delta_k} +\int_{0}^{\tau^*\wedge\sigma\wedge\delta_k} \Gamma(s,x,Z^{*}_s,u^*,v) \, \mathrm{d}{s}\biggr].
    \end{equation*}
    The sequence of stopping times $(\delta_k)_{k\geq1}$ is increasing, of stationary type converging to $T$, then by taking the limit when $k\longrightarrow +\infty$ we deduce that,
     \begin{equation*}
      Y_0^*\geq\mathbb{E}^{(u^*,v)}\biggl[Y^*_{\tau^*\wedge\sigma} +\int_{0}^{\tau^*\wedge\sigma} \Gamma(s,x,Z^{*}_s,u^*,v) \, \mathrm{d}{s}\biggr].
      \end{equation*}
      And as $$Y^*_{\tau^*\wedge\sigma}\geq U_{\tau^*}1_{[\tau^*<\sigma]}+L_{\sigma}1_{[\sigma\leq\tau^*<T]}+g(x)1_{[\tau^*=\sigma=T]},\quad \P^{(u^*,v)}-a.s..$$
      Then
      \begin{equation*}
      Y_0^*\geq J(u^*,\tau^*;v,\sigma).
      \end{equation*}
      In the same way we can show that $Y_0^*\leq J(u,\tau;v^*,\sigma^*)$ for any stopping time $\tau$ and $u$ an admissible control for
      $C_1.$ It follows that $(u^*,\tau^*;v^*,\sigma^*)$ is a saddle-point strategy for the mixed zero-sum stochastic differential game.\\
      The proof is complete.
 \end{proof}
In the previous result we characterized $Y_0$ as the value of the game. However, more
can be said about the process $(Y^{*}_t )_{t\leq T}$. Actually we have the following proposition.
\begin{remark} The process $(Y^{*}_t )_{t\leq T}$ is the value function of the game, i.e., for any stopping time $\upsilon$ we have
	\begin{equation}
	\begin{aligned}
Y^{*}_{\upsilon}&=\essinf_{u\in \U;\tau\in\mathcal{ T}_{\upsilon}}\underset{v\in \V;\sigma\in\mathcal{ T}_{\upsilon}}{\esssup}\, J(u,\tau;v,\sigma),\\
& =\underset{v\in \V;\sigma\in\mathcal{ T}_{\upsilon}}{\esssup}\essinf_{u\in \U;\tau\in\mathcal{ T}_{\upsilon} }\, J(u,\tau;v,\sigma).
	\end{aligned}
	\end{equation}
\end{remark}
 \section{Viscosity solutions  of the mixed zero-sum stochastic differential game problem}
 \subsection{Connection with doubly reflected BSDEs  with double obstacle variational inequality}
 Let $\sigma:[0,T]\times\mathbb{R}^d\to\mathbb{R}^{d\times d}$ is a bounded continuous function, lipschitz w.r.t. the second variable, uniformly w.r.t. $t\in [0,T],$  invertible and its inverse $\sigma^{-1}(t,x)$ is bounded and continuous. For $(t,x)\in[0.T]\times\mathbb{R}^d,$ let $X^{t,x}:=(X^{t,x}_{s})_{s\leq T}$ be the unique $\mathbb{R}^d$-valued process solution of the following standard SDE:
 \begin{equation}
 \begin{cases}
 dX^{t,x}_{s}=\sigma(s,X^{t,x}_{s})dB_s,\qquad t\leq s\leq T;\\
 X^{t,x}_{s}=x, \qquad\qquad\qquad\qquad s \leq t.
 \end{cases}
 \end{equation}
  In the following result we collect some properties of $X^{t,x}.$
  \begin{proposition}\label{estimxx} (see e.g. $\cite{RY}$) The process $X^{t,x}$ satisfies the following estimates:
  \begin{itemize}
  \item [$(i)$] For any $q\geq 2$, there exists a constant $C$ such that
  \begin{equation}\label{estimat1}
  \mathbb{E}[\sup_{0\leq s\leq T}|X^{t,x}_s|^q]\leq C(1+|x|^q).
  \end{equation}
  \item[$(ii)$] There exists a constant $C$ such that for any $t,t'\in [0,T]$ and $x,x'\in \mathbb{R}^d$,
  \begin{equation}\label{estimat2}
  \mathbb{E}[\sup_{0\leq s\leq T}|X^{t,x}_s-X^{t',x'}_s|^2]\leq
  C(1+|x|^2)(|x-x'|^2+|t-t'|).
  \end{equation}$\qquad\qquad\qquad\qquad\qquad\qquad\qquad\qquad\quad\qquad\qquad\qquad\qquad\qquad\qquad\qquad\qquad\qquad\qquad\Box$
  \end{itemize}
  \end{proposition}
 Now let us consider the functions:
\begin{enumerate}
\item [(i)]- $h,h':[0,T]\times \mathbb{R}^d \rightarrow \mathbb{R}$ and $g:\mathbb{R}^d\rightarrow \mathbb{R}$. They are continuous and of polynomial growth i.e., for some positive constants $C$ and $p$ $$
|g(x)|+|h(t,x)|+|h'(t,x)|\le C(1+|x|^p).
$$ In addition for any $(t,x)\in[0,T]\times\mathbb{R}^d$
$$h(t,x)<h'(t,x)\quad and \quad h(T,x)\leq g(x)\leq h'(T,x).$$
\item[(ii)]- $\Gamma:[0,T]\times\mathbb{R}^d\times A \times B \longmapsto\mathbb{R}^{+}$ is continuous and for each $(t,x)$ the mapping $(a,b)\in A \times B\mapsto \Gamma(t,x,a,b)$
is  uniformly continuous, uniformly in $(t,x)$. In addition $\Gamma$ is of polynomial growth i.e., for some positive constants $C$ and $p$
\begin{equation}\label{lnr_growth}
|\Gamma(t,x,a,b)|\le C(1+|x|^p).
\end{equation}
\item[(iii)]- $f:[0,T]\times\mathbb{R}^d\times A \times B \longmapsto\mathbb{R}^{d}$ is continuous and the mapping $(a,b)\in A \times B\mapsto f(t,x,a,b)$ is  continuous,  uniformly w.r.t. $t$  and Lipschitz w.r.t. $x$. Moreover it is of linear growth which means that it satisfies the previous inequality $\eqref{lnr_growth}$ with $p=1$.
\end{enumerate}
Next let us define the quadruple $(Y^{t,x}_s,Z^{t,x}_s,K^{+,t,x}_s,K^{-,t,x}_s)\in \mathbb{R}^{1+d+1+1},$ $\forall s\in[t,T]$, as the unique solution of the following double barrier reflected BSDE: $\forall s\in[t,T]$,
\begin{equation}\label{RBSDExts}
\begin{cases}
\begin{aligned}
Y^{t,x}_s = g(X^{t,x}_T) + \int_{s}^{T}& H^*(r,X^{t,x}_r,Z^{t,x}_r) \, \mathrm{d}{r} + (K_T^{+,t,x} - K_s^{+,t,x}) - ( K_T^{-,t,x} - K_s^{-,t,x})\\ & -  \int_{s}^{T}Z^{t,x}_r \, \mathrm{d}{B_r}, \quad\forall s\leq T;
\end{aligned}\\
\forall s\in[t,T], \quad  h(s,X^{t,x}_s)\leq Y^{t,x}_s\leq h'(s,X^{t,x}_s),\\
\int_{0}^{T} (Y^{t,x}_s - h(s,X^{t,x}_s)) \, \mathrm{d}{K_s^{+,t,x}} = 0 \; and \; \int_{0}^{T} (h'(s,X^{t,x}_s)-Y^{t,x}_s) \, \mathrm{d}{K_s^{-,t,x}} = 0.
\end{cases}
\end{equation}
where $H^*$ is the same as in $\eqref{H*}$ with the new functions $f$, $\sigma$ and $\Gamma$ introduced above. By Theorem $\ref{theo31}$ the quadruple $(Y^{t,x},Z^{t,x},K^{+,t,x},K^{-,t,x})$ exists and is unique. On the other hand, for $(t,x)\in [0,T]\times \mathbb{R}^d$, $Y^{t,x}_t$ is  deterministic (see e.g. \cite{HM}). Then let us define the function $u(t,x)$ by
\begin{equation}\label{uy}
u(t,x)=Y^{t,x}_t.
\end{equation}
Now to proceed we are going to provide some properties of the function $u$ which we need later.
\begin{lemma}\label{lemma41}
The function $u$ is of polynomial growth.
\end{lemma}
\begin{proof}
Let  $(y^{t,x}_s,z^{t,x}_s,k^{+,t,x}_s,k^{-,t,x}_s)_{s\leq T}$ be the processes defined as follows:
\begin{equation}\label{RBSDEyxts}
\begin{cases}
(y^{t,x}, k^{+,t,x},k^{-,t,x})\in \R^{1+1+1}, \mbox{ are continuous} \mbox{ and }
\P-a.s., \,\, z^{t,x}(\omega) \mbox{ is } dt\mbox{-square integrable} ;\\ \,k^{\pm,t,x}
\mbox{ is  non-decreasing and }k^{\pm,t,x}_0=0,\\
\begin{aligned}
y^{t,x}_s = g(X^{t,x}_T)& + \int_{s}^{T} H^*(r,X^{t,x},z^{t,x}_r)1_{[r\geq t]} \, \mathrm{d}{r} + (k_T^{+,t,x} - k_s^{+,t,x})\\ &  - ( k_T^{-,t,x} - k_s^{-,t,x}) -  \int_{s}^{T}z^{t,x}_r \, \mathrm{d}{B_r},\quad\forall s\leq T;
\end{aligned}\\
\forall s\leq T, \quad  h(t\vee s,X_{t\vee s}^{t,x})\leq y^{t,x}_s\leq h'(t\vee s,X_{t\vee s}^{t,x}),\\
\int_{0}^{T} (y^{t,x}_s - h(s,X^{t,x}_s)) \, \mathrm{d}{k_s^{+,t,x}} = 0 \; and \; \int_{0}^{T} (h'(s,X^{t,x}_s)-Y^{t,x}_s) \, \mathrm{d}{k_s^{-,t,x}} = 0.
\end{cases}
\end{equation}
The existence of $(y^{t,x},z^{t,x},k^{+,t,x},k^{-,t,x})$, is obtained by Theorem \ref{theo31}. On the other hand, thanks to uniqueness of the solution of $\eqref{RBSDEyxts}$, for any $s\in [t,T]$, we have $y^{t,x}_s=Y^{t,x}_s$ and for any $s\in [0,t],$ $y^{t,x}_s=y^{t,x}_t=Y^{t,x}_t,$  $z^{t,x}_s=0$ and $k^{\pm,t,x}_s=0$.\\
\par Next let $\P^{(u,v)}$ (one should say
$\P^{t,x,u,v}$ but we omit the dependence w.r.t. $(t,x)$ when there is no confusion) be the probability, equivalent to $\P$, defined as follows: $$d\P^{(u,v)}=M^{t,x,u,v}_T d\P,$$ where for any $t\le T$,
\begin{equation*}
\begin{aligned}
M^{t,x,u,v}_T:=&\exp\biggl\{\int_0^T\sigma^{-1}(r,X^{t,x}_r)f(r,X^{t,x}_r,u_r,v_r)1_{[r\geq t]}dB_r \biggr.\\
& \biggl. -\frac{1}{2}\int_0^T||\sigma^{-1}(r,X^{t,x}_r)f(r,X^{t,x}_r,u_r,v_r)||^2 1_{[r\geq t]}dr \biggr\}.
\end{aligned}
\end{equation*}
Under $\P^{(u,v)}$,  $X^\tx$ is a weak solution of the following SDE:
\begin{equation}\label{dyn}
\left\{
\begin{array}{l}
dX^{t,x}_s=f(s,X^{t,x}_s,u_s,v_s)ds+\sigma(s,X^{t,x}_s)dB^{(u,v)},~~s\in[t,T];\\
X^{t,x}_s=x, \qquad\qquad\qquad\qquad s\leq t.
\end{array}\right.
\end{equation}
Then it verifies the following estimate: $\fr q\ge 2$,
\begin{equation}\label{estimat1x}
\overline{\E}^{(t,x)}[\sup_{0\le s\leq T}|X^{t,x}_s|^q]\leq C(1+|x|^q),
\end{equation}
where $\overline{\E}^{(t,x)}$ is the expectation under $\P^{(u,v)}$ and $C$ is a constant which does not depend on $u$, $v$, $t$ and $x$. Next by $(\ref{y0j})$ we have
\begin{equation}\label{y0tx}
\begin{aligned}
u(t,x)=y^{t,x}_{0}=&\, \inf_{u\in \U,\tau\in\mathcal{T}_0}\sup_{v\in \V,\sigma\in\mathcal{ T}_0}\, \overline{\E}^{(t,x)} \biggl[\int^{\tau\wedge\sigma}_{0}\, \Gamma(s,X^{t,x}_s,u_s,v_s)1_{[s\geq t]}\mathrm{d}{s} + h(\sigma\vee s,X^{t,x}_{\sigma\vee s})1_{[\sigma\leq\tau<T]}\biggr.\\
& \biggl.\, \quad +h'(\tau\vee s,X^{t,x}_{\tau\vee s})1_{[\tau<T]}+ g(X^{t,x}_T)1_{[\sigma=\tau =T]}\biggr].
\end{aligned}
\end{equation}
Now as $\Gamma$, $h$, $h'$ and $g$ are of polynomial growth then we easily deduce that $u(t,x)$ is also of polynomial growth.
\end{proof}
\begin{proposition}\label{continuous}
	The function $u$ is continuous.
\end{proposition}
\begin{proof}
	Let $\epsilon >0$ and $(t',x')\in B((t,x),\epsilon)$  (the ball centered in $(t,x)$ with radius $\epsilon$). Recall $\eqref{y0tx}$ and then
\begin{equation}\label{y0t'x'}
\begin{aligned}
y^{t',x'}_{0}&=\, \inf_{u\in \U,\tau\in\mathcal{T}_0}\sup_{v\in \V,\sigma\in\mathcal{ T}_0}\, \overline{\E}^{t',x'}\biggl[\int^{\tau\wedge\sigma}_{0}\, \Gamma(s,X^{t',x'}_s,u_s,v_s)1_{[s\geq t']}\mathrm{d}{s} + h(\sigma\vee s,X^{t',x'}_{\sigma\vee s})1_{[\sigma\leq\tau<T]}\biggr.\\
& \biggl.\quad +h'(\tau\vee s,X^{t',x'}_{\tau\vee s})1_{[\tau<T]}+ g(X^{t,x}_T)1_{[\sigma=\tau =T]}\biggr].\\ &  =
\, \inf_{u\in \U,\tau\in\mathcal{T}_0}\sup_{v\in \V,\sigma\in\mathcal{ T}_0}\, \mathbb{E} \biggl[M^{t',x',u,v}_T(\int^{\tau\wedge\sigma}_{0}\, \Gamma(s,X^{t',x'}_s,u_s,v_s)1_{[s\geq t']}\mathrm{d}{s} + h(\sigma\vee s,X^{t',x'}_{\sigma\vee s})1_{[\sigma\leq\tau<T]}\biggr.\\
& \biggl.\quad +h'(\tau\vee s,X^{t',x'}_{\tau\vee s})1_{[\tau<T]}+ g(X^{t',x'}_T)1_{[\sigma=\tau =T]})\biggr].
\end{aligned}
\end{equation}
\par First let us show that the function $u$ is upper semi-continuous.\\
Fix an arbitrary $\epsilon>0$, and we assume that $t'<t$ . Pick $(v^{\epsilon},\sigma^{\epsilon})\in \V\times\mathcal{ T}_0$  such that
\begin{equation*}
\begin{aligned}
y^{t',x'}_{0}\leq&\,\mathbb{E} \biggl[M^{t',x',u,v^{\epsilon}}_T(\int^{\tau\wedge\sigma^{\epsilon}}_{0}\, \Gamma(s,X^{t',x'}_s,u_s,v^{\epsilon}_s)1_{[s\geq t']}\mathrm{d}{s} + h(\sigma^{\epsilon}\vee s,X^{t',x'}_{\sigma^{\epsilon}\vee s})1_{[\sigma^{\epsilon}\leq\tau<T]}\biggr.\\
& \biggl.\quad +h'(\tau\vee s,X^{t',x'}_{\tau\vee s})1_{[\tau<T]}+ g(X^{t',x'}_T)1_{[\sigma^{\epsilon}=\tau=T]})\biggr]+ \epsilon.
\end{aligned}
\end{equation*}
where $(u,\tau)\in \U\times\mathcal{ T}_0$  will be chosen later. On the other hand, pick $(u^{\epsilon},\tau^{\epsilon})\in \U\times\mathcal{ T}_0$ such that
\begin{equation*}
\begin{aligned}
y^{t,x}_{0}\geq&\, \mathbb{E} \biggl[M^{t,x,u^{\epsilon},v^{\epsilon}}_T(\int^{\tau^{\epsilon}\wedge\sigma^{\epsilon}}_{0}\, \Gamma(s,X^{t,x}_s,u^{\epsilon}_s,v^{\epsilon}_s)1_{[s\geq t]}\mathrm{d}{s} + h(\sigma^{\epsilon}\vee s,X^{t,x}_{\sigma^{\epsilon}\vee s})1_{[\sigma^{\epsilon}\leq\tau^{\epsilon}<T]}\biggr.\\
& \biggl.\quad +h'(\tau^{\epsilon}\vee s,X^{t,x}_{\tau^{\epsilon}\vee s})1_{[\tau^{\epsilon}<T]}+ g(X^{t,x}_T)1_{[\sigma^{\epsilon}=\tau^{\epsilon} =T]})\biggr]- \epsilon.
\end{aligned}
\end{equation*}
Therefore
\begin{equation*}
\begin{aligned}
y^{t',x'}_0-y^{t,x}_0 &\leq\mathbb{E}[M^{t',x',u^{\epsilon},v^{\epsilon}}_T \int_{0}^{\tau^{\epsilon}\wedge\sigma^{\epsilon}}\{| \Gamma(s,X^{t',x'}_s,u^{\epsilon}_s,v^{\epsilon}_s)-\Gamma(s,X^{t,x}_s,u^{\epsilon}_s,v^{\epsilon}_s)|1_{[s\geq
t]}\\{}&+|\Gamma(s,X^{t',x'}_s,u^{\epsilon}_s,v^{\epsilon}_s)|1_{[t'\leq s< t]}\}ds+M^{t',x',u^{\epsilon},v^{\epsilon}}_T|h(\sigma^{\epsilon}\vee s,X^{t,x}_{\sigma^{\epsilon}\vee s})-h(\sigma^{\epsilon}\vee s,X^{t',x'}_{\sigma^{\epsilon}\vee s})|1_{[\sigma^{\epsilon}\leq\tau^{\epsilon}<T]}\\{}&+M^{t',x',u^{\epsilon},v^{\epsilon}}_T|h'(\tau^{\epsilon}\vee s,X^{t,x}_{\tau^{\epsilon}\vee s})-h'(\tau^{\epsilon}\vee s,X^{t',x'}_{\tau^{\epsilon}\vee s})|1_{[\tau^{\epsilon}<T]} +M^{t',x',u^{\epsilon},v^{\epsilon}}_T|g(X^{t,x}_T)\\&-g(X^{t',x'}_T)|1_{[\sigma^{\epsilon}=\tau^{\epsilon}=T]}+|M^{t,x,u^{\epsilon},v^{\epsilon}}_T-M^{t',x',u^{\epsilon},v^{\epsilon}}_T||(\int_{0}^{\tau^{\epsilon}\wedge\sigma^{\epsilon}}\Gamma(s,X^{t,x}_s,u^{\epsilon}_s,v^{\epsilon}_s)1_{[s\geq t]}ds\\&+h(\sigma^{\epsilon}\vee s,X^{t,x}_{\sigma^{\epsilon}\vee s})1_{[\sigma^{\epsilon}\leq\tau^{\epsilon}<T]}+h'(\tau^{\epsilon}\vee s,X^{t,x}_{\tau^{\epsilon}\vee s})1_{[\tau^{\epsilon}<T]}+g(X^{t,x}_T)1_{[\tau^{\epsilon}=\sigma^{\epsilon}=T]})|]+2\epsilon,\\& \leq I_1 + I_2 + I_3+2\epsilon,
\end{aligned}
\end{equation*}
where
\begin{equation}\label{inegalite}
\begin{aligned}
&I_1:=\mathbb{E}[M^{t',x',u^{\epsilon},v^{\epsilon}}_T \int_{0}^{T}\{| \Gamma(s,X^{t',x'}_s,u^{\epsilon}_s,v^{\epsilon}_s)-\Gamma(s,X^{t,x}_s,u^{\epsilon}_s,v^{\epsilon}_s)|+|\Gamma(s,X^{t',x'}_s,u^{\epsilon}_s,v^{\epsilon}_s)|1_{[t'\leq s< t]}\}ds].\\{}& I_2:=\mathbb{E}[M^{t',x',u^{\epsilon},v^{\epsilon}}_T\sup_{s\leq T}|h(t\vee s,X^{t,x}_{t\vee s})-h(t'\vee s,X^{t',x'}_{t'\vee s})|\\&\qquad\qquad+M^{t',x',u^{\epsilon},v^{\epsilon}}_T\sup_{s\leq T}|h'(t\vee s,X^{t,x}_{t\vee s})-h'(t'\vee s,X^{t',x'}_{t'\vee s})|+M^{t',x',u^{\epsilon},v^{\epsilon}}_T|g(X^{t,x}_T)-g(X^{t',x'}_T)|].\\{}& I_3:=\mathbb{E}[|M^{t,x,u^{\epsilon},v^{\epsilon}}_T-M^{t',x',u^{\epsilon},v^{\epsilon}}_T||(\int_{0}^{\tau^{\epsilon}\wedge\sigma^{\epsilon}}\Gamma(s,X^{t,x}_s,u^{\epsilon}_s,v^{\epsilon}_s)1_{[s\geq t]}ds+h(\sigma^{\epsilon}\vee s,X^{t,x}_{\sigma^{\epsilon}\vee s})1_{[\sigma^{\epsilon}\leq\tau^{\epsilon}<T]}\\&\qquad\qquad +h'(\tau^{\epsilon}\vee s,X^{t,x}_{\tau^{\epsilon}\vee s})1_{[\tau^{\epsilon}<T]}+g(X^{t,x}_T)1_{[\tau^{\epsilon}=\sigma^{\epsilon}=T]})|].
\end{aligned}
\end{equation}
In the right-hand side of $\eqref{inegalite}$ the first term $(I_1)$ converges to $0$ as $(t',x')\rightarrow (t,x)$ since $\Gamma$ is continuous and of polynomial growth and there exists $p>1$ such that $\E[(M^{t',x',u^{\epsilon},v^{\epsilon}}_T)^p]\le C$ where $C$ does not depend on $t',x'$ (give a sharp estimate of this latter).\\ Next let us deal with the term $I_2$ and let us show first that,
\begin{equation*}
\mathbb{E}[M^{t',x',u^{\epsilon},v^{\epsilon}}_T\sup_{s\leq T}|h(t\vee s,X^{t,x}_{t\vee s})-h(t'\vee s,X^{t',x'}_{t'\vee s})|]\rightarrow 0 \mbox{ as }(t',x')\rightarrow (t,x),
\end{equation*}
\begin{equation*}
\mathbb{E}[M^{t',x',u^{\epsilon},v^{\epsilon}}_T\sup_{s\leq T}|h'(t\vee s,X^{t,x}_{t\vee s})-h'(t'\vee s,X^{t',x'}_{t'\vee s})|]\rightarrow 0 \mbox{ as }(t',x')\rightarrow (t,x),
\end{equation*}
and
\begin{equation*}
\mathbb{E}[M^{t',x',u^{\epsilon},v^{\epsilon}}_T|g(X^{t,x}_T)-g(X^{t',x'}_T)|]\rightarrow 0 \mbox{ as }(t',x')\rightarrow (t,x).
\end{equation*}
By Holder's inequality we have,
\begin{equation*}
\begin{aligned}
\mathbb{E}[M^{t',x',u^{\epsilon},v^{\epsilon}}_T&\sup_{s\leq T}|h(t\vee s,X^{t,x}_{t\vee s})-h(t'\vee s,X^{t',x'}_{t'\vee s})|],\\ & \leq \mathbb{E}[(M^{t',x',u^{\epsilon},v^{\epsilon}}_T)^{p}]^{\frac{1}{p}}\mathbb{E}[\sup_{s\leq T}|h(t\vee s,X^{t,x}_{t\vee s})-h(t'\vee s,X^{t',x'}_{t'\vee s})|^{q}]^{\frac{1}{q}},
\\ & \leq C\mathbb{E}[\sup_{s\leq T}|h(t\vee s,X^{t,x}_{t\vee s})-h(t'\vee s,X^{t',x'}_{t'\vee s})|^{q}]^{\frac{1}{q}},
\end{aligned}
\end{equation*}
where $\frac{1}{p}+\frac{1}{q}=1$ and $C$ a constant which does not depend on $t,x$ ($p>1$ is the constant which appears in Lemma $\ref{lemma21}$). Now for any $\rho>0$ we have:
\begin{equation*}
\begin{aligned}
|h(t'&\vee s,X^{t',x'}_{t'\vee s})-h(t\vee s, X^{t,x}_{t\vee s})|^{q}\leq
|h(t'\vee s ,X^{t',x'}_{t'\vee s }) -h(t\vee s
,X^{t',x'}_{t'\vee s})|^q1_{[|X^{t',x'}_{t'\vee s }|\leq \rho]}\\ &+
|h(t'\vee s ,X^{t',x'}_{t'\vee s }) -h(t\vee
s ,X^{t',x'}_{t'\vee s})|^q1_{[|X^{t',x'}_{t'\vee s }|\geq
\rho]}+|h(t\vee s ,X^{t',x'}_{t'\vee s }) -h(t\vee s
,X^{t,x}_{t\vee s})|^q.
\end{aligned}
\end{equation*}
Therefore we have:
\begin{equation*}
\begin{aligned}
\mathbb{E}[\sup_{s\leq T}|h(t'\vee s,&X^{t',x'}_{t'\vee s})-h(t\vee s, X^{t,x}_{t\vee s})|^{q}]\leq\mathbb{E}[\sup_{s\leq T}\{
|h(t'\vee s ,X^{t',x'}_{t'\vee s }) -h(t\vee s
,X^{t',x'}_{t'\vee s})|^q1_{[|X^{t',x'}_{t'\vee s }|\leq \rho]}\}]\\ &
+ \mathbb{E} [\sup_{s\leq T}\{|h(t'\vee s ,X^{t',x'}_{t'\vee s })
-h(t\vee s ,X^{t',x'}_{t'\vee
s})|^q\}1_{[\sup_{s\leq T}|X^{t',x'}_{t'\vee s }|\geq \rho]}]\\ &
+ \mathbb{E} [\sup_{s\leq T}\{|h(t\vee s ,X^{t',x'}_{t'\vee s })
-h(t\vee s ,X^{t,x}_{t\vee s})|^q\}1_{[\sup_{s\leq
T}|X^{t',x'}_{t'\vee s }|+\sup_{s\leq T}|X^{t,x}_{t\vee s }|\geq
\rho]}]\\ & +\mathbb{E}[\sup_{s\leq T}\{|h(t\vee s
,X^{t',x'}_{t'\vee s }) -h(t\vee s ,X^{t,x}_{t\vee
s})|^q\}1_{[\sup_{s\leq T}|X^{t',x'}_{t'\vee s }|+\sup_{s\leq
T}|X^{t,x}_{t\vee s }|\leq \rho]}].
\end{aligned}
\end{equation*}
But since $h$ is continuous then it is uniformly continuous on
$[0,T]\times \{x\in \mathbb{R}^k, |x|\leq \rho\}.$ Henceforth for any
$\epsilon_1>0$ there exists $\eta_{\epsilon_1}>0$ such that for any
$|t-t'|<\eta_{\epsilon_1}$ we have:
\begin{equation}\label{term1}
\sup_{s\leq T}\{|h(t'\vee s ,X^{t',x'}_{t'\vee s }) -h(t\vee s
,X^{t',x'}_{t'\vee s})|1_{[|X^{t',x'}_{t'\vee s }|\leq
\rho]}\}\leq \epsilon_1
\end{equation}
Next using Cauchy-Schwarz's inequality and then Markov's one with the second term we obtain:
\begin{equation}\label{term2}
\mathbb{E}[\sup_{s\leq T}\{|h(t'\vee s ,X^{t',x'}_{t'\vee s })
-h(t\vee s ,X^{t',x'}_{t'\vee s})|^q\}1_{[\sup_{s\leq
T}|X^{t',x'}_{t'\vee s }|\geq \rho]}]\leq
C(1+|x'|^{q_1})\rho^{-\frac{1}{2}},
\end{equation}
where $C$ and $q_1$ are real constants (due to the the polynomial growth of $h$ and estimate $(\ref{estimat1}))$. In the same way we have:
\begin{equation}\label{I2term3}
\mathbb{E}[\sup_{s\leq T}\{|h(t\vee s ,X^{t',x'}_{t'\vee s })
-h(t\vee s ,X^{t,x}_{t\vee s})|^q\}1_{[\sup_{s\leq
T}|X^{t',x'}_{t'\vee s }|+\sup_{s\leq T}|X^{t,x}_{t\vee s }|\geq
\rho]}]\leq C(1+|x|^p+|x'|^{q_1}) \rho^{-\frac{1}{2}}.
\end{equation}
Finally
using the uniform continuity of $h$ on compact subsets, the
continuity property (\ref{estimat2}) and the Lebesgue dominated
convergence theorem we obtain that
\begin{equation}\label{term4}
\mathbb{E}[\sup_{s\leq
T}\{|h(t\vee s ,X^{t',x'}_{t'\vee s }) -h(t\vee s
,X^{t,x}_{t\vee s})|^q\}1_{[\sup_{s\leq T}|X^{t',x'}_{t'\vee s
}|+\sup_{s\leq T}|X^{t,x}_{t\vee s }|\leq \rho]}] \rightarrow 0 \mbox{ as }(t',x')\rightarrow(t,x).
\end{equation}
Taking now into account
$\eqref{term1}$-$\eqref{term4}$ we have:
\begin{equation*}
\limsup_{(t',x')\rightarrow (t,x)}\mathbb{E}[\sup_{s\leq T}|h(t'\vee s ,X^{t',x'}_{t'\vee s })
-h(t\vee s ,X^{t,x}_{t\vee s})|^q]\leq
\epsilon_1+C(1+|x|^{q_1})\rho^{-\frac{1}{2}}.
\end{equation*}
 As $\epsilon_1$ and $\rho$ are arbitrary then making $\epsilon_1 \rightarrow  0$ and
$\rho \rightarrow +\infty$ to obtain that:
\begin{equation*}
\lim_{(t',x')\rightarrow (t,x)}\mathbb{E}[\sup_{s\leq T}|h(t'\vee s ,X^{t',x'}_{t'\vee s })
-h(t\vee s ,X^{t,x}_{t\vee s})|^q]=0.
\end{equation*}
Now since there exists a constant $C$  such that:
\begin{equation}\label{MM}
\mathbb{E}[(M^{t',x',u^{\epsilon},v^{\epsilon}}_T)^p])^{\frac{1}{p}}\leq C,
\end{equation}
then
\begin{equation*}
\mathbb{E}[M^{t',x',u^{\epsilon},v^{\epsilon}}_T\sup\limits_{s\leq T}| h(t\vee s,X_{t\vee
s}^{t,x})-h(t'\vee s,X_{t'\vee s}^{t',x'})|]\rightarrow 0 \mbox{ as }(t',x')\rightarrow (t,x).
\end{equation*}
By the same argument as above, we have
\begin{equation*}
\mathbb{E}[M^{t',x',u^{\epsilon},v^{\epsilon}}_T\{\sup_{s\leq T}|h'(t\vee s,X^{t,x}_{t\vee s})-h'(t'\vee s,X^{t',x'}_{t'\vee s})|+|g(X^{t,x}_T)-g(X^{t',x'}_T)|\}]\rightarrow 0 \mbox{ as }(t',x')\rightarrow (t,x).
\end{equation*}
Thus the claim is proved.\\
Finally let us focus on the last term $(I_3)$ in $\eqref{inegalite}$. From  the polynomial growth of $\Gamma$, $h$, $h'$ and $g$ we have:
\begin{equation}\label{term3}
I_3\leq C\mathbb{E}[|M^{t,x,u^{\epsilon},v^{\epsilon}}_T-M^{t',x',u^{\epsilon},v^{\epsilon}}_T|(1+ \sup\limits_{s\leq T}|X^{t,x}_s|^p)],
\end{equation}
where $C$ and $p$ are appropriate constants which stem from the polynomial growth of $\Gamma$, $h$, $h'$ and $g$. Now let us define the  stopping time $\delta$ by:
$$\delta=\inf\{s\in [0,T], \,\sup\limits_{r\leq s}|X^{t,x}_r|+\sup\limits_{r\leq s}|X^{t',x'}_r|\geq \rho \}\wedge T.$$
Note that $\delta$ depends actually on $\rho$, $t$, $x$, $t'$ and $x'$. On the other hand since $(t',x')\in B((t,x),\epsilon)$, then taking into account of $\eqref{estimat1}$ it holds
\begin{equation}\label{estg}
\P[\delta <T]\leq \P[\sup\limits_{r\leq T}|X^{t,x}_r|+\sup\limits_{r\leq T}|X^{t',x'}_r|\geq \rho]\le C_{t,x,\epsilon}\rho^{-1}.
\end{equation}
where $C_{t,x,\epsilon}$ is a constant which depends on $t,$ $x,$ $\epsilon.$ Next we deal with the right-hand side of \eqref{term3}, i.e.,
 We have,
 \begin{equation*}
 \begin{aligned}
 &\mathbb{E}[|M^{t,x,u^{\epsilon},v^{\epsilon}}_T-M^{t',x',u^{\epsilon},v^{\epsilon}}_T|(1+ \sup\limits_{s\leq T}|X^{t,x}_s|^p)]\;\leq \mathbb{E}[(|M^{t,x,u^{\epsilon},v^{\epsilon}}_T-M^{t',x',u^{\epsilon},v^{\epsilon}}_T|\\ &-|M^{t,x,u^{\epsilon},v^{\epsilon}}_{\delta}-M^{t',x',u^{\epsilon},v^{\epsilon}}_{\delta}|)(1+ \sup\limits_{s\leq T}|X^{t,x}_s|^p)1_{\delta< T}]+  \mathbb{E}[|M^{t,x,u^{\epsilon},v^{\epsilon}}_{\delta}-M^{t',x',u^{\epsilon},v^{\epsilon}}_{\delta}|(1+ \sup\limits_{s\leq T}|X^{t,x}_s|^p)],\\ & \leq C_{t,x}(\P[\delta<T])^{q_0}+C(1+|x|^p)\E[|M^{t,x,u}_\delta-M^{t',x',u}_\delta|^{q}]^{\frac{1}{q}}.
 \end{aligned}
 \end{equation*}
 where $C_{\tx}$ is a constant which depend on $\tx$, $q_0>0$. To obtain the last inequality we have used Holder's inequalities (twice). Therefore using \eqref{estg} yields:
 \begin{equation}\label{etapeinterm1}
 \E[|M^{t,x,u^{\epsilon},v^{\epsilon}}_T-M^{t',x',u^{\epsilon},v^{\epsilon}}_T|(1+ \sup\limits_{s\leq T}|X^{t,x}_s|^p)]\leq \frac{C_{t,x,\epsilon}}{\rho^{q_0}}+C(1+|x|^p)\E[|M^{t,x,u^{\epsilon},v^{\epsilon}}_\delta-M^{t',x',u^{\epsilon},v^{\epsilon}}_\delta|^{q}]^{\frac{1}{q}}.
 \end{equation}
 We are now going to show that
 \begin{equation*}
 \E[|M^{t,x,u^{\epsilon},v^{\epsilon}}_\delta-M^{t',x',u^{\epsilon},v^{\epsilon}}_\delta|^{q}]^{\frac{1}{q}}\rw 0 \mx{ as }\ttrw.
 \end{equation*}
For any $s\le T$,
$$M^{t,x,u^{\epsilon},v^{\epsilon}}_{ s\wedge \delta }= 1+ \int_0^{ s\wedge \delta } M^{t,x,u^{\epsilon},v^{\epsilon}}_r \sigma^{-1}(r,X^{t,x}_r)f(r,X^{t,x}_r,u^{\epsilon}_r,v^{\epsilon}_r)1_{[r\geq t]}dB_r,$$
and  a similar equation for $M^{t',x',u^{\epsilon},v^{\epsilon}}$.
By using the Itô's formula and taking expectation we obtain: $\fr s\le T$,
\begin{equation*}
\begin{aligned}
\E[|M^{t,x,u^{\epsilon},v^{\epsilon}}_{ s\wedge \delta }-M^{t',x',u^{\epsilon},v^{\epsilon}}_{ s\wedge \delta }|^2]=& \E\biggl[\int_0^{s\wedge\delta} (M^{t,x,u^{\epsilon},v^{\epsilon}}_r \sigma^{-1}(r,X^{t,x}_r)f(r,X^{t,x}_r,u^{\epsilon}_r,v^{\epsilon}_r)1_{[r\geq t]}\\ &-M^{t',x',u^{\epsilon},v^{\epsilon}}_r \sigma^{-1}(r,X^{t',x'}_r)f(r,X^{t',x'}_r,u^{\epsilon}_r,v^{\epsilon}_r)1_{[r\geq t']})^2\;dr\biggr].
\end{aligned}
\end{equation*}
We then have,
\begin{equation*}
\begin{aligned}
&\E[|M^{t,x,u^{\epsilon},v^{\epsilon}}_{ s\wedge \delta }-M^{t',x',u^{\epsilon},v^{\epsilon}}_{ s\wedge \delta }|^2]\\&= \E\biggl[\int_0^{s\wedge\delta} \{(M^{t,x,u^{\epsilon},v^{\epsilon}}_r-M^{t',x',u^{\epsilon},v^{\epsilon}}_r) \sigma^{-1}(r,X^{t,x}_r)f(r,X^{t,x}_r,u^{\epsilon}_r,v^{\epsilon}_r)1_{[r\geq t]}\\& +M^{t',x',u^{\epsilon},v^{\epsilon}}_r\sigma^{-1}(r,X^{t,x}_r)f(r,X^{t,x}_r,u^{\epsilon}_r,v^{\epsilon}_r)1_{[r\geq t]}-M^{t',x',u^{\epsilon},v^{\epsilon}}_r \sigma^{-1}(r,X^{t',x'}_r)f(r,X^{t',x'}_r,u^{\epsilon}_r,v^{\epsilon}_r)1_{[r\geq t']}\}^2\;dr\biggr],
\\ & \leq 2\E\biggl[\int_0^{s\wedge\delta} \{(M^{t,x,u^{\epsilon},v^{\epsilon}}_r-M^{t',x',u^{\epsilon},v^{\epsilon}}_r) \sigma^{-1}(r,X^{t,x}_r)f(r,X^{t,x}_r,u^{\epsilon}_r,v^{\epsilon}_r)1_{[r\geq t]}\}^2\;dr\biggr]\\& +2\E\biggl[\int_0^{s\wedge\delta} \{M^{t',x',u^{\epsilon},v^{\epsilon}}_r(\sigma^{-1}(r,X^{t,x}_r)f(r,X^{t,x}_r,u^{\epsilon}_r,v^{\epsilon}_r)1_{[r\geq t]}- \sigma^{-1}(r,X^{t',x'}_r)f(r,X^{t',x'}_r,u^{\epsilon}_r,v^{\epsilon}_r)1_{[r\geq t']})\}^2\;dr\biggr],
\\ & \leq 2\E\biggl[\int_0^{s} \{(M^{t,x,u^{\epsilon},v^{\epsilon}}_{r\wedge\delta}-M^{t',x',u^{\epsilon},v^{\epsilon}}_{r\wedge\delta}) \sigma^{-1}({r\wedge\delta},X^{t,x}_{r\wedge\delta})f({r\wedge\delta},X^{t,x}_{r\wedge\delta},u^{\epsilon}_{r\wedge\delta},v^{\epsilon}_{r\wedge\delta})1_{[{r\wedge\delta}\geq t]}\}^2\;dr\biggr]\\& +2\E\biggl[\int_0^{s\wedge\delta} \{M^{t',x',u^{\epsilon},v^{\epsilon}}_r(\sigma^{-1}(r,X^{t,x}_r)f(r,X^{t,x}_r,u^{\epsilon}_r,v^{\epsilon}_r)1_{[r\geq t]}- \sigma^{-1}(r,X^{t',x'}_r)f(r,X^{t',x'}_r,u^{\epsilon}_r,v^{\epsilon}_r)1_{[r\geq t']})\}^2\;dr\biggr],
\\ & \leq 2C_{\rho}\E\biggl[\int_0^{s} (M^{t,x,u^{\epsilon},v^{\epsilon}}_{r\wedge\delta}-M^{t',x',u^{\epsilon},v^{\epsilon}}_{r\wedge\delta})^2\;dr\biggr] +2\E\biggl[\int_0^{\delta} \{M^{t',x',u^{\epsilon},v^{\epsilon}}_r(\sigma^{-1}(r,X^{t,x}_r)f(r,X^{t,x}_r,u^{\epsilon}_r,v^{\epsilon}_r)1_{[r\geq t]}\\&\qquad-\sigma^{-1}(r,X^{t',x'}_r)f(r,X^{t',x'}_r,u^{\epsilon}_r,v^{\epsilon}_r)1_{[r\geq t']})\}^2\;dr\biggr].
\end{aligned}
\end{equation*}
Using now Gronwall's inequality to obtain: $\fr s\in [0,T]$,
\begin{equation}\label{limmu}
\E[|M^{t,x,u^{\epsilon},v^{\epsilon}}_{ s\wedge \delta }-M^{t',x',u^{\epsilon},v^{\epsilon}}_{ s\wedge \delta }|^2]\leq \Psi(t',x',\eps,t,x,u^{\epsilon},v^{\epsilon})\times \exp({C_\rho. s}),
\end{equation}
where
\begin{equation*}
\begin{aligned}
\Psi(t',x',\eps,t,x,u^{\epsilon},v^{\epsilon})=&2\E\biggl[\int_0^{\delta}\{M^{t',x',u^{\epsilon},v^{\epsilon}}_r(\sigma^{-1}(r,X^{t,x}_r)f(r,X^{t,x}_r,u^{\epsilon}_r,v^{\epsilon}_r)1_{[r\geq t]}\\&\qquad-\sigma^{-1}(r,X^{t',x'}_r)f(r,X^{t',x'}_r,u^{\epsilon}_r,v^{\epsilon}_r)1_{[r\geq t']})\}^2\;dr\biggr].
\end{aligned}
\end{equation*}
Finally by Fatou's Lemma we have:
\begin{equation}\label{limmu2}
\E[|M^{t,x,u^{\epsilon},v^{\epsilon}}_{\delta }-M^{t',x',u^{\epsilon},v^{\epsilon}}_{\delta }|^2]\leq \Psi(t',x',\eps,t,x,u^{\epsilon},v^{\epsilon})\times \exp({C_\rho. T}).
\end{equation}
But
\begin{equation*}
\begin{aligned}
&\Psi(t',x',\eps,t,x,u^{\epsilon},v^{\epsilon})\\&=2\E\biggl[\int_0^{\delta}\{M^{t',x',u^{\epsilon},v^{\epsilon}}_r(\sigma^{-1}(r,X^{t,x}_r)f(r,X^{t,x}_r,u^{\epsilon}_r,v^{\epsilon}_r)1_{[r\geq t]}-\sigma^{-1}(r,X^{t',x'}_r)f(r,X^{t',x'}_r,u^{\epsilon}_r,v^{\epsilon}_r)1_{[r\geq t']})\}^2\;dr\biggr]\\&\qquad+2\E\biggl[\int_0^{\delta}\{M^{t',x',u^{\epsilon},v^{\epsilon}}_r \sigma^{-1}(r,X^{t',x'}_r)f(r,X^{t',x'}_r,u^{\epsilon}_r,v^{\epsilon}_r)1_{[t>r\geq t']}\}^2\;dr\biggr].
\end{aligned}
\end{equation*}
Making use now of the Cauchy-Schwarz's inequality to obtain:
\begin{equation*}
\begin{aligned}
\Psi(t',x',\eps,t,&x,u^{\epsilon},v^{\epsilon})\leq 2\E\biggl[\int_0^{\delta}(M^{t',x',u^{\epsilon},v^{\epsilon}}_r)^41_{[r\geq t]}\;dr\biggr]^{\frac{1}{2}}\times \E\biggl[\int_0^{\delta}(\sigma^{-1}(r,X^{t,x}_r)f(r,X^{t,x}_r,u^{\epsilon}_r,v^{\epsilon}_r)\\&-\sigma^{-1}(r,X^{t',x'}_r)f(r,X^{t',x'}_r,u^{\epsilon}_r,v^{\epsilon}_r))^41_{[r\geq t]}\;dr\biggr]^{\frac{1}{2}}+2\E\biggl[\int_0^{\delta}(M^{t',x',u^{\epsilon},v^{\epsilon}}_r)^41_{[t>r\geq t']}\;dr\biggr]^{\frac{1}{2}}\times \\&\; \E\biggl[\int_0^{\delta}(\sigma^{-1}(r,X^{t',x'}_r)f(r,X^{t',x'}_r,u^{\epsilon}_r,v^{\epsilon}_r))^41_{[t>r\geq t']}\;dr\biggr]^{\frac{1}{2}}.
\end{aligned}
\end{equation*}
Next for $r\in [0,\delta]$, $|X^{t,x}_r|+|X^{t',x'}_r|\le \rho$, then $M^{t',x',u^{\epsilon}_r,v^{\epsilon}_r}_{r\wedge \delta}$ has moments of any order uniformly  (see Lemma $\ref{lemma21}$), which implies the existence of a constant $\kappa_\rho$ such that $$
\sup_{r,t',x'u^{\epsilon}_r,v^{\epsilon}_r}\E[(M^{t',x'u^{\epsilon}_r,v^{\epsilon}_r}_{r\wedge \delta})^4]\le \kappa_\rho.$$
On the other hand we have:
\begin{equation*}
\begin{aligned}
\E&\biggl[\int_0^{\delta}(\sigma^{-1}(r,X^{t,x}_r)f(r,X^{t,x}_r,u^{\epsilon}_r,v^{\epsilon}_r)-\sigma^{-1}(r,X^{t',x'}_r)f(r,X^{t',x'}_r,u^{\epsilon}_r,v^{\epsilon}_r))^41_{[r\geq t]}\;dr\biggr]\\&= \E\biggl[\int_0^{\delta}(\sigma^{-1}(r,X^{t,x}_r)f(r,X^{t,x}_r,u^{\epsilon}_r,v^{\epsilon}_r)-\sigma^{-1}(r,X^{t',x'}_r)f(r,X^{t',x'}_r,u^{\epsilon}_r,v^{\epsilon}_r))^41_{\{|X^{t,x}_r|+|X^{t',x'}_r|\le \rho\}}1_{[r\geq t]}\;dr\biggr].
\end{aligned}
\end{equation*}
Now the continuity property of $\sigma^{-1}f$ implies the existence of a continuous increasing bounded function $\Phi_\rho$ with $\Phi_\rho(0)=0$ (typically the modulus of continuity of $\sigma^{-1}f$ on $[0,T]\times \bar B(0,\rho)\times A\times B$) such that
\begin{equation*}
\begin{aligned}
 \E\biggl[\int_0^{\delta}(\sigma^{-1}(r,X^{t,x}_r)&f(r,X^{t,x}_r,u^{\epsilon}_r,v^{\epsilon}_r)-\sigma^{-1}(r,X^{t',x'}_r)f(r,X^{t',x'}_r,u^{\epsilon}_r,v^{\epsilon}_r))^41_{\{|X^{t,x}_r|+|X^{t',x'}_r|\le \rho\}}1_{[r\geq t]}\;dr\biggr]\\ &\leq \E[\int_0^{\delta}(\Phi_\rho(X^{t,x}_r-X^{t',x'}_r))^41_{\{|X^{t,x}_r|+|X^{t',x'}_r|\le \rho\}}1_{[r\geq t]}dr]\rw 0 \mx{ as }\ttrw.
\end{aligned}
\end{equation*}
Thus,
\begin{equation*}
\begin{aligned}
&\E\biggl[\int_0^{\delta}(M^{t',x',u^{\epsilon},v^{\epsilon}}_r)^41_{[r\geq t]}\;dr\biggr]^{\frac{1}{2}}\times \E\biggl[\int_0^{\delta}(\sigma^{-1}(r,X^{t,x}_r)f(r,X^{t,x}_r,u^{\epsilon}_r,v^{\epsilon}_r)\\&\qquad-\sigma^{-1}(r,X^{t',x'}_r)f(r,X^{t',x'}_r,u^{\epsilon}_r,v^{\epsilon}_r))^41_{[r\geq t]}\;dr\biggr]^{\frac{1}{2}}\rw 0 \mx{ as }\ttrw.
\end{aligned}
\end{equation*}
Next by similar arguments we show that the quantity
\begin{equation*}
\E\biggl[\int_0^{\delta}(M^{t',x',u^{\epsilon},v^{\epsilon}}_r)^41_{[t>r\geq t']}\;dr\biggr]^{\frac{1}{2}}\times  \E\biggl[\int_0^{\delta}(\sigma^{-1}(r,X^{t',x'}_r)f(r,X^{t',x'}_r,u^{\epsilon}_r,v^{\epsilon}_r))^41_{[t>r\geq t']}\;dr\biggr]^{\frac{1}{2}},
\end{equation*}
converges likewise to 0 as $(t',x')\rightarrow (t,x)$. Thus
$$\Psi(t',x',\eps,t,x,u^{\epsilon},v^{\epsilon})\rw 0 \mx{ as }\ttrw.$$
And then going back to \eqref{limmu2} to obtain that
\begin{equation*}\label{limmu3}
\E[|M^{t,x,u^{\epsilon},v^{\epsilon}}_{\delta }-M^{t',x',u^{\epsilon},v^{\epsilon}}_{\delta }|^2]\rw 0 \mx{ as }\ttrw.
\end{equation*}
Next from $\eqref{etapeinterm1}$, we obtain
\begin{equation*}
\limsup_{(t',x')\rightarrow (t,x)}\E[|M^{t,x,u^{\epsilon},v^{\epsilon}}_T-M^{t',x',u^{\epsilon},v^{\epsilon}}_T|(1+ \sup\limits_{s\leq T}|X^{t,x}_s|^p)]\leq \frac{C_{t,x,\epsilon}}{\rho^{q_0}}.
\end{equation*}
which implies that, since $\rho$ is arbitrary,
\begin{equation*}
\limsup_{(t',x')\rightarrow (t,x)}\E[|M^{t,x,u^{\epsilon},v^{\epsilon}}_T-M^{t',x',u^{\epsilon},v^{\epsilon}}_T|(1+ \sup\limits_{s\leq T}|X^{t,x}_s|^p)]=0.
\end{equation*}
Finally from $\eqref{term3}$, the term $I_3$ converges to $0$ as $\ttrw$.
It implies that
\begin{equation*}
\limsup_{(t',x')\rightarrow (t,x)}y^{t',x'}_0\leq y^{t,x}_0+2\epsilon.
\end{equation*}
As $\epsilon$ are arbitrary then sending  $\epsilon\rightarrow 0$
\begin{equation*}
\limsup_{(t',x')\rightarrow (t,x)}y^{t',x'}_0\leq y^{t,x}_0.
\end{equation*}
Therefore the function $u$ is upper semi-continuous.
\par Now we show that the function $u$ is lower semi-continuous.
Fix an arbitrary $\epsilon>0$, and we assume that $t'<t$ . Pick $(u^{\epsilon},\tau^{\epsilon})\in \U\times\mathcal{T}_0$  such that
\begin{equation*}
\begin{aligned}
y^{t',x'}_{0}\geq&\,\mathbb{E} \biggl[M^{t',x',u^{\epsilon},v}_T(\int^{\tau^{\epsilon}\wedge\sigma}_{0}\, \Gamma(s,X^{t',x'}_s,u^{\epsilon}_s,v_s)1_{[s\geq t']}\mathrm{d}{s} + h(\sigma\vee s,X^{t',x'}_{\sigma^\vee s})1_{[\sigma\leq\tau^{\epsilon}<T]}\biggr.\\
& \biggl.\quad +h'(\tau^{\epsilon}\vee s,X^{t',x'}_{\tau^{\epsilon}\vee s})1_{[\tau^{\epsilon}<T]}+ g(X^{t',x'}_T)1_{[\sigma=\tau^{\epsilon}=T]})\biggr]- \epsilon.
\end{aligned}
\end{equation*}
where $(v,\sigma)\in \V\times\mathcal{T}_0$  will be chosen later. On the other hand, pick $(v^{\epsilon},\sigma^{\epsilon})\in \V\times\mathcal{T}_0$ such that
\begin{equation*}
\begin{aligned}
y^{t,x}_{0}\leq&\, \mathbb{E} \biggl[M^{t,x,u^{\epsilon},v^{\epsilon}}_T(\int^{\tau^{\epsilon}\wedge\sigma^{\epsilon}}_{0}\, \Gamma(s,X^{t,x}_s,u^{\epsilon}_s,v^{\epsilon}_s)1_{[s\geq t]}\mathrm{d}{s} + h(\sigma^{\epsilon}\vee s,X^{t,x}_{\sigma^{\epsilon}\vee s})1_{[\sigma^{\epsilon}\leq\tau^{\epsilon}<T]}\biggr.\\
& \biggl.\quad +h'(\tau^{\epsilon}\vee s,X^{t,x}_{\tau^{\epsilon}\vee s})1_{[\tau^{\epsilon}<T]}+ g(X^{t,x}_T)1_{[\sigma^{\epsilon}=\tau^{\epsilon} =T]})\biggr]+ \epsilon.
\end{aligned}
\end{equation*}
Therefore
\begin{equation*}
\begin{aligned}
y^{t',x'}_0-y^{t,x}_0 &\geq\mathbb{E}[-M^{t',x',u^{\epsilon},v^{\epsilon}}_T \int_{0}^{\tau^{\epsilon}\wedge\sigma^{\epsilon}}\{| \Gamma(s,X^{t',x'}_s,u^{\epsilon}_s,v^{\epsilon}_s)-\Gamma(s,X^{t,x}_s,u^{\epsilon}_s,v^{\epsilon}_s)|1_{[s\geq
	t]}\\{}&+|\Gamma(s,X^{t',x'}_s,u^{\epsilon}_s,v^{\epsilon}_s)|1_{[t'\leq s< t]}\}ds-M^{t',x',u^{\epsilon},v^{\epsilon}}_T|h(\sigma^{\epsilon}\vee s,X^{t,x}_{\sigma^{\epsilon}\vee s})-h(\sigma^{\epsilon}\vee s,X^{t',x'}_{\sigma^{\epsilon}\vee s})|1_{[\sigma^{\epsilon}\leq\tau^{\epsilon}<T]}\\{}&-M^{t',x',u^{\epsilon},v^{\epsilon}}_T|h'(\tau^{\epsilon}\vee s,X^{t,x}_{\tau^{\epsilon}\vee s})-h'(\tau^{\epsilon}\vee s,X^{t',x'}_{\tau^{\epsilon}\vee s})|1_{[\tau^{\epsilon}<T]} -M^{t',x',u^{\epsilon},v^{\epsilon}}_T|g(X^{t,x}_T)\\&-g(X^{t',x'}_T)|1_{[\sigma^{\epsilon}=\tau^{\epsilon}=T]}-|M^{t,x,u^{\epsilon},v^{\epsilon}}_T-M^{t',x',u^{\epsilon},v^{\epsilon}}_T||(\int_{0}^{\tau^{\epsilon}\wedge\sigma^{\epsilon}}\Gamma(s,X^{t,x}_s,u^{\epsilon}_s,v^{\epsilon}_s)1_{[s\geq t]}ds\\&+h(\sigma^{\epsilon}\vee s,X^{t,x}_{\sigma^{\epsilon}\vee s})1_{[\sigma^{\epsilon}\leq\tau^{\epsilon}<T]}+h'(\tau^{\epsilon}\vee s,X^{t,x}_{\tau^{\epsilon}\vee s})1_{[\tau^{\epsilon}<T]}+g(X^{t,x}_T)1_{[\tau^{\epsilon}=\sigma^{\epsilon}=T]})|]-2\epsilon,\\& \geq -I_1 - I_2 - I_3-2\epsilon,
\end{aligned}
\end{equation*}
where $I_1$, $I_2$ and $I_3$ they are the same as above.\\
By the same argument as above, we have
\begin{equation*}
\liminf_{(t',x')\rightarrow (t,x)}y^{t',x'}_0\geq y^{t,x}_0-2\epsilon.
\end{equation*}
As $\epsilon$ are arbitrary then sending  $\epsilon\rightarrow 0$
\begin{equation*}
\liminf_{(t',x')\rightarrow (t,x)}y^{t',x'}_0\geq y^{t,x}_0.
\end{equation*}
Therefore the function $u$ is lower semi-continuous. We then proved that $u$ is continuous in $(t,x)$.
\end{proof}
We now consider the following double obstacle problem which is parabolic PDE. Roughly speaking, a solution of the double obstacle problem is a function $u:[0,T]\times \mathbb{R}^d\rightarrow \mathbb{R}$ which satisfies:
\begin{equation}\label{dovi}
\begin{cases}
\begin{aligned}
&\min\biggl\{u(t,x)-h(t,x),\max
\biggl[-\frac{\partial u}{\partial t}(t,x)-\mathcal{L}u(t,x)\\ &\;\; \qquad-H^*(t,x,\sigma(t,x)\nabla u(t,x)),u(t,x)-h'(t,x)\biggr]\biggr\}=0 ,\;(t,x)\in[0,T[ \times \mathbb{R} ^d;\\
&u(T,x)=g(x),\qquad x\in\mathbb{R}^d,
\end{aligned}
\end{cases}
\end{equation}
where $\mathcal{L}$ is second order partial differential operator associated with $X^{t,x}_.$ i.e.
$$\mathcal{L}=\frac{1}{2}\sum_{i,j=1}^{d}\left(\left(\sigma\sigma^*\right)\left(t,x\right)\right)_{i,j}\frac{\partial^2}{\partial x_i\partial x_j}.$$
\par First we start by the definition of the viscosity solution for $\eqref{dovi}.$
\begin{definition}
	Let $u$ be a function which belongs to $\mathcal{C}\left([0,T]\times \mathbb{R}^d\right)$. It is said to be a viscosity:
	\begin{itemize}
		\item[$(i)$] subsolution  (resp. supersolution) of $\eqref{dovi}$ if $u(T,x) \leq g(x)$ (resp., $u(T,x) \geq g(x)$) and for any $\phi\in\mathcal{C}^{1,2}([0,T]\times\mathbb{R}^d)$ and any local maximum (resp. minimum) point $(t,x)\in[0,T]\times\mathbb{R}^d$ of $v-\phi,$ we have
		\begin{equation*}
		\begin{aligned}
		&\min\biggl\{u(t,x)-h(t,x),\max
		\biggl[-\frac{\partial u}{\partial t}(t,x)-\mathcal{L}u(t,x)\\ &\;\;\qquad\qquad -H^*(t,x,\sigma(t,x)\nabla u(t,x)),u(t,x)-h'(t,x)\biggr]\biggr\}\leq0\,\,(resp.\;\; \geq0).
		\end{aligned}
		\end{equation*}
		\item[$(ii)$]solution of $\eqref{dovi}$ if it is both a viscosity subsolution and supersolution. $\qquad\qquad\qquad\quad\Box$
	\end{itemize}
\end{definition}
We now give an equivalent definition of a viscosity solution of the parabolic equation $\eqref{dovi}$, (see $\cite{CIL}$ of  definition of the notions of superjet $J^{2,+}$ and subjet $J^{2,-}$ of a continuous function $u$).
\begin{definition}Let $u$ be a  continuous functions defined
	on $[0,T]\times \R^d$, $\R$-valued and such that
	$u(T,x)=g(x)$ for any $x\in \R^d$.  It can be said that $u$ is a viscosity supersolution (resp.
	subsolution) of (\ref{dovi}) if for any, $(t,x)\in
	(0,T)\times \R^d$ and $(p,q,X)\in J^{2,-} u(t,x)$ (resp. $J^{2,+}
	u (t,x)$),
	\begin{equation}\label{dovi2}
	\begin{aligned}
	&\min\biggl\{u(t,x)-h(t,x),\max
	\biggl[-p -\frac{1}{2}Tr[\sigma^{*} X \sigma] \\ &\;\; \qquad\qquad-H^*(t,x,q\sigma(t,x)),u(t,x)-h'(t,x)\biggr]\biggr\}\geq 0 \,\,(resp,\;\; \leq 0).
	\end{aligned}
	\end{equation}
	It is called a viscosity solution if it is both a viscosity subsolution and
	supersolution .$\qquad\qquad\Box$
\end{definition}
\subsection{Existence of a continuous solution of the double obstacle problem}
\begin{proposition}
The function $ u(t,x)=Y^{t,x}_t$  defined in $\eqref{uy}$ is a viscosity solution of the obstacle problem $\eqref{dovi}$.
\end{proposition}
\begin{proof}
Let, $\forall n,m\geq1$
\begin{equation*}
\begin{aligned}
&\overline{H}^{*n,m}(t,x,z)={H^*}^+(t,x,z)(1_{\{|x| \leq n\}}+(x+1+n)1_{\{-n-1\leq x \leq -n\}}+(-x+1+n)1_{\{n\leq x \leq n+1\}}) \\ & \quad\qquad\qquad-{H^*}^-(t,x,z)(1_{\{|x| \leq m\}}+(x+1+m)1_{\{-m-1\leq x \leq -m\}}+(-x+1+m)1_{\{m\leq x \leq m+1\}}),
\end{aligned}
\end{equation*}
$\overline{H}^{*n,m}$ is also a truncation of $H^*$ which is continuous w.r.t. $(t,x,z)$ and it is non-decreasing in $n,$ decreasing in $m.$ Now let us consider the double barrier reflected BSDE associated with the quadruple $(g,\overline{H}^{*n,m},h,h')$
\begin{equation}\label{RBSDExtsnm}
\begin{cases}
\begin{aligned}
\overline{Y}^{t,x,n,m}_s = g(X^{t,x}_T) &+ \int_{s}^{T} \overline{H}^{*n,m}(r,X^{t,x},\overline{Z}^{t,x,n,m}_r) \, \mathrm{d}{r} + (\overline{K}_T^{+,t,x,n,m} - \overline{K}_s^{+,t,x,n,m}) \\ &- ( \overline{K}_T^{-,t,x,n,m} - \overline{K}_s^{-,t,x,n,m}) -  \int_{s}^{T}\overline{Z}^{t,x,n,m}_r \, \mathrm{d}{B_r}, \quad\forall s\leq T;
\end{aligned}\\
\forall s\in[t,T], \quad  h(s,X^{t,x}_s)\leq \overline{Y}^{t,x,n,m}_s\leq h'(s,X^{t,x}_s),\\
\int_{0}^{T} (\overline{Y}^{t,x,n,m}_s - h(s,X^{t,x}_s)) \, \mathrm{d}{\overline{K}_s^{+,t,x,n,m}} = 0 \; and \; \int_{0}^{T} (h'(s,X^{t,x}_s)-\overline{Y}^{t,x,n,m}_s) \, \mathrm{d}{\overline{K}_s^{-,t,x,n,m}} = 0.
\end{cases}
\end{equation}
Since $\overline{H}^{*n,m}$ is Lipschitz w.r.t. $z$, and thanks to $\cite{HM}$, the equation $(\ref{RBSDExtsnm})$ has a solution. Once more due to \cite{HM} we know that
$$\overline{u}_{n,m}(t,x)=\overline{Y}^{t,x,n,m}_t,\qquad (t,x)\in[0,T]\times \R^d,$$
is deterministic and is a viscosity solution of the following parabolic PDE: $\forall(t,x)\in[0,T]\times \R^d$
\begin{equation}\label{dovinm}
\begin{aligned}
&\min\biggl\{\overline{u}_{n,m}(t,x)-h(t,x),\max
\biggl[-\frac{\partial \overline{u}_{n,m}}{\partial t}(t,x)-\mathcal{L}\overline{u}_{n,m}(t,x)\\ &\;\; \qquad\qquad-\overline{H}^{*n,m}(t,x,\sigma(t,x)\nabla \overline{u}_{n,m}(t,x)),\overline{u}_{n,m}(t,x)-h'(t,x)\biggr]\biggr\}=0.
\end{aligned}
\end{equation}
Moreover it satisfies the following estimate $$|\overline{u}_{n,m}(t,x)|\leq C'(1+|x|^p).$$
Now since $\overline{H}^{*n,m}$ is Lipschitz w.r.t. $z$ and increasing w.r.t. $n$ then from the comparison Theorem (see \cite{HM}), $(\overline{Y}^{t,x,n,m})_{n\geq1}\nearrow \overline{Y}^{t,x,m}$. we then obtain:
$$\lim_{n\rightarrow +\infty} \overline{u}_{n,m}(t,x)= \overline{u}_{m}(t,x)=\overline{Y}^{t,x,m}_t,\quad \forall(t,x)\in[0,T]\times \R^d.$$
where $(\overline{Y}^{t,x,m},\overline{Z}^{t,x,m}, \overline{K}^{+,t,x,m},\overline{K}^{-,t,x,m})$ be the solution of the double barrier reflected BSDE associated with $(g,\overline{H}^{*m},h,h')$.\\
We now show that $\overline{u}_{m}$, is the viscosity solution of the following parabolic PDE, for every $(t,x)\in[0,T]\times \R^d$
\begin{equation}\label{dovim}
\begin{aligned}
&\min\biggl\{\overline{u}_{m}(t,x)-h(t,x),\max
\biggl[-\frac{\partial \overline{u}_{m}}{\partial t}(t,x)-\mathcal{L}\overline{u}_{m}(t,x)\\ &\;\; \qquad\qquad-\overline{H}^{*m}(t,x,\sigma(t,x)\nabla \overline{u}_{m}(t,x)),\overline{u}_{m}(t,x)-h'(t,x)\biggr]\biggr\}=0,
\end{aligned}
\end{equation}
and $\overline{u}_{m}(T,x)=g(x).$\\
We can easily show the continuity of $\overline{u}_{m}$. Indeed,  let $\epsilon> 0$ and $(t',x') \in B((t,x),\epsilon)$ let $t'\in[t,T]$, we have:
	\begin{eqnarray*}
		\overline{u}_{m}(t',x')-\overline{u}_{m}(t,x) &= & \overline{u}_{m}(t',x')-\overline{u}_{n,m}(t',x')\\& &
+\overline{u}_{n,m}(t',x')-\overline{u}_{n,m}(t,x)
		\\& & +\overline{u}_{n,m}(t,x)-\overline{u}_{m}(t,x).
	\end{eqnarray*}
	As $\overline{u}_{n,m}\nearrow \overline{u}_m$, then we have $\overline{u}_{n,m}\leq \overline{u}_m$.
	Therefore
\begin{equation}
		\overline{u}_{m}(t',x')-\overline{u}_{m}(t,x) \leq  \overline{u}_{m}(t',x')-\overline{u}_{n,m}(t',x')
+\overline{u}_{n,m}(t',x')-\overline{u}_{n,m}(t,x).\label{unm}
	\end{equation}
As $n$ is arbitrary then putting $n\rightarrow +\infty$ and using $\lim\limits_{n\rightarrow +\infty}\overline{u}_{n,m}(t,x)=\overline{u}_{m}(t,x)$ and
the continuity of $\overline{u}_{n,m}(t,x)$ in $t$ and $x$, we get that the right hand side terms of (\ref{unm})
converge to 0 as $(t', x') \rightarrow (t, x)$. Then we obtain:
	$$\limsup_{(t',x')\to (t,x)}\overline{u}_{m}(t',x')\leq\overline{u}_{m}(t,x).$$
	Therefore $\overline{u}_{m}$ is upper semi-continuous. But $\overline{u}_{m}$ is also lower semi-continuous, therefore it is continuous.

Now since $\overline{u}_{n,m}$ and $\overline{u}_{m}$ are continuous, it follows from Dini's Theorem that the above convergence is uniform on compact. \\
Let us show that $\overline{u}_{m}$ is a viscosity subsolution  of $(\ref{dovim})$. Since $\overline{u}_{m}(T,x)=g(x)$ and $h(t,x)\leq \overline{u}_{m}(t,x)\leq h'(t,x)$, it is sufficient to prove that for any $(t,x)\in
(0,T)\times \R^d$ and $(p,q,X)\in J^{2,+}\overline{u}_{m}(t,x)$, such that $\overline{u}_{m}(t,x)\geq h(t,x),$ we have
$$-p -\frac{1}{2}Tr[\sigma^{*} X \sigma]-\overline{H}^{*m}(t,x,q\sigma(t,x))\leq 0.$$
From Lemma $6.1$ in $\cite{CIL},$ there exists sequences:
$$
\begin{array}{l}
n_j \rightarrow+\infty,\quad (t_j,x_j)\rightarrow (t,x), \quad (p_j,q_j,X_j)\in
J^{2,+}\overline{u}_{n_j,m}(t_j,x_j),
\end{array}
$$
such that $$(p_j,q_j,X_j)\rightarrow (p,q,X).$$ But for any $j$,
$$
\begin{array}{l}
-p_j -\frac{1}{2}Tr[\sigma^{*}(t_j,x_j) X_j \sigma(t_j,x_j)] -\overline{H}^{*n_j,m}(t_j,x_j,q_j\sigma(t_j,x_j))\leq 0.
\end{array}
$$
From the assumption that $\overline{u}_{m}(t,x)\geq h(t,x)$ and the uniform convergence of
$\overline{u}_{n,m},$ it follows that for $j$ large enough $\overline{u}_{n_j,m}(t_j,x_j)\geq h(t_j,x_j)$.
Hence, taking the limit as $j\rightarrow +\infty$ in the above
inequality yields:
$$-p -\frac{1}{2}Tr[\sigma^{*} X \sigma] -\overline{H}^{*m}(t,x,q\sigma(t,x))\leq 0.$$
Thus that $\overline{u}_{m}$ is a subsolution of (\ref{dovim}).
\par We now show that $\overline{u}_{m}$ is a supersolution of
$(\ref{dovim}).$ Let $(t,x)$ be arbitrary in $[0,T)\times \R^d$, and
$(p,q,X)\in J^{2,-}\overline{u}_{m}(t,x).$ We already know that $\overline{u}_{m}(t,x)\geq h(t,x).$ By the
same argument as above, there exist sequences:
$$
\begin{array}{l}
n_j \rightarrow+\infty,\quad (t_j,x_j)\rightarrow (t,x),\quad (p_j,q_j,X_j)\in
J^{2,-}\overline{u}_{n_j,m}(t_j,x_j),
\end{array}
$$
such that $$(p_j,q_j,X_j)\rightarrow (p,q,X).$$ But for any $j$,
$$
\begin{array}{l}
-p_j -\frac{1}{2}Tr[\sigma^{*}(t_j,x_j) X_j \sigma(t_j,x_j)] -\overline{H}^{*n_j,m}(t_j,x_j,q_j\sigma(t_j,x_j))\geq 0.
\end{array}
$$
Hence, taking the limit as $j\rightarrow +\infty$, we conclude that:
$$-p -\frac{1}{2}Tr[\sigma^{*} X \sigma] -\overline{H}^{*m}(t,x,q\sigma(t,x))\geq 0. $$
And we have proved that $\overline{u}_{m}$ is a supersolution of $(\ref{dovim}).$
\par Next we have
$$|\overline{Y}^{t,x,m}_t|\leq C'(1+|x|^p).$$
From the previous results we have
$$\lim_{m\rightarrow +\infty} \overline{u}_{m}(t,x)= u(t,x)=Y^{t,x}_t,\quad \forall(t,x)\in[0,T]\times \R^d,$$
since $\overline{Y}^{t,x,m}\searrow Y^{t,x}$ as $m\rightarrow \infty$. Since $\overline{u}_{m}$ and $u$ are continuous, it follows from Dini's Theorem that the above convergence
is uniform on compact. We now show that $u$ is a subsolution  of (\ref{dovi}). Let $(t,x)$ be a point at which
$u(t, x) \geq h(t,x)$, and let $(p,q,X)\in J^{2,+} u(t,x)$. From Lemma 6.1 in \cite{CIL}, there exists sequences:
$$
\begin{array}{l}
m_j \rightarrow+\infty,\quad (t_j,x_j)\rightarrow (t,x), \quad (p_j,q_j,X_j)\in
J^{2,+}\overline{u}_{m_j}(t_j,x_j),
\end{array}
$$
such that $$(p_j,q_j,X_j)\rightarrow (p,q,X).$$ But for any $j$,
$$
\begin{array}{l}
-p_j -\frac{1}{2}Tr[\sigma^{*}(t_j,x_j) X_j \sigma(t_j,x_j)] -\overline{H}^{*m_j}(t_j,x_j,q_j\sigma(t_j,x_j))\leq 0.
\end{array}
$$
From the assumption that $u(t,x)\geq h(t,x)$ and the uniform convergence of
$\overline{u}_{m},$ it follows that for $j$ large enough $\overline{u}_{m_j}(t_j,x_j)\geq h(t_j,x_j)$.
Hence, taking the limit as $j\rightarrow +\infty$ in the above
inequality yields:
$$-p -\frac{1}{2}Tr[\sigma^{*} X \sigma] -H^{*}(t,x,q\sigma(t,x))\leq 0. $$
And we have proved that $u$ is a subsolution of (\ref{dovi}).\par
We now show that $u$ is a supersolution of
(\ref{dovi}). Let $(t,x)$ be arbitrary in $[0,T]\times \R^d$, and
$(p,q,X)\in J^{2,-}u(t,x).$ We already know that $u(t,x)\geq h(t,x).$ By the
same argument as above, there exist sequences:
$$
\begin{array}{l}
m_j \rightarrow+\infty,\quad (t_j,x_j)\rightarrow (t,x),\quad (p_j,q_j,X_j)\in
J^{2,-}\overline{u}_{m_j}(t_j,x_j),
\end{array}
$$
such that $$(p_j,q_j,X_j)\rightarrow (p,q,X).$$ But for any $j$,
$$
\begin{array}{l}
-p_j -\frac{1}{2}Tr[\sigma^{*}(t_j,x_j) X_j \sigma(t_j,x_j)] -\overline{H}^{*m_j}(t_j,x_j,q_j\sigma(t_j,x_j))\geq 0.
\end{array}
$$
Hence, taking the limit as $j\rightarrow +\infty$, we conclude that:
$$-p -\frac{1}{2}Tr[\sigma^{*} X \sigma] -H^*(t,x,q\sigma(t,x))\geq 0. $$
And we have proved that $u$ is a supersolution of (\ref{dovi}).
\end{proof}
\subsection{Uniqueness}
 In this section we are going to show the uniqueness of the viscosity solution of the system . We first need the following lemma.
\begin{lemma}\label{lemma-supersol}
	Let $v$ be a supersolution of the system $(\ref{dovi}),$ then for any $\gamma\geq 0$ there exists $\alpha>0$ such that for any $\lambda \geq \alpha$ and $\theta >0$,  $(v(t,x)+\theta e^{-\lambda t}\mid x\mid^{2\gamma+2})$ is a supersolution for $(\ref{dovi})$.
\end{lemma}
\begin{proof}
	 We assume w.l.o.g. that the function
	$v(t,x)$ is lsc. Let  $\varphi \in \mathcal{C}^{1,2}$ be such that the
	function $\varphi-(v+\theta e^{-\lambda t}\mid
	x\mid^{2\gamma+2})$ has a local maximum in $(t,x)$ which is equal to
	0. Since $v(t,x)$ is a supersolution for
	$(\ref{dovi}),$ then we have:
	\begin{equation*}
	\begin{aligned}
	&\min\biggl\{v(t,x)-h(t,x),\max
	\biggl[-\partial_{t}\Big(\varphi(t,x)-\theta e^{-\lambda t}\mid x\mid^{2\gamma+2}\Big) -\frac{1}{2}Tr\Big[\sigma.\sigma^{*}(t,x)D^{2}_{x}\Big(\varphi(t,x)\\ &\qquad-\theta e^{-\lambda t}\mid x\mid^{2\gamma+2}\Big)\Big]-H^*(t,x,\nabla (\varphi(t,x)-\theta e^{-\lambda t}\mid x\mid^{2\gamma+2})\sigma(t,x)),v(t,x)-h'(t,x)\biggr]\biggr\}\geq0.
	\end{aligned}
	\end{equation*}
	Hence
	\begin{eqnarray}\label{supersol}
	(v(t,x)+\theta e^{-\lambda t}\mid
	x\mid^{2\gamma+2})-h(t,x) \geq v(t,x)-h(t,x) \geq 0.
	\end{eqnarray}
	On the other hand:
	\begin{equation*}
	\begin{array}{l}
	\left.-\partial_{t}\Big(\varphi(t,x)-\theta e^{-\lambda t}\mid x\mid^{2\gamma+2}\Big)-\frac{1}{2}Tr\Big[\sigma.\sigma^{*}(t,x)D^{2}_{x}\Big(\varphi(t,x)-\theta e^{-\lambda t}\mid x\mid^{2\gamma+2}\Big)\Big] \right.\\\qquad
	\left.-H^*(t,x,\nabla (\varphi(t,x)-\theta
	e^{-\lambda t}\mid x\mid^{2\gamma+2})\sigma(t,x))\right. \geq 0.
	\end{array}
	\end{equation*}
	Thus
	\begin{equation*}
	\begin{aligned}
	-\partial_{t}\varphi(t,x)&-\frac{1}{2}Tr\Big[\sigma.\sigma^{*}(t,x)D^{2}_{x}\varphi(t,x)\Big]-H^*(t,x,\nabla \varphi(t,x)\sigma(t,x)) \geq \theta \lambda e^{-\lambda t}\mid
	x\mid^{2\gamma+2}\\ &-\frac{1}{2}\theta e^{-\lambda t}
	Tr\Big[\sigma.\sigma^{*}(t,x)D^{2}_{x}\mid x\mid^{2\gamma+2}\Big]
	+H^*(t,x,\nabla (\varphi(t,x)-\theta
	e^{-\lambda t}\mid x\mid^{2\gamma+2})\sigma(t,x))\\ &- H^*(t,x,\nabla \varphi(t,x)\sigma(t,x)).
	\end{aligned}
	\end{equation*}
	From the definition of $H^*$ we have,
	\begin{equation}\label{supersol1}
	\begin{aligned}
	-\partial_{t}&\varphi(t,x)-\frac{1}{2}Tr\Big[\sigma.\sigma^{*}(t,x)D^{2}_{x}\varphi(t,x)\Big]-H^*(t,x,\nabla \varphi(t,x)\sigma(t,x))
	\geq \theta \lambda e^{-\lambda t}\mid
	x\mid^{2\gamma+2}\\ &-\frac{1}{2}\theta e^{-\lambda t}
	Tr\Big[\sigma.\sigma^{*}(t,x)D^{2}_{x}\mid x\mid^{2\gamma+2}\Big]
	-\inf\limits_{u\in A}\sup\limits_{v\in B}\Big[\nabla (\theta
	e^{-\lambda t}\mid x\mid^{2\gamma+2})f(t,x,u,v)\Big].
	\end{aligned}
    \end{equation}
	Therefore taking into account the growth conditions on $f$ and $\sigma$ and setting $\theta>0$, there exists two positive constants $C_{1}$ and $C_{2}$ such that:\\
	$\frac{1}{2}\theta e^{-\lambda t} Tr[\sigma.\sigma^{*}(t,x)D^{2}_{x}\mid x\mid^{2\gamma+2}]\leq C_{1}\theta e^{-\lambda t} \mid x\mid^{2\gamma+2}$ and $\nabla (\theta
	e^{-\lambda t}\mid x\mid^{2\gamma+2})f(t,x,u,v)\leq C_{2}\theta e^{-\lambda t}\mid x\mid^{2\gamma+2}$. \\
	Then by $(\ref{supersol1})$ we get
	\begin{equation}\label{supersol2}
	\left.-\partial_{t}\varphi(t,x)-\frac{1}{2}Tr\Big[\sigma.\sigma^{*}(t,x)D^{2}_{x}\varphi(t,x)\Big]-H^*(t,x,\nabla \varphi(t,x)\sigma(t,x)) \right.\\
	\left.\geq \theta (\lambda-\alpha)  e^{-\lambda t}\mid
	x\mid^{2\gamma+2},\right.
	\end{equation}
	where $\alpha=C_{1}+C_{2}$.\\
	Now since $\theta>0$, then we conclude that for $\lambda \geq \alpha$ the right hand side of $(\ref{supersol2})$ is non-negative.\\
	Finally, $(v(t,x)+\theta e^{-\lambda
		t}\mid x\mid^{2\gamma+2})$ is a viscosity
	supersolution for $(\ref{dovi}).$
\end{proof}
	Now we give an equivalent of quasi-variational inequality
	(\ref{dovi}). In this section, we consider the new function
	$\overline{v}$ given by the classical change of variable $\overline{v}(t,x) =
	\exp(t)v(t, x)$, for any $t\in[0,T ]$ and $x\in \R^d$. Of course,
	the function $\overline{v}$ is continuous and of polynomial growth with respect to its
	second argument.\\ A second property is given by the following proposition.
	\begin{proposition}
		$v$ is a viscosity solution of $(\ref{dovi})$ if and only if
		$\overline{v}$ is a viscosity solution to the following quasi-variational
		inequality in $[0,T[\times \R^d$,
		\begin{equation}\label{dovi1}
		\begin{aligned}
		&\min\biggl\{\overline{v}(t,x)-\exp(t)h(t,x),\max
		\biggl[-\frac{\partial \overline{v}}{\partial t}(t,x)+\overline{v}(t,x)-\mathcal{L}\overline{v}(t,x)\\ &\qquad\qquad-\exp(t)H^*(t,x,\sigma(t,x)\nabla (\exp(-t)\overline{v}(t,x)),\overline{v}(t,x)-\exp(t)h'(t,x)\biggr]\biggr\}=0.
        \end{aligned}
		\end{equation}
		The terminal condition for $\overline{v}$ is: $\overline{v}(T,x)=\exp(T)g(x)$ in
		$\R^d.$ \qed
	\end{proposition}
	We are going now to address the
	question of uniqueness of the viscosity solution of
	quasi-variational inequality $(\ref{dovi1}).$ We have the following theorem:
	\begin{theorem}\label{uni}
		The solution in viscosity sense of quasi-variational inequality $(\ref{dovi1})$ is unique in the space of continuous functions on $[0,T]\times \R^d$ which satisfy a polynomial growth condition, i.e., in the space
		\begin{equation*}
		\begin{aligned}
		\cal C:=&\{\varphi: [0,T]\times \R^d\rightarrow \R, \mbox{ continuous and for any } (t,x), \\ &\qquad
		|\varphi(t,x)|\leq C(1+|x|^\gamma) \mbox{ for some constants } C
		\mbox{ and }\gamma\}.
		\end{aligned}
	   \end{equation*}
	\end{theorem}
\begin{proof}
	We will show by contradiction that if $u$ and $v$ is a subsolution and a supersolution respectively for $(\ref{dovi})$ then  $u\leq v$.
	Therefore if we have two solutions of $(\ref{dovi1})$ then they are obviously equal. Next according to Lemma \ref{lemma-supersol}, it is enough to show that $\forall (t,x) \in [0,T]\times \mathbb{R}^{d}$,
	$$ u(t,x)\leq v(t,x)+\kappa e^{(1-\lambda) t}\mid x\mid^{2\gamma+2} \ ,$$
	since in taking the limit as $\kappa\rightarrow 0$ we get the
	desired result. \\ So let us set  $w(t,x)=v(t,x)+\kappa
	e^{(1-\lambda) t}\mid x\mid^{2\gamma+2}$, $(t,x) \in [0,T]\times \mathbb{R}^{d}$.
	Then, there exists $R > 0$ such that
	\begin{equation}\label{comp_uni}
		\begin{aligned}
			\max_{(t,x) \in [0,T]\times
				\mathbb{R}^{d}}(u(t,x)-w(t,x))&=
			\max_{(t,x) \in [0,T[\times
				B_R}(u(t,x)-w(t,x))\\ &=u(\overline{t},\overline{x})-w(\overline{t},\overline{x})=\eta>0,
		\end{aligned}
	\end{equation}
	where $B_R := \{x\in \mathbb{R}^{d}; |x|<R\}$ and $(\overline{t},\overline{x}) \in [0,T[\times B_R$.\\
	Now let us take
	$\theta$ and $\beta \in ]0,1]$ small enough. Here $\gamma$ is the growth exponent of the functions
	which w.l.o.g. we assume integer and $\geq 2$. Then, for a small
	$\epsilon>0$, let us define:
	
	\begin{equation}
		\label{phi}
		\Phi_{\epsilon}(t,x,y)=u(t,x)-w(t,y)-\frac{1}{2\epsilon}|x-y|^{2\gamma}
		-\theta( |x-\overline{x}|^{2\gamma+2}+|y-\overline{x}|^{2\gamma+2})-\beta
		(t-\overline{t})^2.
	\end{equation}
	By the growth assumption on $u$ and $w$, there exists a
	$(t_{\epsilon},x_{\epsilon},y_{\epsilon})\in [0,T]\times B_R
	\times B_R $, for $R$ large enough, such that:
	$$\Phi_{\epsilon}(t_{\epsilon},x_{\epsilon},y_{\epsilon})=\max\limits_{(t,x,y)\in [0,T]\times \overline{B}_R\times
		\overline{B}_R}\Phi_{\epsilon}(t,x,y).$$
	On the other hand, from
	$2\Phi_{\epsilon}(t_{\epsilon},x_{\epsilon},y_{\epsilon})\geq
	\Phi_{\epsilon}(t_{\epsilon},x_{\epsilon},x_{\epsilon})+\Phi_{\epsilon}(t_{\epsilon},y_{\epsilon},y_{\epsilon})$,
	we have
	\begin{equation}
		\frac{1}{\epsilon}|x_{\epsilon} -y_{\epsilon}|^{2\gamma} \leq
		(u(t_{\epsilon},x_{\epsilon})-u(t_{\epsilon},y_{\epsilon}))+(w(t_{\epsilon},x_{\epsilon})-w(t_{\epsilon},y_{\epsilon})),
	\end{equation}
	and consequently $\frac{1}{\epsilon}|x_{\epsilon}
	-y_{\epsilon}|^{2\gamma}$ is bounded, and as $\epsilon\rightarrow 0$,
	$|x_{\epsilon} -y_{\epsilon}|\rightarrow 0$. Since $u$ and $w$ are
	uniformly continuous on $[0,T]\times \overline{B}_R$, then
	$\frac{1}{2\epsilon}|x_{\epsilon} -y_{\epsilon}|^{2\gamma}\rightarrow 0$
	as
	$\epsilon\rightarrow 0.$\\
	Since
	$$u(\overline{t},\overline{x})-w(\overline{t},\overline{x}) \leq
	\Phi_{\epsilon}(t_{\epsilon},x_{\epsilon},y_{\epsilon})\leq u(t_{\epsilon},x_\epsilon)-w(t_{\epsilon},y_\epsilon),$$
	from the continuity of $u$ and
	$w$, we deduce that,
	\begin{equation}\label{subsequence}
		(t_\epsilon,x_\epsilon,y_\epsilon)\rightarrow (\overline{t},\overline{x},\overline{x}).
	\end{equation}
Now we claim that
\begin{equation}\label{ineg}u(t_\epsilon,x_\epsilon)-\exp(t_\epsilon)h(t_\epsilon,x_\epsilon)>0.
\end{equation}If not, there exist a subsequence such that $   u(t_\epsilon,x_\epsilon)-\exp(t_\epsilon)h(t_\epsilon,x_\epsilon)\leq 0$, then  we have, $u(\overline{t},\overline{x})-\exp(\overline{t})h(\overline{t},\overline{x})\leq 0$ but from the assumption $u(\overline{t},\overline{x})-v(\overline{t},\overline{x})>0$, we deduce that $0\geq u(\overline{t},\overline{x})-\exp(\overline{t})h(\overline{t},\overline{x})>w(\overline{t},\overline{x})-\exp(\overline{t})h(\overline{t},\overline{x})$. Therefore we have $w(\overline{t},\overline{x})-\exp(\overline{t})h(\overline{t},\overline{x})<0$, which leads to a contradiction with \eqref{dovi1}.
Next let us denote
\begin{equation}
	\varphi_{\epsilon}(t,x,y)=\frac{1}{2\epsilon}|x-y|^{2\gamma}
	+\theta( |x-\overline{x}|^{2\gamma+2}+|y-\overline{x}|^{2\gamma+2})+\beta
	(t-\overline{t})^2.
\end{equation}
Then we have:
\begin{equation}\label{derive}
	\begin{cases}
		\begin{aligned}
			&D_{t}\varphi_{\epsilon}(t,x,y)=2\beta(t-\overline{t}),\\ &
			D_{x}\varphi_{\epsilon}(t,x,y)= \frac{\gamma}{\epsilon}(x-y)|x-y|^{2\gamma-2} +\theta({2\gamma}+2)
			(x-\overline{x})|x-\overline{x}|^{2\gamma}, \\ &
			D_{y}\varphi_{\epsilon}(t,x,y)= -\frac{\gamma}{\epsilon}(x-y)|x-y|^{{2\gamma}-2} +
			\theta({2\gamma}+2)(y-\overline{x})|y-\overline{x}|^{2\gamma},\\ &
			B(t,x,y)=D_{x,y}^{2}\varphi_{\epsilon}(t,x,y)=\frac{1}{\epsilon}
			\begin{pmatrix}
				a_1(x,y)&-a_1(x,y) \\
				-a_1(x,y)& a_1(x,y)
			\end{pmatrix}+ \begin{pmatrix}
				a_2(x)&0 \\
				0&a_2(y)
			\end{pmatrix}, \\ &
			\mbox{ with } a_1(x,y)=\gamma|x-y|^{{2\gamma}-2}I+{\gamma}({2\gamma}-2)(x-y)(x-y)^*|x-y|^{{2\gamma}-2}  \\ & \mbox{ and }
			\, a_2(x)=\theta({2\gamma}+2)|x-\overline{x}|^{2\gamma}I+2\theta\gamma({2\gamma}+2)
			(x-\overline{x})(x-\overline{x})^*|x-\overline{x}|^{2\gamma-2} .
		\end{aligned}
	\end{cases}
\end{equation}
Taking into account $(\ref{ineg})$ then applying the result by Crandall et al. (Theorem $8.3,$ $\cite{CIL}$) to the function
\begin{equation*}
u(t,x)-w(t,x)-\varphi_{\epsilon}(t,x,y)
\end{equation*}
at
the point $(t_\epsilon,x_\epsilon,y_\epsilon)$, for any $\epsilon_1 >0$, we can find
$c,c_1 \in \R$, $q_1,q_2 \in \R^d$ and $X,Y \in S_d$ (the set of symmetric matrices), such that:
\begin{equation}\label{lemmeishii}
\begin{cases}
\begin{aligned}
&(c,q_1,X)
\in J^{2,+}(u(t_\epsilon,x_\epsilon)),\\ &
(-c_1,q_2,Y)\in J^{2,-}
(w(t_\epsilon,y_\epsilon)),\\ &
q_1=\frac{{\gamma}}{\epsilon}(x_\epsilon-y_\epsilon)|x_\epsilon-y_\epsilon|^{2\gamma-2} +\theta({2\gamma}+2)
(x_\epsilon-\overline{x})|x_\epsilon-\overline{x}|^{2\gamma},\\ &
q_2= \frac{{\gamma}}{\epsilon}(x_\epsilon-y_\epsilon)|x_\epsilon-y_\epsilon|^{2\gamma-2} -
\theta({2\gamma}+2)(y_\epsilon-\overline{x})|y_\epsilon-\overline{x}|^{{2\gamma}},\\ &
c+c_1=D_{t}\varphi_{\epsilon}(t_\epsilon,x_\epsilon,y_\epsilon)=2\beta(t_\epsilon-\overline{t}) \mbox{ and finally }\\ &
-(\frac{1}{\epsilon_1}+||B(t_\epsilon,x_\epsilon,y_\epsilon)||)I\leq
\begin{pmatrix}
X&0 \\
0&-Y
\end{pmatrix}\leq B(t_\epsilon,x_\epsilon,y_\epsilon)+\epsilon_1 B(t_\epsilon,x_\epsilon,y_\epsilon)^2.
\end{aligned}
\end{cases}
\end{equation}
Taking now into account (\ref{ineg}), and the
definition of viscosity solution, we get:
\begin{equation*}
-c-\frac{1}{2}Tr[\sigma^{*}(t_\epsilon,x_\epsilon)X\sigma(t_\epsilon,x_\epsilon)]
+u(t_\epsilon,x_\epsilon)-\exp(t_\epsilon)H^*(t_\epsilon,x_\epsilon,\exp(-t_\epsilon)q_1\sigma(t_\epsilon,x_\epsilon))\leq 0
\end{equation*}
and
\begin{equation*}
c_1-\frac{1}{2}Tr[\sigma^{*}(t_\epsilon,y_\epsilon)Y\sigma(t_\epsilon,y_\epsilon)]
+ w(t_\epsilon,y_\epsilon)-\exp(t_\epsilon)H^*(t_\epsilon,y_\epsilon,\exp(-t_\epsilon)q_2\sigma(t_\epsilon,y_\epsilon))\geq
0.
\end{equation*}
 Then
\begin{equation}\label{viscder}
\begin{aligned}
u(t_\epsilon,x_\epsilon)-w(t_\epsilon,y_\epsilon)-c-c_1\leq& \frac{1}{2}Tr[\sigma^{*}(t_\epsilon,x_\epsilon)X\sigma(t_\epsilon,x_\epsilon)-\sigma^{*}(t_\epsilon,y_\epsilon)Y\sigma(t_\epsilon,y_\epsilon)] \\ &\; +\exp(t_\epsilon)H^*(t_\epsilon,x_\epsilon,\exp(-t_\epsilon)q_1\sigma(t_\epsilon,x_\epsilon))\\ &\;-\exp(t_\epsilon)H^*(t_\epsilon,
y_\epsilon,\exp(-t_\epsilon)q_2\sigma(t_\epsilon,y_\epsilon)).
\end{aligned}
\end{equation}
But from (\ref{derive}) there exist two constants $C$ and $C_1$ such
that:
$$||a_1(x_\epsilon,y_\epsilon)||\leq C|x_\epsilon- y_\epsilon|^{2\gamma-2} \mbox{ and }(||a_2(x_\epsilon)||\vee ||a_2(y_\epsilon)||)\leq C_1 \theta.$$
As
$$B= B(t_\epsilon,x_\epsilon,y_\epsilon)= \frac{1}{\epsilon}
\begin{pmatrix}
a_1(x_\epsilon,y_\epsilon)&-a_1(x_\epsilon,y_\epsilon) \\
-a_1(x_\epsilon,y_\epsilon)&a_1(x_\epsilon,y_\epsilon)
\end{pmatrix}+ \begin{pmatrix}
a_2(x_\epsilon)&0 \\
0&a_2(y_\epsilon)
\end{pmatrix},$$
then
$$B\leq \frac{C}{\epsilon}|x_\epsilon - y_\epsilon|^{2\gamma-2}
\begin{pmatrix}
I&-I \\
-I&I
\end{pmatrix}+ C_1 \theta I.$$
It follows that:
\begin{equation}
B+\epsilon_1 B^2 \leq C(\frac{1}{\epsilon}|x_\epsilon - y_\epsilon|^{2\gamma-2}+
\frac{\epsilon_1}{\epsilon^2}|x_\epsilon - y_\epsilon|^{4\gamma-4})\begin{pmatrix}
I&-I \\
-I&I
\end{pmatrix}+ C_1\theta I,
\end{equation}
where $C$ and $C_1$ which hereafter may change from line to line.
Choosing now $\epsilon_1=\epsilon$, yields the relation:
\begin{equation}
\label{ineg_matreciel}
B+\epsilon_1 B^2 \leq \frac{C}{\epsilon}(|x_\epsilon - y_\epsilon|^{2\gamma-2}+|x_\epsilon - y_\epsilon|^{4\gamma-4})\begin{pmatrix}
I&-I \\
-I&I
\end{pmatrix}+ C_1\theta I.
\end{equation}
Now, from the Lipschitz continuity of $\sigma$, (\ref{lemmeishii}) and
(\ref{ineg_matreciel}) we get:
$$\frac{1}{2}Tr[\sigma^{*}(t_\epsilon,x_\epsilon)X\sigma(t_\epsilon,x_\epsilon)-\sigma^{*}(t_\epsilon,y_\epsilon)
Y\sigma(t_\epsilon,y_\epsilon)]\leq \frac{C}{\epsilon}(|x_\epsilon - y_\epsilon|^{2\gamma}+|x_\epsilon - y_\epsilon|^{4\gamma-2}) +C_1 \theta.$$ Next from the Lipschitz condition on $f$ we get:
\begin{equation}
\begin{array}{llll}
\label{fin}&\exp(t_\epsilon)H^*(t_\epsilon,x_\epsilon,\exp(-t_\epsilon)q_1\sigma(t_\epsilon,x_\epsilon))-\exp(t_\epsilon)H^*(t_\epsilon,y_\epsilon,\exp(-t_\epsilon)q_2\sigma(t_\epsilon,y_\epsilon)),\\
& \leq  \inf\limits_{u\in A}\sup\limits_{v\in B}[\exp(t_\epsilon)(H(t_\epsilon,x_\epsilon,\exp(-t_\epsilon)q_1\sigma(t_\epsilon,x_\epsilon),u,v)-H(t_\epsilon,y_\epsilon,\exp(-t_\epsilon)q_2\sigma(t_\epsilon,y_\epsilon),u,v))],\\

& \leq
\inf\limits_{u\in A}\sup\limits_{v\in B}[|q_1f(t_\epsilon,x_\epsilon,u,v)-q_2f(t_\epsilon,y_\epsilon,u,v)\mid+\exp(t_\epsilon)\mid\Gamma(t_\epsilon,x_\epsilon,u,v)
-\Gamma(t_\epsilon,y_\epsilon,u,v)\mid],\\
& \leq \inf\limits_{u\in A}\sup\limits_{v\in B}[|f(t_\epsilon,x_\epsilon,u,v)||q_1-q_2|+|q_2||f(t_\epsilon,x_\epsilon,u,v)-f(t_\epsilon,y_\epsilon,u,v)|\\ & \quad +\exp(t_\epsilon)\mid\Gamma(t_\epsilon,x_\epsilon,u,v)
-\Gamma(t_\epsilon,y_\epsilon,u,v)|],\\
& \leq C(1+|x_\epsilon|)\big|\theta({2\gamma}+2)
(x_\epsilon-\overline{x})|x_\epsilon-\overline{x}|^{2\gamma}+\theta({2\gamma}+2)(y_\epsilon-\overline{x})|y_\epsilon-\overline{x}|^{{2\gamma}}\big|\\ & \quad +
|q_2||x_\epsilon-y_\epsilon| + \exp(t_\epsilon)\inf\limits_{u\in A}\sup\limits_{v\in B}[\mid\Gamma(t_\epsilon,x_\epsilon,u,v)
-\Gamma(t_\epsilon,y_\epsilon,u,v)|].
\end{array}
\end{equation}
By plugging into (\ref{viscder})  we
obtain:
\begin{equation*}
\begin{aligned}
u(t_\epsilon,x_\epsilon)-&w(t_\epsilon,y_\epsilon)-2\beta(t_\epsilon-\overline{t}) \leq \frac{C}{\epsilon}(|x_\epsilon - y_\epsilon|^{2\gamma}+|x_\epsilon - y_\epsilon|^{4\gamma-2}) +C_1 \theta\\ &
+ C(1+|x_\epsilon|)\big|\theta({2\gamma}+2)
(x_\epsilon-\overline{x})|x_\epsilon-\overline{x}|^{2\gamma}+\theta({2\gamma}+2)(y_\epsilon-\overline{x})|y_\epsilon-\overline{x}|^{{2\gamma}}\big|\\ &+
|q_2||x_\epsilon-y_\epsilon|
+ \inf\limits_{u\in A}\sup\limits_{v\in B}[\mid\Gamma(t_\epsilon,x_\epsilon,u,v)
-\Gamma(t_\epsilon,y_\epsilon,u,v)|].
\end{aligned}
\end{equation*}
By sending $\epsilon\rightarrow0$,  $\theta
\rightarrow0$ and taking into account of the uniform continuity  of $\Gamma$, we obtain $\eta < 0$ which is
a contradiction. The proof of Theorem $\ref{uni}$ is now complete.
\end{proof}

\section*{ Appendix A: Proof of Proposition \ref{prop2.1}}

 For $n\geq0$ let us set
 \begin{equation}
\varphi^n(t,x,z)=[C(1+||x||_t)\wedge n]|z|+C(1+||x||^{p}_{t}), \qquad \forall (t,x,z)\in [0,T]\times\Omega\times\mathbb{R}^{d}.
\end{equation}
Then $\varphi^n$ is Lipschitz with respect to $z$. Therefore by Hamadène and Hassani's result in $\cite{HM}$, there exists a quadruple of processes $(Y^n,Z^n,K^{+,n},K^{-,n})\in\mathcal{S}^{2}\times\mathcal{H}^{2,d}\times\mathcal{S}_{ci}^2\times\mathcal{S}_{ci}^2$ ( $K^{\pm,n}$  is non-decreasing and $K^{\pm,n}_0=0$) such that,
\begin{equation}\label{RBSDEn}
\begin{cases}
\begin{aligned}
Y^n_t = g(x) + \int_{t}^{T}& \varphi^n(s,x,Z^n_s) \, \mathrm{d}{s} + (K_T^{+,n} - K_t^{+,n}) - ( K_T^{-,n} - K_t^{-,n})\\ & -  \int_{t}^{T}Z^n_s \, \mathrm{d}{B_s};
\end{aligned}\\
\forall t \leq T, \quad L_t \leq Y^n_t \leq U_t,\\
\int_{0}^{T} (Y^n_s - L_s) \, \mathrm{d}{K_s^{+,n}} = 0 \; and \; \int_{0}^{T} (U_s-Y^n_s) \, \mathrm{d}{K_s^{-,n}} = 0.
\end{cases}
\end{equation}
As $\varphi^n\leq \varphi^{n+1}$ and the barriers $L$ and $U$ are fixed, then by the comparison Theorem (see, $\cite{HM}$) we have, for any $n\geq0$, $Y^n\leq Y^{n+1}\leq U$, $K_s^{+,n}\geq K_s^{+,n+1}$ and $K_s^{-,n}\leq K_s^{-,n+1}$. So for any $t\leq T$ let us set $Y_t$  the pointwise limit of the sequence $(Y^n_t)_{n\geq0}$. \\
\par We now show the following Lemmas (which are steps forward in the proof of Proposition $\ref{prop2.1}.$) related
to the estimation of the processes $Y^n,$ $n\geq0$ and the convergence of the sequence $(Y^n)_{n\geq0}$.
\begin{lemma}\label{lemma22}
There exists a $\mathcal{P}$-measurable RCLL process $Y:=(Y_t )_{t\leq T}$, $\mathbb{R}$-valued and such that $\P$-a.s. for any $t\leq T,$ $Y^n_t\nearrow Y_t.$ Moreover for any constant $\gamma\geq 1$ and any stopping times $\tau\in[0,T],$
\begin{equation}
\mathbb{E}[|Y_{\tau}|^{\gamma}]\leq C,
\end{equation}
where $C$ is a constant independent of $\tau.$
\end{lemma}
\begin{proof}
Let $\P^n$ be the probability, equivalent to $\P$, defined as follows:
 $$ \mathrm{d}{\P^n}=L_T^n\,\mathrm{d}{\P},$$
where for any $t\leq T,$
\begin{equation*}
L_t^n:=\exp\left\{\int_{0}^{t}[C(1+||x||_s)\wedge n]\frac{|Z_s^n|}{Z_s^n}1_{\{Z_s^n\neq0\}}\, \mathrm{d}{B_s}-\frac{1}{2}\int_{0}^{t}||[C(1+||x||_s)\wedge n]\frac{|Z_s^n|}{Z_s^n}1_{\{Z_s^n\neq0\}}||^2 \, \mathrm{d}{s} \right\}.
\end{equation*}
By Girsanov's Theorem, the process $(B_t^n:=B_t-\int_{0}^{t}[C(1+||x||_s)\wedge n]\frac{|Z_s^n|}{Z_s^n}1_{\{Z_s^n\neq0\}}\, \mathrm{d}{s})_{t\leq T}$ is a Brownian motion under $\P^n$ and the quadruple $(Y^n,Z^n,K^{+,n},K^{-,n})$ verifies: For any $t\in[0,T],$
\begin{equation*}
\begin{cases}
Y^n_t = g(x) + \int_{t}^{T} C(1+||x||^{p}_{s}) \, \mathrm{d}{s} + (K_T^{+,n} - K_t^{+,n}) - ( K_T^{-,n} - K_t^{-,n}) -  \int_{t}^{T}Z^n_s \, \mathrm{d}{B_s^n};\\
\forall t \leq T, \quad L_t \leq Y^n_t \leq U_t,\\
\int_{0}^{T} (Y^n_s - L_s) \, \mathrm{d}{K_s^{+,n}} = 0 \; and \; \int_{0}^{T} (U_s-Y^n_s) \, \mathrm{d}{K_s^{-,n}} = 0.
\end{cases}
\end{equation*}
Therefore for any stopping time $\tau\leq T,$ $\P$-a.s., we have:
\begin{equation}
\begin{aligned}
Y^n_{\tau}=\, \essinf_{\nu\in\mathcal{ T}_{\tau}}\underset{\sigma\in\mathcal{ T}_{\tau}}{\esssup}\, \mathbb{E}^{\P^n} \biggl[&\int^{\nu\wedge\sigma}_{\tau} C(1+||x||^{p}_{s}) \, \mathrm{d}{s} + L_{\sigma}1_{[\sigma\leq\nu<T]}\biggr.\\
&\; \biggl.\, +U_{\nu}1_{[\nu<T]}+ g(x)1_{[\sigma=\nu =T]}\mid\mathcal{F}_{\tau} \biggr].
\end{aligned}
\end{equation}
(one can see $\cite{CK}$ for this characterization). Next by using assumptions $\mathbf{[H3]}$, $\mathbf{[H4]}$ and the fact that $\P$ and $\P^n$ are equivalent, we deduce that:
\begin{equation*}
\begin{aligned}
|Y^n_{\tau}|& \, \leq \, \mathbb{E}^{\P^n} \left[C\int^{T}_{\tau} (1+||x||^{p}_{s}) \, \mathrm{d}{s} + 3C(1+||x||^{p}_{T})\mid\mathcal{F}_{\tau}\right],\\
               &\, \leq C(1+\mathbb{E}^{\P^n} [||x||^{p}_{T}\mid\mathcal{F}_{\tau}]).
\end{aligned}
\end{equation*}
Next let $\gamma\geq1$, then by conditional Jensen's inequality we have:
\begin{equation}\label{estmY}
\begin{aligned}
|Y^n_{\tau}|^{\gamma}& \, \leq \, C(1+\mathbb{E}^{\P^n} [||x||^{\gamma p}_{T}\mid\mathcal{F}_{\tau}]),\\
               & \, = C(1+\mathbb{E} [L^n_{\tau,T}||x||^{\gamma p}_{T}\mid\mathcal{F}_{\tau}]).
\end{aligned}
\end{equation}
where for any $t\leq T,$
\begin{equation*}
\begin{aligned}
L^{n}_{t,T}\, & := \, \frac{L^n_{T}}{L^{n}_{t}},\\
             & \, =\exp\left\{\int_{t}^{T}[C(1+||x||_s)\wedge n]\frac{|Z_s^n|}{Z_s^n}1_{\{Z_s^n\neq0\}}\, \mathrm{d}{B_s}\right.\\
             & \left.\qquad\, -\frac{1}{2}\int_{t}^{T}||[C(1+||x||_s)\wedge n]\frac{|Z_s^n|}{Z_s^n}1_{\{Z_s^n\neq0\}}||^2 \, \mathrm{d}{s} \right\},\\
             & \, =\exp\left\{\int_{0}^{T}1_{[t\leq s\leq T]}\biggl([C(1+||x||_s)\wedge n]\frac{|Z_s^n|}{Z_s^n}1_{\{Z_s^n\neq0\}}\, \mathrm{d}{B_s} \right.\\
             & \left. \qquad\, -\frac{1}{2}||[C(1+||x||_s)\wedge n]\frac{|Z_s^n|}{Z_s^n}1_{\{Z_s^n\neq0\}}||^2 \, \mathrm{d}{s} \biggr)\right\}.
\end{aligned}
\end{equation*}
Now by Lemma $\ref{lemma21}$, there exist constants $C$ and $q> 1,$ which do not depend on $n$ and $\tau$ such that $\mathbb{E}[(L^n_{\tau,T})^q]\leq C.$
 Next by Young's inequality we get from $(\ref{estmY}):$
 \begin{equation*}
 |Y^n_{\tau}|^{\gamma} \, \leq \,  C\left(1+\mathbb{E} \left[\frac{1}{q} (L^n_{\tau,T})^q+\frac{1}{q'}||x||^{\gamma pq'}_{T}\mid\mathcal{F}_{\tau}\right]\right).
 \end{equation*}
where $\frac{1}{q}+\frac{1}{q'}=1.$  Next taking into account of $(\ref{estmx})$ we deduce that:
\begin{equation}\label{estmYn}
\mathbb{E}[|Y^n_{\tau}|^{\gamma}]\leq C, \qquad\qquad \forall n\geq0.
\end{equation}
where $C$ is a constant which do not depend on $n$. As for any $n\geq0$, $Y^n\leq Y^{n+1}\leq U$, since $Y^n$ is a continuous supermartingale which converges increasingly and pointwise to $Y$, then the sequence of processes $(Y^n)_{n\geq0}$ converges to a RCLL process $Y:=(Y_t)_{t\leq T}$ (see \cite{DMM}, pp. $86$). which satisfies $L_t \leq Y_t \leq U_t.$ Therefore $\P$-a.s. $Y^n_{\tau}\longrightarrow Y_{\tau},$ as $n \longrightarrow +\infty.$ And by $(\ref{estmYn})$ and  Fatou's Lemma we have
 \begin{equation*}
\mathbb{E}[|Y_{\tau}|^{\gamma}]\leq C.
\end{equation*}
\end{proof}
\par Next let $\zeta = |x_0|+|L_0|+|U_0|+|Y^0_0|+|Y_0|$ ($\zeta$ is a constant) and for $k\geq1$ let us define the sequence of stopping times $(\tau_k )_{k\geq1}$ by:
 \begin{equation*}
 \tau_k := \inf\{t\geq 0, ||x||_t+|L_t|+|U_t|+|Y^0_t|+|Y_t|\geq \zeta+k  \}\wedge T.
 \end{equation*}
 The sequence of stopping times  $(\tau_k )_{k\geq1}$ is non-decreasing, of stationary type (i.e. constant after some rank $k_0(\o)$) converging to $T$ since the process $Y$ is RCLL and $Y^0$, $x$, $U$ and $L$ are continuous. Moreover for any $k\geq1,$
\begin{equation*}
\max\{\sup_{t\leq\tau_k}|L_t|, \sup_{t\leq\tau_k}|U_t|, \sup_{t\leq\tau_k}|Y_t|, \sup_{t\leq\tau_k}|Y^n_t|\}\leq \zeta+k := \zeta_k.
\end{equation*}
Next we have the following result:
\begin{lemma}\label{lemma23}${}$\\
i) The process $Y$ is continuous.\\
ii) There exist $\mathcal{P}$-measurable processes $K^+$, $K^-$ and $Z$ valued in $\mathbb{R}^{1+1+d}$ such that $(Y,Z,K^+,K^-)$ is a solution of the double barrier reflected BSDE associated with $(g, \varphi,L, U)$ and verifies \eqref{cond3}.
\end{lemma}
\ni \textbf{Proof}:
For any $k\geq1$ and $n\geq0$ we have,  for any $t\leq T,$
\begin{equation}\label{eqn}
\begin{cases}
\begin{aligned}
Y^n_{t\wedge\tau_k} = Y^n_{\tau_k} + \int_{t\wedge\tau_k}^{\tau_k} &\varphi^n(s,x,Z^n_s) \, \mathrm{d}{s} +( K_{\tau_k}^{+,n} - K_{t\wedge\tau_k}^{+,n}) - ( K_{\tau_k}^{-,n} - K_{t\wedge\tau_k}^{-,n})\\ & -  \int_{t\wedge\tau_k}^{\tau_k}Z^n_s \, \mathrm{d}{B_s};
\end{aligned}\\
\forall t \leq T, \quad L_{t\wedge\tau_k} \leq Y^n_{t\wedge\tau_k} \leq U_{t\wedge\tau_k},\\
\int_{0}^{\tau_k} (Y^n_s - L_s) \, \mathrm{d}{K_s^{+,n}} = 0 \; and \; \int_{0}^{\tau_k} (U_s-Y^n_s) \, \mathrm{d}{K_s^{-,n}} = 0.
\end{cases}
\end{equation}
By using Itô's formula with $(Y^n_{t\wedge\tau_k})^2$ and taking into account of \eqref{estmx}, we classicaly deduce the existence of a constant $C_k$, which depends on $k$, such that uniformly on $n,$ we have
\begin{equation}\label{zntau}
\mathbb{E}\left[\int^{\tau_k}_{0}|Z^n_s|^2ds\right]<C_k,
\end{equation}
since $|Y^n_{t\wedge\tau_k}|\le \zeta_k$, for any $t\le T$, and $Y^0\le Y^n\le Y$.\\
Under $\eqref{eqn}$, we have,
\begin{equation}\label{eqkn-+}
K_{\tau_k}^{-,n} \leq K_{\tau_k}^{+,n}+|Y_{0}^n| + |Y_{\tau_k}^n| + \int_{0}^{\tau_k} (1+||x||_s)|Z^n_s| \, \mathrm{d}{s}+  \biggr|\int_{0}^{\tau_k}Z^n_s \, \mathrm{d}{B_s}\biggl|.
\end{equation}
The sequence of increasing processes $(K^{+,n})_{n\geq1}$ is non-increasing then it is convergent to a process $(K_{t\wedge\tau_k}^{k+})_{t\leq T}$ which moreover is increasing, upper semi-continuous and integrable since $\E[(K_{\tau_k}^{k+})^{\gamma}]\leq \E[(K_{\tau_k}^{+,0})^{\gamma}]<+\infty$.\\
Now, inequality $\eqref{eqkn-+}$ implies that for any $n\geq1,$ $\E[(K_{\tau_k}^{-,n})^{\gamma}]<C'_k$ for some constant $C'_k$, which depends on $k$.  On the
other hand the sequence of increasing processes $(K^{-,n})_{n\geq1}$ is increasing then, in
combination with Fatou's Lemma, it converges also to a process $(K_{t\wedge\tau_k}^{k-})_{t\leq T}$ which moreover is lower semi-continuous and satisfies $\E[(K_{\tau_k}^{k-})^{\gamma}]<+\infty$.
\\
\par Next once more by Itô's formula with $(Y^n_{t\wedge\tau_k}-Y^m_{t\wedge\tau_k})^2$ yields, for any $t\leq T,$
\begin{equation*}
\begin{aligned}
(Y^n_{t\wedge\tau_k}-Y^m_{t\wedge\tau_k})^2\,  = &\, (Y^n_{\tau_k}-Y^m_{\tau_k})^2 + 2\int_{t\wedge\tau_k}^{\tau_k}(Y^n_{s}-Y^m_{s})( \varphi^n(s,x,Z^n_s)-\varphi^m(s,x,Z^m_s))\mathrm{d}{s}  \\
             & \, + 2\int_{t\wedge\tau_k}^{\tau_k}(Y^n_{s}-Y^m_{s})d(K_{s}^{+,n} - K_{s}^{+,m} -  K_{s}^{-,n} + K_{s}^{-,m}) \\
             & \, -2\int_{t\wedge\tau_k}^{\tau_k}(Y^n_{s}-Y^m_{s})(Z^n_s-Z^m_s) \, \mathrm{d}{B_s}
             - \int_{t\wedge\tau_k}^{\tau_k}|Z^n_s-Z^m_s|^2 \, \mathrm{d}{s} .
\end{aligned}
\end{equation*}
The definition of $\varphi^n$ and the fact that,
\begin{equation*}
\begin{aligned}
 \int_{t\wedge\tau_k}^{\tau_k}&(Y^n_{s}-Y^m_{s})d(K_{s}^{+,n} - K_{s}^{+,m} -  K_{s}^{-,n} + K_{s}^{-,m})\\
& =\int_{t\wedge\tau_k}^{\tau_k}(Y^n_{s}-Y^m_{s})d(K_{s}^{+,n}-  K_{s}^{+,m})
 - \int_{t\wedge\tau_k}^{\tau_k}(Y^n_{s}-Y^m_{s})d(K_{s}^{-,n}- K_{s}^{-,m})  \\
 & = \int_{t\wedge\tau_k}^{\tau_k}(L_s-Y^m_{s})dK_{s}^{+,n}-\int_{t\wedge\tau_k}^{\tau_k}(Y^n_{s}-L_s)dK_{s}^{+,m}- \int_{t\wedge\tau_k}^{\tau_k}(U_s-Y^m_{s})dK_{s}^{-,n}\\
 & \qquad +\int_{t\wedge\tau_k}^{\tau_k}(Y^n_{s}-U_s)dK_{s}^{-,m} \leq 0.
\end{aligned}
\end{equation*}
Imply that, for any $t\leq T,$
\begin{equation}\label{itos2}
\begin{aligned}
(Y^n_{t\wedge\tau_k} & -Y^m_{t\wedge\tau_k})^2+ \int_{t\wedge\tau_k}^{\tau_k}|Z^n_s-Z^m_s|^2 \, \mathrm{d}{s} \\
    & \leq \, (Y^n_{\tau_k}-Y^m_{\tau_k})^2 + 2C\int_{t\wedge\tau_k}^{\tau_k}|Y^n_{s}-Y^m_{s}|(1+||x||_s)(|Z^n_s|+|Z^m_s|)\mathrm{d}{s} \\
             & \, -2\int_{t\wedge\tau_k}^{\tau_k}(Y^n_{s}-Y^m_{s})(Z^n_s-Z^m_s) \, \mathrm{d}{B_s}.
\end{aligned}
\end{equation}
But the process $(\int_{0}^{t\wedge\tau_k}(Y^n_{s}-Y^m_{s})(Z^n_s-Z^m_s) \, \mathrm{d}{B_s})_{t\leq T}$ is a martingale, then take expectation on both hand-sides to deduce that,
\begin{equation*}
\begin{aligned}
\mathbb{E}\left[\int_{t\wedge\tau_k}^{\tau_k}|Z^n_s-Z^m_s|^2 \, \mathrm{d}{s}\right] & \leq \mathbb{E}\biggl[(Y^n_{\tau_k}-Y^m_{\tau_k})^2\biggr]\\ & \qquad +2C\mathbb{E}\left[\int_{t\wedge\tau_k}^{\tau_k}|Y^n_{s}-Y^m_{s}|(1+||x||_s)(|Z^n_s|+|Z^m_s|)\mathrm{d}{s} \right].
\end{aligned}
\end{equation*}
Now the definition of $\tau_k$, estimate $(\ref{zntau})$ and the Cauchy-Schwarz's inequality yield:
\begin{equation}\label{Znm}
\mathbb{E}\left[\int_{t\wedge\tau_k}^{\tau_k}|Z^n_s-Z^m_s|^2 \, \mathrm{d}{s}\right]\longrightarrow 0,\,\,  as \,\, n,\,\, m \longrightarrow +\infty.
\end{equation}
Consequently the sequence $((Z^n_t1_{\{t\leq \tau_k\}})_{t\leq T})_{n\geq0}$ converges in $\mathcal{H}^{2,d}$ to a process which we denote $(Z_t^k)_{t\leq T} .$
\par Now going back to $(\ref{itos2}),$ take the supremum over $t$, make use of  Burkholder-Davis-Gundy (see e.g. $\cite{KS}$, $\cite{RY}$ and BDG for short) inequality  and finally take the expectation to deduce that
\begin{equation}
\mathbb{E}[\sup_{s\leq T}|Y^n_{s\wedge\tau_k}-Y^m_{s\wedge\tau_k}|^2 ]\longrightarrow 0, \,\, as \,\, n,\,\, m \longrightarrow +\infty.
\end{equation}
It follows that the process $(Y_{t\wedge\tau_k})_{t\leq T}$ is continuous for any $k\geq1.$ As $(\tau_k)_{k\geq1}$ is a stationary sequence then the process $(Y_{t})_{t\leq T}$ is continuous.\\
Next we focus on the continuity of the processes $K^{\pm}.$ According to $(\ref{eqn}),$ for any $t\leq T$ we have,
\begin{equation}\label{eqkn}
K_{t\wedge\tau_k}^{+,n} - K_{t\wedge\tau_k}^{-,n} = Y_{0}^n - Y_{t\wedge\tau_k}^n - \int_{0}^{t\wedge\tau_k} \varphi^n(s,x,Z^n_s) \, \mathrm{d}{s}+  \int_{0}^{t\wedge\tau_k}Z^n_s \, \mathrm{d}{B_s}.
\end{equation}
But by $(\ref{zntau})$ and $(\ref{Znm})$, for any $k \geq0$, the sequence of process $((\varphi^n(t,x,Z^n_t)1_{\{t\leq \tau_k\}})_{t\leq T})_{n\geq0}$ converges in $\mathcal{H}^{2,d}$ to $(\varphi(t,x,Z^k_t)1_{\{t\leq \tau_k\}})_{t\leq T}$.\\
Therefore from $(\ref{eqkn})$ there exists a subsequence of $(K^{+,n} - K^{-,n})_{n\geq0}$ such that:
\begin{equation}\label{convkn}
\mathbb{E}[\sup_{s\leq T}|(K^{+,n}_{s\wedge\tau_k}-K^{-,n}_{s\wedge\tau_k})-(K^{+,m}_{s\wedge\tau_k}-K^{-,m}_{s\wedge\tau_k})|^2 ]\longrightarrow 0, \,\, as \,\, n,\,\, m \longrightarrow +\infty.
\end{equation}
Indeed, for any $t\leq T$, we have
\begin{equation*}
\begin{aligned}
\mathbb{E}\left[\sup_{s\leq T}|(K^{+,n}_{s\wedge\tau_k}-K^{-,n}_{s\wedge\tau_k})-(K^{+,m}_{s\wedge\tau_k}-K^{-,m}_{s\wedge\tau_k})|^2\right] & \leq \mathbb{E}\biggl[(Y^n_{0}-Y^m_{0})^2\biggr]+\mathbb{E}\left[\sup_{s\leq T}|Y^n_{s\wedge\tau_k}-Y^m_{s\wedge\tau_k}|^2\right]\\ & +C\mathbb{E}\left[\int_{0}^{t\wedge\tau_k}|\varphi^n(s,x,Z^n_s)-\varphi^m(s,x,Z^m_s)|^2 \, \mathrm{d}{s}\right].\\ & +\mathbb{E}\left[\int_{0}^{t\wedge\tau_k}|Z^n_s-Z^m_s|^2 \, \mathrm{d}{s}\right].
\end{aligned}
\end{equation*}
Consequently the sequence $(((K^{+,n} - K^{-,n})1_{\{t\leq \tau_k\}})_{t\leq T})_{n\geq0}$ converges to a process which we denote $(K_t^{k+} - K_t^{k-})_{t\leq T} .$ In addition, by $(\ref{convkn})$, $K^{k+} - K^{k-}$ is continuous. Therefore for any $k \geq0$ we have
\begin{equation}\label{eqzdovi}
\begin{aligned}
Y_{t\wedge\tau_k} = & Y_{\tau_k} + \int_{t\wedge\tau_k}^{\tau_k} \varphi(s,x,Z^k_s) \, \mathrm{d}{s} + (K_{\tau_k}^{k+} - K_{t\wedge\tau_k}^{k+}) \\
                &  - ( K_{\tau_k}^{k-} - K_{t\wedge\tau_k}^{k-}) -  \int_{t\wedge\tau_k}^{\tau_k}Z^k_s \, \mathrm{d}{B_s}, \qquad \forall t\leq T.
\end{aligned}
\end{equation}
Implies also that,
\begin{equation*}
K_{t\wedge\tau_k}^{k+}-K_{t\wedge\tau_k}^{k-} =  Y_{0}-Y_{t\wedge\tau_k} - \int_{0}^{t\wedge\tau_k} \varphi(s,x,Z^k_s) \, \mathrm{d}{s}+ \int_{0}^{t\wedge\tau_k}Z^k_s \, \mathrm{d}{B_s}, \qquad \forall t\leq T.
\end{equation*}
and then
\begin{equation*}
K_{t\wedge\tau_k}^{k+}=K_{t\wedge\tau_k}^{k-}+Y_{0}-Y_{t\wedge\tau_k} - \int_{0}^{t\wedge\tau_k} \varphi(s,x,Z^k_s) \, \mathrm{d}{s}+ \int_{0}^{t\wedge\tau_k}Z^k_s \, \mathrm{d}{B_s}, \qquad \forall t\leq T.
\end{equation*}
As $K^{k+}$ is upper semi-continuous and $K^{k-}$ is lower semi-continuous. It means that $K^{k+}$ and $K^{k-}$ are lower and upper semi-continuous in the same time then the processes $(K_{t\wedge\tau_k}^{k+})_{t\leq T}$ and $(K_{t\wedge\tau_k}^{k-})_{t\leq T}$ are continuous. Henceforth $K^{k+}$ and $K^{k-}$ are continuous since $(\tau_k)_{k\geq1}$ is of stationary type. In addition, from Dini's theorem, the sequences $(K^{n,+})_{n\geq0}$ and $(K^{n,-})_{n\geq0}$ converge $\P$-a.s. uniformly to $K^{k+}$ and $K^{k-}$ respectively.
\par Finally, the uniform convergence of $Y^n$ (resp. $K^{+,n}$; resp. $K^{-,n}$) to $Y$ (resp. $K^{k+}$; resp. $K^{k-}$) and the facts that $\int_{0}^{T} (Y^n_s - L_s) \, \mathrm{d}{K^{+,n}} = \int_{0}^{T} (U_s-Y^n_s) \, \mathrm{d}{K^{-,n}} = 0$,  imply that, in view of Helly's Convergence Theorem (see e.g. $\cite{KF}$, pp, $370$),
$$
\int_{0}^{\tau_k} (Y_t - L_t) \, \mathrm{d}{K_t^{k+}} = 0 \; and \; \int_{0}^{\tau_k} (U_t-Y_t) \, \mathrm{d}{K_t^{k-}} = 0.
$$
It means that for any $k\geq1$ we have: $\forall t\leq T$,

\begin{equation*}
\begin{cases}
\begin{aligned}
Y_{t\wedge\tau_k} = Y_{\tau_k} + \int_{t\wedge\tau_k}^{\tau_k} &\varphi(s,x,Z^k_s) \, \mathrm{d}{s} +( K_{\tau_k}^{k+} - K_{t\wedge\tau_k}^{k+}) - ( K_{\tau_k}^{k-} - K_{t\wedge\tau_k}^{k-})\\ & -  \int_{t\wedge\tau_k}^{\tau_k}Z^k_s \, \mathrm{d}{B_s};
\end{aligned}\\
\forall t \leq T, \quad L_{t\wedge\tau_k} \leq Y_{t\wedge\tau_k} \leq U_{t\wedge\tau_k},\\
\int_{0}^{\tau_k} (Y_s - L_s) \, \mathrm{d}{K_s^{k+}} = 0 \; and \; \int_{0}^{\tau_k} (U_s-Y_s) \, \mathrm{d}{K_s^{k-}} = 0.
\end{cases}
\end{equation*}
Take now the double barrier reflected BSDE satisfied by $(Y,Z^{k+1},K^{k+1})$ on $[0,\tau_k]$ (since $\tau_k\leq \tau_{k+1}$) yields:  $\forall t\leq T$,
\begin{equation*}
\begin{cases}
\begin{aligned}
Y_{t\wedge\tau_k} = Y_{\tau_k} + \int_{t\wedge\tau_k}^{\tau_k} &\varphi(s,x,Z^{k+1}_s) \, \mathrm{d}{s} +( K_{\tau_k}^{(k+1)+} - K_{t\wedge\tau_k}^{(k+1)+}) - ( K_{\tau_k}^{(k+1)-} - K_{t\wedge\tau_k}^{(k+1)-})\\ & -  \int_{t\wedge\tau_k}^{\tau_k}Z^{k+1}_s \, \mathrm{d}{B_s};
\end{aligned}\\
\forall t \leq T, \quad L_{t\wedge\tau_k} \leq Y_{t\wedge\tau_k} \leq U_{t\wedge\tau_k},\\
\int_{0}^{\tau_k} (Y_s - L_s) \, \mathrm{d}{K_s^{(k+1)+}} = 0 \; and \; \int_{0}^{\tau_k} (U_s-Y_s) \, \mathrm{d}{K_s^{(k+1)-}} = 0.
\end{cases}
\end{equation*}
By uniqueness  (which holds since on $[0,\tau_k]$, $\varphi$ is Lipschitz w.r.t. $z$), we have for any $k\geq1$ :
$$Z_{t}^{k+1}1_{\{t\leq\tau_k\}}=Z_{t}^{k}1_{\{t\leq\tau_k\}},dt\otimes d\P-a.s. \mbox{ and }K_{t\wedge\tau_k}^{(k+1)\pm}=K_{t\wedge\tau_k}^{k\pm},\,\,\forall t\leq T.$$
\par Finally let us define (by concatenation) the processes $Z$ and $K^{\pm}$ as follows: $\forall t\le T$,
\begin{equation}
Z_t= Z^1_t1_{\{t\leq \tau_k\}}+\sum_{k\geq2} Z^k_t1_{\{\tau_{k-1}<t\leq\tau_{k}\}}.
\end{equation}
and
\begin{equation*}
K^{\pm}_t=
\begin{cases}
K_{t}^{1\pm} \mbox{ if }t\leq \t_1,\\
(K_{t}^{(k+1)\pm}-K_{\t_k}^{(k+1)\pm}) +K_{\t_k}^{\pm}\mbox{ for }\t_k<t\leq \t_{k+1}, \,\,k\geq 1.
\end{cases}
\end{equation*}
Note that the processes $Z$ and $K^{\pm}$ are well-defined since the sequence ($\t_k)_{k\ge 1}$ is of stationary type.
On the other hand, $K^{\pm}$ is continuous non-decreasing and $\P$-a.s., $K^{\pm}_T(\omega)<+\infty$ and $(Z_t(\omega))_{t\leq T}$ is $ds-$square integrable. Finally for any $k\geq1$, it holds: $\forall t\leq T$,
\begin{equation}\label{eqttau}
\begin{cases}
\begin{aligned}
Y_{t\wedge\tau_k} = Y_{\tau_k} + \int_{t\wedge\tau_k}^{\tau_k} &\varphi(s,x,Z_s) \, \mathrm{d}{s} +( K_{\tau_k}^{+} - K_{t\wedge\tau_k}^{+}) - ( K_{\tau_k}^{-} - K_{t\wedge\tau_k}^{-})\\ & -  \int_{t\wedge\tau_k}^{\tau_k}Z_s \, \mathrm{d}{B_s};
\end{aligned}\\
\forall t \leq T, \quad L_{t\wedge\tau_k} \leq Y_{t\wedge\tau_k} \leq U_{t\wedge\tau_k},\\
\int_{0}^{\tau_k} (Y_s - L_s) \, \mathrm{d}{K_s^{+}} = 0 \; and \; \int_{0}^{\tau_k} (U_s-Y_s) \, \mathrm{d}{K_s^{-}} = 0.
\end{cases}
\end{equation}
Take now $k$ great enough and since once more $(\tau_k)_{k\geq 1}$ is of stationary type to obtain: $\forall t\leq T $,
\begin{equation}\label{drbsde}
\begin{cases}
\begin{aligned}
Y_{t} = g(x) + \int_{t}^{T} &\varphi(s,x,Z_s) \, \mathrm{d}{s} +( K_{T}^{+} - K_{t}^{+}) - ( K_{T}^{-} - K_{t})\\ & - \int_{t}^{T}Z_s \, \mathrm{d}{B_s};
\end{aligned}\\
\forall t \leq T, \quad L_{t} \leq Y_{t} \leq U_{t},\\
\int_{0}^{T} (Y_s - L_s) \, \mathrm{d}{K_s^{+}} = 0 \; and \; \int_{0}^{T} (U_s-Y_s) \, \mathrm{d}{K_s^{-}} = 0.
\end{cases}
\end{equation}

The proof of Proposition $\ref{prop2.1}.$ is now complete.$\Box$
\section*{ Appendix B: Proof of Proposition $\ref{prop31}$}

\par
For any $m,n\geq1$, let us set,
\begin{equation}\label{Hnm}
H^{*n,m}(t,x,z)=H^{*+}(t,x,z)1_{\{||x||_t\leq n\}}-H^{*-}(t,x,z)1_{\{||x||_t\leq m\}}.
\end{equation}
Then $H^{*n,m}$ is Lipschitz w.r.t. z, therefore by Hamadène and Hassani result $\cite{HM}$, there exists a quadruple of processes $(Y^{*n,m},Z^{*n,m},K^{*+,n,m},K^{*-,n,m})\in\mathcal{S}^{2}\times\mathcal{H}^{2,d}\times\mathcal{S}_{ci}^2\times\mathcal{S}_{ci}^2$ that  satisfies:
\begin{equation}\label{RBSDEnm}
\begin{cases}
\begin{aligned}
Y^{*n,m}_t = \quad & g(x) + \int_{t}^{T} H^{*n,m}(s,x,Z^{*n,m}_s) \, \mathrm{d}{s} + (K_T^{*+,n,m} - K_t^{*+,n,m})\\ & - ( K_T^{*-,n,m} - K_t^{*-,n,m}) -  \int_{t}^{T}Z^{*n,m}_s \, \mathrm{d}{B_s};
\end{aligned}\\
\forall t \leq T, \quad L_t \leq Y^{*n,m}_t \leq U_t,\\
\int_{0}^{T} (Y^{*n,m}_s - L_s) \, \mathrm{d}{K_s^{*+,n,m}} = 0 \; and \; \int_{0}^{T} (U_s-Y^{*n,m}_s) \, \mathrm{d}{K_s^{*-,n,m}} = 0.
\end{cases}
\end{equation}
Thus from the definition of $H^{*n,m}$, we can clearly see that it is a non-decreasing (resp. non-increasing) sequence of functions w.r.t. $n$ (resp. $m$). Then according to comparison theorem (see, $\cite{HM}$), we have
\begin{equation}\label{compar}
\begin{aligned}
&Y^{*n,m}\leq Y^{*n+1,m},\; K^{*+,n,m}\geq K^{*+,n+1,m},\; K^{*-,n,m}\leq K^{*-,n+1,m},\\
&Y^{*n,m}\geq Y^{*n,m+1},\; K^{*+,n,m}\leq K^{*+,n,m+1}\;  and \; K^{*-,n,m}\geq K^{*-,n,m+1}.
\end{aligned}
\end{equation}
Now since for any $t\leq T,$ $L_t \leq Y^{*n,m}_t \leq U_t$ then there exists a $\mathcal{P}$-measurable process $Y^{*m}:=(Y^{*m}_t)_{t\leq T}$ such that $\P$-a.s. for any $t\leq T$ the sequence $(Y^{*n,m}_t)_{n\geq1}$ converges pointwisely to $Y^{*m}_t.$
 \par On the other hand, let $(Y,Z,K^+,K^-)$ be the solution of the double barrier reflected BSDE associated with $(g,\varphi,L,U).$ The function $H^{*n,m}$ satisfies: \begin{equation}\label{Hphi}
 |H^{*n,m}(t,x,z)| \leq \varphi(t,x,z), \qquad \forall (t,x,z)\in [0,T]\times\Omega\times\mathbb{R}^d.
 \end{equation}
 Under $\eqref{Hphi}$, once again the comparison theorem implies that, for any $t\leq T $
\begin{equation}\label{ynkp}
 Y^{*n,m}\leq Y,\quad K^{*+,n,m}\geq K^{+},\quad and \quad K^{*-,n,m}\leq K^{-}, \, \; \, \forall n,m\geq 1.
\end{equation}
Next we will divide the proof into two steps. In the first (resp. second) one we will prove
$i)$ (resp. $ii)$).
\\ \underline{\bf{Step 1}}: From Lemma $\ref{lemma22}$ we have that for any $\gamma\geq 1$ and any stopping times $\tau\in[0,T],$  $\mathbb{E}[|Y_{\tau}|^{\gamma}]\leq C,$ where $C$ is a constant independent of $\tau.$ Then from $(\ref{ynkp}),$ we have for any constant $\gamma\geq 1$ and any stopping times $\tau\in[0,T],$
\begin{equation}\label{estimynm}
\mathbb{E}[|Y^{*n,m}_{\tau}|^{\gamma}]\leq C,\, \, \, \forall n,m\geq1.
\end{equation}
where $C$ is a constant that does not depend neither on $n$ nor $m.$ Therefore from Fatou's Lemma and $\eqref{estimynm}$, we get that, for any $\gamma\geq 1$ and stopping times $\tau\in[0,T],$
\begin{equation}\label{estimym}
\forall m\geq1,\;\; \mathbb{E}[|Y^{*m}_{\tau}|^{\gamma}]\leq C, \, \, and \, \, then \, \, |Y^{*m}|<+\infty,\; \P-a.s..
\end{equation}
On the other hand from $\eqref{compar}$ it holds that: $\forall m\geq1$
\begin{equation}\label{ymym+1}
L\leq Y^{*(m+1)}\leq Y^{*m}\leq Y,\;\; K^{+}\leq K^{*+,m}\leq K^{*+,m+1}\;  and \; K^{*-,m+1}\leq K^{*-,m}\leq K^{-}.
\end{equation}
 \\ \underline{\bf{Step 2}}:
  Let $(\tau^{*}_k )_{k\geq1}$ be the sequence of stopping times defined as follow:
 \begin{equation*}
\forall k\geq1,\;\; \tau^{*}_k := \inf\{t\geq0, ||x||_t+|L_t|+|U_t|+|Y_t|\geq \zeta^{*}_{k};= \zeta^{*}+k  \}\wedge T,
 \end{equation*}
 where $\zeta^{*} = |x_0|+|L_0|+|U_0|+|Y_0|.$ \\
 The sequence of stopping times  $(\tau^{*}_k )_{k\geq1}$ is non-decreasing, of stationary type converging to $T$ since the process $Y$, $x$, $L$ and $U$ are continuous. Moreover for any $k\geq1,$
\begin{equation*}
\max\{\sup_{t\leq\tau^{*}_k}|L_t|, \sup_{t\leq\tau^{*}_k}|U_t|, \sup_{t\leq\tau^{*}_k}|Y_t|, \sup_{t\leq\tau^{*}_k}|Y^{*n,m}_t|\}\leq \zeta^{*}_k.
\end{equation*}
Let us now consider the following double barrier reflected BSDE. For any $n,m\geq 1$ and any $k\geq1$ we have: $\forall t\in[0,T]$
\begin{equation}\label{eqn*}
\begin{cases}
\begin{aligned}
Y^{*n,m}_{t\wedge\tau^{*}_k}= \quad & Y^{*n,m}_{\tau^{*}_k} + \int_{{t\wedge\tau^{*}_k}}^{\tau^{*}_k} H^{*n,m}(s,x,Z^{*n,m}_s) \, \mathrm{d}{s} + (K_{\tau^{*}_k}^{*+,n,m} - K_{t\wedge\tau^{*}_k}^{*+,n,m})  \\ &-( K_{\tau^{*}_k}^{*-,n,m} - K_{t\wedge\tau^{*}_k}^{*-,n,m}) -  \int_{{t\wedge\tau^{*}_k}}^{\tau^{*}_k}Z^{*n,m}_s \, \mathrm{d}{B_s};
\end{aligned}\\
\forall t \leq T, \quad L_{t\wedge\tau^{*}_k} \leq Y^{*n,m}_{t\wedge\tau^{*}_k} \leq U_{t\wedge\tau^{*}_k},\\
\int_{0}^{\tau^{*}_k} (Y^{*n,m}_s - L_s) \, \mathrm{d}{K_s^{*+,n,m}} = 0 \; and \; \int_{0}^{\tau^{*}_k} (U_s-Y^{*n,m}_s) \, \mathrm{d}{K_s^{*-,n,m}} = 0.
\end{cases}
\end{equation}
Under $(\ref{ynkp})$ and as, $\forall m,\gamma\geq1,$ $\mathbb{E}[(K_{\tau^{*}_k}^{*+,1,m})^{\gamma}+(K_{\tau^{*}_k}^{-})^{\gamma}]<\infty,$ then $\P$-a.s., for any $t\leq T$ the sequence $(K_t^{*+,n,m})_{n\geq1}$ (resp. $(K_t^{*-,n,m})_{n\geq1}$) converges to $K_t^{*+,m}$ (resp. $K_t^{*-,m})$).
In addition, the process  $K^{*+,m}:=(K_t^{*+,m})_{t\leq T}$ (resp. $K^{*-,m}:=(K_t^{*-,m})_{t\leq T}$) is non-decreasing upper and (resp. lower) semi-continuous and $\mathbb{E}[(K_{\tau^{*}_k}^{*+,m})^{\gamma}]<\infty$ (resp. $\mathbb{E}[(K_{\tau^{*}_k}^{*-,m})^{\gamma}]<\infty$).
\par Now we take into account $({\ref{estmx}})$ and the fact that $\forall t\in[0,T]$, $|Y^{*n,m}_{t\wedge \tau^{*}_k}| \leq \zeta^{*}_k$ and  $Y^{*1,m}\leq Y^{*n,m}\leq Y^{*m}$, then using Itô's formula with $(Y^{*n,m}_{t\wedge\tau^{*}_k})^2$ to conclude that there exists a constant $C_k$ such that:
\begin{equation}\label{zn*tau*}
\mathbb{E}\left[\int^{\tau^{*}_k}_{0}|Z^{*n,m}_s|^2ds\right]<C_k, \, \, \, \forall n,m\geq1.
\end{equation}
\par Let us show that the process $(Y_t^{*m})_{t\leq T}$ is continuous and for any $t\leq T$ the sequence $(Z_t^{*n,m})_{n\geq1}$ is convergent in $\mathcal{H}^{2,d}$.
Using Itô's formula with $(Y^{*n,m}_{t\wedge\tau^{*}_k}-Y^{*p,m}_{t\wedge\tau^{*}_k})^2$ yields, for any $t\leq T,$
\begin{equation*}
\begin{aligned}
(Y^{*n,m}_{t\wedge\tau^{*}_k}-Y^{*p,m}_{t\wedge\tau^{*}_k})^2\,  = &\, (Y^{*n,m}_{\tau^{*}_k}-Y^{*p,m}_{\tau^{*}_k})^2 + 2\int_{t\wedge\tau^{*}_k}^{\tau^{*}_k}(Y^{*n,m}_{s}-Y^{*p,m}_{s})( H^{*n,m}(s,x,Z^{*n,m}_s)\\ &-H^{*p,m}(s,x,Z^{*p,m}_s))\mathrm{d}{s}  +2\int_{t\wedge\tau^{*}_k}^{\tau^{*}_k}(Y^{*n,m}_{s}-Y^{*p,m}_{s})d(K_{s}^{*+,n,m} -K_{s}^{*+,p,m} \\ & - K_{s}^{*-,n,m} + K_{s}^{*-,p,m}) -2\int_{t\wedge\tau^{*}_k}^{\tau^{*}_k}(Y^{*n,m}_{s}-Y^{*p,m}_{s})(Z^{*n,m}_s-Z^{*p,m}_s) \, \mathrm{d}{B_s} \\ &-\int_{t\wedge\tau^{*}_k}^{\tau^{*}_k}|Z^{*n,m}_s-Z^{*p,m}_s|^2 \, \mathrm{d}{s} .
\end{aligned}
\end{equation*}
Obviously, for any $t\leq T$,
\begin{equation*}
 \int_{t\wedge\tau^{*}_k}^{\tau^{*}_k}(Y^{*n,m}_{s}-Y^{*p,m}_{s})d(K_{s}^{*+,n,m} -K_{s}^{*+,p,m} -  K_{s}^{*-,n,m} + K_{s}^{*-,p,m})\leq 0.
\end{equation*}
and we take into account $\eqref{Hphi},$
then for any $t\leq T,$
\begin{equation}\label{itos2*}
\begin{aligned}
(Y^{*n,m}_{t\wedge\tau^{*}_k} & -Y^{*p,m}_{t\wedge\tau^{*}_k})^2+ \int_{t\wedge\tau^{*}_k}^{\tau^{*}_k}|Z^{*n,m}_s-Z^{*p,m}_s|^2 \, \mathrm{d}{s} \leq \, (Y^{*n,m}_{\tau^{*}_k}-Y^{*p,m}_{\tau^{*}_k})^2\\
& + 2\int_{t\wedge\tau^{*}_k}^{\tau^{*}_k}|Y^{*n,m}_{s}-Y^{*p,m}_{s}|[C(1+||x||_s)(|Z^{*n,m}_s|+|Z^{*p,m}_s|)+2C(1+||x||_s^p)]\mathrm{d}{s} \\ & \, -2\int_{t\wedge\tau^{*}_k}^{\tau^{*}_k}(Y^{*n,m}_{s}-Y^{*p,m}_{s})(Z^{*n,m}_s-Z^{*p,m}_s) \, \mathrm{d}{B_s}.
\end{aligned}
\end{equation}
Take now expectation in $\eqref{itos2*}$ to deduce that:
\begin{equation*}
\begin{aligned}
\mathbb{E}&\left[\int_{t\wedge\tau^{*}_k}^{\tau^{*}_k}|Z^{*n,m}_s-Z^{*p,m}_s|^2 \, \mathrm{d}{s}\right]  \leq \mathbb{E}\biggl[(Y^{*n,m}_{\tau^{*}_k}-Y^{*p,m}_{\tau^{*}_k})^2\biggr]\\ &\qquad\qquad +2\mathbb{E}\left[\int_{t\wedge\tau^{*}_k}^{\tau^{*}_k}|Y^{*n,m}_{s}-Y^{*p,m}_{s}|[C(1+||x||_s)(|Z^{*n,m}_s|+|Z^{*p,m}_s|)+2C(1+||x||_s^p)]\mathrm{d}{s} \right].
\end{aligned}
\end{equation*}
Now the definition of $\tau^{*}_k$, estimate $(\ref{zn*tau*})$ and the Cauchy-Schwarz's inequality yield:
\begin{equation}\label{Znm*}
\mathbb{E}\left[\int_{t\wedge\tau^{*}_k}^{\tau^{*}_k}|Z^{*n,m}_s-Z^{*p,m}_s|^2 \, \mathrm{d}{s}\right]\longrightarrow 0, \,\, as \,\, n,\, p \longrightarrow +\infty.
\end{equation}
Consequently the sequence $((Z^{*n,m}_t1_{\{t\leq \tau^{*}_k\}})_{t\leq T})_{n\geq1}$ converges in $\mathcal{H}^{2,d}$ to a process which we denote $(Z_t^{*k,m})_{t\leq T} .$ In addition it satisfies, for any $k\geq1$, $Z_{t\wedge\tau^{*}_k}^{*k+1,m}=Z_{t\wedge\tau^{*}_{k}}^{*k,m},$ $dt\otimes d\P-a.s.$.
  Next by $(\ref{itos2*})$ and the use of BDG inequality we obtain
\begin{equation}
\mathbb{E}[\sup_{s\leq T}|Y^{*n,m}_{s\wedge\tau^{*}_k}-Y^{*p,m}_{s\wedge\tau^{*}_k}|^2 ]\longrightarrow 0, \quad as \,\, n,\, p \longrightarrow +\infty.
\end{equation}
It follows that the process $(Y^{*m}_{t\wedge\tau^{*}_k})_{t\leq T}$ is continuous for any $k\geq1.$ As $(\tau^{*}_k)_{k\geq1}$ is a stationary sequence then the process $(Y^{*m}_{t})_{t\leq T}$ is continuous.
\par Now by $(\ref{zn*tau*})$ and $(\ref{Znm*})$, for any $k,m\geq1$, the sequence of process\\ $((H^{*n,m}(t,x,Z^{*n,m}_t)1_{\{t\leq \tau^{*}_k\}})_{t\leq T})_{n\geq1}$ converges in $\mathcal{H}^2$ to $(H^{*m}(t,x,Z^{*k,m}_t)1_{\{t\leq\tau^{*}_k\}})_{t\leq T}$. Therefore for any $k\geq1$ we have,
\begin{equation}\label{eqk*}
\begin{aligned}
Y^{*m}_{t\wedge\tau^{*}_k} = & Y^{*m}_{\tau^{*}_k} + \int_{t\wedge\tau^{*}_k}^{\tau^{*}_k} H^{*m}(s,x,Z^{*k,m}_s) \, \mathrm{d}{s} + (K_{\tau^{*}_k}^{*+,m} - K_{t\wedge\tau^{*}_k}^{*+,m}) \\
                 &  - ( K_{\tau^{*}_k}^{*-,m} - K_{t\wedge\tau^{*}_k}^{*-,m}) -  \int_{t\wedge\tau^{*}_k}^{\tau^{*}_k}Z^{*k,m}_s \, \mathrm{d}{B_s}, \qquad \forall t\leq T.
\end{aligned}
\end{equation}
Now $(\ref{eqk*})$ implies also that,
\begin{equation*}
Y^{*m}_{t\wedge\tau^{*}_k} = Y^{*m}_{0} -\int_{0}^{t\wedge\tau^{*}_k} H^{*m}(s,x,Z^{*k}_s) \, \mathrm{d}{s}- K_{t\wedge\tau^{*}_k}^{*+,m}+K_{t\wedge\tau^{*}_k}^{*-,m}+  \int_{0}^{t\wedge\tau^{*}_k}Z^{*k,m}_s \, \mathrm{d}{B_s}, \qquad \forall t\leq T.
\end{equation*}
And then,
\begin{equation*}
K_{t\wedge\tau^{*}_k}^{*+,m}= K_{t\wedge\tau^{*}_k}^{*-,m}+Y^{*m}_{0} -\int_{0}^{t\wedge\tau^{*}_k} H^{*m}(s,x,Z^{*k,m}_s) \, \mathrm{d}{s}+  \int_{0}^{t\wedge\tau^{*}_k}Z^{*k,m}_s \, \mathrm{d}{B_s}, \qquad \forall t\leq T.
\end{equation*}
As $K^{*+,m}$ is upper semi-continuous and $K^{*-,m}$ is lower semi-continuous. It means that $K^{*+,m}$ and $K^{*-,m}$ are lower and upper semi-continuous in the same time then the processes $(K_{t\wedge\tau^{*}_k}^{*+,m})_{t\leq T}$ and $(K_{t\wedge\tau^{*}_k}^{*-,m})_{t\leq T}$ are continuous. Henceforth $K^{*+,m}$ and $K^{*-,m}$ are continuous since $(\tau^{*}_k)_{k\geq1}$ is of stationary type. In addition, from Dini's theorem, the sequences $(K^{*n,+})_{n\geq1}$ and $(K^{*n,-})_{n\geq1}$ converge $\P$-a.s. uniformly to $K^{*+,m}$ and $K^{*-,m}$ respectively.
\par Now let us define $Z^{*m}$ by concatenation as follows: $\forall t\leq T$ and $m\geq1$,
\begin{equation}
Z^{*m}_t= Z^{*1,m}_t1_{\{t\leq\tau^{*m}_1\}}+\sum_{k\geq 2} Z^{*k,m}_t1_{\{\tau^{*m}_{k-1}<t\leq\tau^{*m}_{k}\}}.
\end{equation}
As $\mathbb{E}\left[\int^{\tau^{*}_k}_{0}|Z^{*k,m}_s|^2ds\right]<+\infty,$ for any $k\geq1$ and the sequence $(\tau^{*}_k)_{k\geq1}$ is of stationary type then $\int^{T}_{0}|Z^{*m}_s|^2ds<+\infty,$ $\P$-a.s.. On the other hand, with the definition of $Z^{*m}$ and $(\ref{eqk*})$ we have,
\begin{equation}\label{eqtau*}
\begin{aligned}
Y^{*m}_{t\wedge\tau^{*}_k} = & Y^{*m}_{\tau^{*}_k} + \int_{t\wedge\tau^{*}_k}^{\tau^{*}_k} H^{*m}(s,x,Z^{*m}_s) \, \mathrm{d}{s} + (K_{\tau^{*}_k}^{*+,m} - K_{t\wedge\tau^{*}_k}^{*+,m}) \\
                 &  - ( K_{\tau^{*}_k}^{*-,m} - K_{t\wedge\tau^{*}_k}^{*-,m}) -  \int_{t\wedge\tau^{*}_k}^{\tau^{*}_k}Z^{*m}_s \, \mathrm{d}{B_s}, \qquad \forall t\leq T.
\end{aligned}
\end{equation}
Now taking $k$ great enough in $(\ref{eqtau*})$ yields, $\forall t\leq T $
\begin{equation*}
Y^{*m}_t = g(x) + \int_{t}^{T} H^{*m}(s,x,Z^{*m}_s) \, \mathrm{d}{s} + (K^{*+,m}_T- K^{*+,m}_t) - (K^{*-,m}_T- K^{*-,m}_t) -  \int_{t}^{T}Z^{*m}_s \, \mathrm{d}{B_s}.
\end{equation*}
\par Finally it remains to show that $\int_{0}^{T} (Y^{*m}_s - L_s) \, \mathrm{d}{K_s^{*+,m}} = \int_{0}^{T} (U_s-Y^{*m}_s) \, \mathrm{d}{K_s^{*-,m}} = 0.$ But this is a direct consequence of the $\P$-a.s. uniform convergence of $Y^{*n,m}$ (resp. $K^{*+,n,m}$; resp. $K^{*-,n,m}$) to $Y^{*m}$ (resp. $K^{*+,m}$; resp. $K^{*-,m}$) and the facts that $\int_{0}^{T} (Y^{*n,m}_s - L_s) \, \mathrm{d}{K^{*+,n,m}} = \int_{0}^{T} (U_s-Y^{*n,m}_s) \, \mathrm{d}{K^{*-,n,m}} = 0.$
%(see e.g the Helly's Theorem in $\cite{KF}$,pp,$370$).\\
\\ The proof is now complete.$\Box$

\end{document}